%
\newif\ifloadreferences\loadreferencestrue
%
%
%
%
%
\let\myfrac=\frac%
\input eplain %
\let\frac=\myfrac%
\let\myfootnote=\footnote%
\input amstex \input epsf %
\let\footnote=\myfootnote%
%
%
\loadeufm\loadmsam\loadmsbm\message{symbol names}\UseAMSsymbols\message{,}%
\magnification 1200 %
\font\myfontdefault=cmr10%
\newif\ifmakebiblio%
\newif\ifinappendices%
\newif\ifundefinedreferences%
\newif\ifchangedreferences%
\makebibliofalse%
\undefinedreferencesfalse%
\changedreferencesfalse%
%
%
%
%
%
\def\setcatcodes{\catcode`\!=0 \catcode`\\=11}%
{\global\let\noe=\noexpand%
\catcode`\@=11 \catcode`\_=11 \setcatcodes%
!global!def!_@@internal@@makeref#1{%
!global!expandafter!def!csname #1ref!endcsname##1{%
!csname _@#1@##1!endcsname%
!expandafter!ifx!csname _@#1@##1!endcsname!relax%
    !write16{#1 ##1 not defined - run saving references}%
    !undefinedreferencestrue%
!fi}}%
!global!def!_@@internal@@makelabel#1{%
!global!expandafter!def!csname #1label!endcsname##1{%
!edef!temptoken{!csname #1info!endcsname}%
!ifloadreferences%
!expandafter!ifx!csname _@#1@##1!endcsname!relax%
!write16{#1 ##1 not hitherto defined - rerun saving references}%
!changedreferencestrue%
!else%
!expandafter!ifx!csname _@#1@##1!endcsname!temptoken%
!else%
!write16{#1 ##1 reference has changed - rerun saving references}%
!changedreferencestrue%
!fi%
!fi%
!else%
!expandafter!edef!csname _@#1@##1!endcsname{!temptoken}%
!edef!textoutput{!write!references{\global\def\_@#1@##1{!temptoken}}}%
!textoutput%
!fi}}%
!global!def!makecounter#1{!_@@internal@@makelabel{#1}!_@@internal@@makeref{#1}}%
!unsetcatcodes%
}
%
%
%
%
%
\def\turnintolatin#1{\ifcase #1 _\or i\or ii\or iii\or iv\or v\or vi\or vii\or viii\or ix\or x\or xi\or xii\or xiii\or xiv\or xv\or xvi\or xvii\or xviii\or xix\or xx\or xxi\or xxii\or xxiii\or xxiv\or xxv\or xxvi\fi}%
\def\alphanum#1{\ifcase #1 _\or A\or B\or C\or D\or E\or F\or G\or H\or I\or J\or K\or L\or M\or N\or O\or P\or Q\or R\or S\or T\or U\or V\or W\or X\or Y\or Z\fi}%
\newwrite\references%
\ifloadreferences{\catcode`\@=11 \catcode`\_=11 \global\def\_@citation@Aronszajn{1}
\global\def\_@citation@Berger{2}
\global\def\_@citation@CaffNirSprV{3}
\global\def\_@citation@ColdingDeLellis{4}
\global\def\_@citation@CrandhallIshiiLyons{5}
\global\def\_@citation@ElworthyTromba{6}
\global\def\_@citation@EspRos{7}
\global\def\_@citation@GilbTrud{8}
\global\def\_@citation@GuillemanPollack{9}
\global\def\_@citation@Hamilton{10}
\global\def\_@citation@Huisken{11}
\global\def\_@citation@Kato{12}
\global\def\_@citation@LabA{13}
\global\def\_@citation@LabB{14}
\global\def\_@citation@MaximoNunesSmith{15}
\global\def\_@citation@Milnor{16}
\global\def\_@citation@PacardXu{17}
\global\def\_@citation@Robeday{18}
\global\def\_@citation@RosSchneid{19}
\global\def\_@citation@SchneiderI{20}
\global\def\_@citation@SchneiderII{21}
\global\def\_@citation@SimonSmith{22}
\global\def\_@citation@SmiAAT{23}
\global\def\_@citation@SmiEC{24}
\global\def\_@citation@SmiPPG{25}
\global\def\_@citation@SmiSLC{26}
\global\def\_@citation@SmiAC{27}
\global\def\_@citation@Spivak{28}
\global\def\_@citation@Tromba{29}
\global\def\_@citation@WhiteI{30}
\global\def\_@citation@WhiteII{31}
\global\def\_@citation@WhiteIII{32}
\global\def\_@citation@WhiteIV{33}
\global\def\_@citation@Ye{34}
\global\def\_@head@Introduction{1}
\global\def\_@subhead@General{1.1}
\global\def\_@subhead@Background{1.2}
\global\def\_@subhead@Applications{1.3}
\global\def\_@proc@MainTheoremMeanCurvature{1.3.1}
\global\def\_@proc@MainTheoremExtrinsicCurvature{1.3.2}
\global\def\_@proc@MainTheoremSpecialLagrangianCurvature{1.3.3}
\global\def\_@proc@MainTheoremExtrinsicCurvatureTwoD{1.3.5}
\global\def\_@subhead@Acknowledgements{1.4}
\global\def\_@head@DegreeTheory{2}
\global\def\_@subhead@TheManifoldOfImmersions{2.1}
\global\def\_@eqn@DefinitionOfI{\relax \unhbox \voidb@x \hbox {{\relax \tenrm (2.1)}}}
\global\def\_@subhead@CurvatureAsAVectorField{2.2}
\global\def\_@eqn@DefinitionOfKi{\relax \unhbox \voidb@x \hbox {{\relax \tenrm (2.2)}}}
\global\def\_@eqn@DefinitionOfTheSectionK{\relax \unhbox \voidb@x \hbox {{\relax \tenrm (2.3)}}}
\global\def\_@eqn@CovariantDerivativeOfK{\relax \unhbox \voidb@x \hbox {{\relax \tenrm (2.4)}}}
\global\def\_@eqn@DefinitionOfEi{\relax \unhbox \voidb@x \hbox {{\relax \tenrm (2.5)}}}
\global\def\_@eqn@DefinitionOfTheSectionE{\relax \unhbox \voidb@x \hbox {{\relax \tenrm (2.6)}}}
\global\def\_@eqn@CovariantDerivativeOfE{\relax \unhbox \voidb@x \hbox {{\relax \tenrm (2.7)}}}
\global\def\_@eqn@DefinitionOfKHati{\relax \unhbox \voidb@x \hbox {{\relax \tenrm (2.8)}}}
\global\def\_@eqn@CovariantDerivativeOfKhati{\relax \unhbox \voidb@x \hbox {{\relax \tenrm (2.9)}}}
\global\def\_@eqn@JacobiOperatorOfKHat{\relax \unhbox \voidb@x \hbox {{\relax \tenrm (2.10)}}}
\global\def\_@eqn@DefinitionOfZ{\relax \unhbox \voidb@x \hbox {{\relax \tenrm (2.11)}}}
\global\def\_@eqn@SolutionSetOverF{\relax \unhbox \voidb@x \hbox {{\relax \tenrm (2.12)}}}
\global\def\_@subhead@Simplicity{2.3}
\global\def\_@proc@PropUniqueContinuation{2.3.1}
\global\def\_@proc@PropInjectivePointsAreDense{2.3.2}
\global\def\_@proc@CorEquivalenceOfNotionsOfSimplicity{2.3.3}
\global\def\_@subhead@Surjectivity{2.4}
\global\def\_@proc@NonDegeneracyI{2.4.1}
\global\def\_@proc@PropGeneralSurj{2.4.2}
\global\def\_@subhead@FiniteDimensionalSections{2.5}
\global\def\_@eqn@DefinitionOfKf{\relax \unhbox \voidb@x \hbox {{\relax \tenrm (2.13)}}}
\global\def\_@eqn@DefinitionOfZf{\relax \unhbox \voidb@x \hbox {{\relax \tenrm (2.14)}}}
\global\def\_@eqn@DefinitionOfRestrictionOfKf{\relax \unhbox \voidb@x \hbox {{\relax \tenrm (2.15)}}}
\global\def\_@eqn@DefinitionOfRestrictionOfZf{\relax \unhbox \voidb@x \hbox {{\relax \tenrm (2.16)}}}
\global\def\_@proc@SurjectivityInGraphChart{2.5.2}
\global\def\_@eqn@TangentSpaceOfZ{\relax \unhbox \voidb@x \hbox {{\relax \tenrm (2.17)}}}
\global\def\_@proc@NonDegeneracyYieldsSmoothness{2.5.5}
\global\def\_@subhead@Extensions{2.6}
\global\def\_@proc@Surjectivity{2.6.1}
\global\def\_@proc@Transversality{2.6.2}
\global\def\_@subhead@DefiningTheSignatureFD{2.7}
\global\def\_@proc@CanonicalLifting{2.7.1}
\global\def\_@eqn@FormulaForOrientationForm{\relax \unhbox \voidb@x \hbox {{\relax \tenrm (2.18)}}}
\global\def\_@eqn@FormulaForOrientation{\relax \unhbox \voidb@x \hbox {{\relax \tenrm (2.19)}}}
\global\def\_@proc@Extensions{2.7.4}
\global\def\_@proc@CorExtensions{2.7.5}
\global\def\_@subhead@DefiningTheSignatureID{2.8}
\global\def\_@eqn@RealEigenvaluesBoundedBelow{\relax \unhbox \voidb@x \hbox {{\relax \tenrm (2.20)}}}
\global\def\_@eqn@DefinitionOfSignature{\relax \unhbox \voidb@x \hbox {{\relax \tenrm (2.21)}}}
\global\def\_@eqn@DefinitionOfParity{\relax \unhbox \voidb@x \hbox {{\relax \tenrm (2.22)}}}
\global\def\_@proc@LemmaDefinitionOfParity{2.8.1}
\global\def\_@subhead@ConstructingTheDegree{2.9}
\global\def\_@proc@LemmaLowerBoundOfRealPart{2.9.1}
\global\def\_@eqn@DefinitionOfOrientationForm{\relax \unhbox \voidb@x \hbox {{\relax \tenrm (2.23)}}}
\global\def\_@proc@Orientation{2.9.2}
\global\def\_@proc@DefinitionOfTheDegree{2.9.3}
\global\def\_@subhead@VaryingTheMetric{2.10}
\global\def\_@eqn@DefinitionOfKCurvatureWithVaryingMetric{\relax \unhbox \voidb@x \hbox {{\relax \tenrm (2.24)}}}
\global\def\_@eqn@DefinitionOfCurvatureSectionWithVaryingMetric{\relax \unhbox \voidb@x \hbox {{\relax \tenrm (2.25)}}}
\global\def\_@eqn@DefinitionOfEvalutionSectionWithVaryingMetric{\relax \unhbox \voidb@x \hbox {{\relax \tenrm (2.26)}}}
\global\def\_@proc@PreliminaryToSurjectivityForMetrics{2.10.1}
\global\def\_@proc@PropGeneralSurjMetrics{2.10.2}
\global\def\_@head@Applications{3}
\global\def\_@subhead@GeneralisedSimonsFormula{3.1}
\global\def\_@proc@SimonsFormulaI{3.1.1}
\global\def\_@eqn@VanishingDerivatives{\relax \unhbox \voidb@x \hbox {{\relax \tenrm (3.1)}}}
\global\def\_@eqn@DefinitionOfFunctionA{\relax \unhbox \voidb@x \hbox {{\relax \tenrm (3.2)}}}
\global\def\_@proc@FirstOrderDifferentiabilityOfLambda{3.1.3}
\global\def\_@proc@SecondDerivativeOfA{3.1.4}
\global\def\_@proc@HopfTypeTheorem{3.1.5}
\global\def\_@eqn@BoundOnDeltaA{\relax \unhbox \voidb@x \hbox {{\relax \tenrm (3.3)}}}
\global\def\_@eqn@LaplacianOfFirstEigenvalue{\relax \unhbox \voidb@x \hbox {{\relax \tenrm (3.4)}}}
\global\def\_@fig@QualitativeBehaviour{3.1.1}
\global\def\_@subhead@SpheresOfPrescribedMeanCurvature{3.2}
\global\def\_@eqn@DefnMeanCurvature{\relax \unhbox \voidb@x \hbox {{\relax \tenrm (3.5)}}}
\global\def\_@eqn@ConditionsOnFPMC{\relax \unhbox \voidb@x \hbox {{\relax \tenrm (3.6)}}}
\global\def\_@eqn@DefinitionOfM{\relax \unhbox \voidb@x \hbox {{\relax \tenrm (3.7)}}}
\global\def\_@eqn@DefinitionOfTraceFreeCurvature{\relax \unhbox \voidb@x \hbox {{\relax \tenrm (3.8)}}}
\global\def\_@eqn@DefinitionOfNorm{\relax \unhbox \voidb@x \hbox {{\relax \tenrm (3.9)}}}
\global\def\_@eqn@KeyPolynomialPEC{\relax \unhbox \voidb@x \hbox {{\relax \tenrm (3.10)}}}
\global\def\_@eqn@SolutionSpacePMC{\relax \unhbox \voidb@x \hbox {{\relax \tenrm (3.11)}}}
\global\def\_@eqn@CodazziMainardi{\relax \unhbox \voidb@x \hbox {{\relax \tenrm (3.12)}}}
\global\def\_@eqn@Curvature{\relax \unhbox \voidb@x \hbox {{\relax \tenrm (3.13)}}}
\global\def\_@eqn@HessianOfRestriction{\relax \unhbox \voidb@x \hbox {{\relax \tenrm (3.14)}}}
\global\def\_@eqn@KeyInequalityCMC{\relax \unhbox \voidb@x \hbox {{\relax \tenrm (3.15)}}}
\global\def\_@proc@EstimateOfDeltaPinchingFactorPMC{3.2.2}
\global\def\_@proc@SolutionsArePinchedPMC{3.2.3}
\global\def\_@proc@ConnectedComponentsPMC{3.2.4}
\global\def\_@proc@PrimePMC{3.2.5}
\global\def\_@proc@UniformDiameterBoundsPMC{3.2.6}
\global\def\_@proc@ProperPMC{3.2.7}
\global\def\_@subhead@CalculatingTheDegree{3.3}
\global\def\_@eqn@AsymptoticFormulaForHI{\relax \unhbox \voidb@x \hbox {{\relax \tenrm (3.16)}}}
\global\def\_@eqn@AsymptoticFormulaForHII{\relax \unhbox \voidb@x \hbox {{\relax \tenrm (3.17)}}}
\global\def\_@proc@TheRadiiTendLinearlyToZero{3.3.1}
\global\def\_@eqn@MeanCurvatureOfNearbySpheres{\relax \unhbox \voidb@x \hbox {{\relax \tenrm (3.18)}}}
\global\def\_@proc@ExistencePMCI{3.3.2}
\global\def\_@proc@ExistencePMCII{3.3.3}
\global\def\_@subhead@ExtrinsicCurvature{3.4}
\global\def\_@eqn@QuarterPinched{\relax \unhbox \voidb@x \hbox {{\relax \tenrm (3.19)}}}
\global\def\_@eqn@PointwiseHalfPinched{\relax \unhbox \voidb@x \hbox {{\relax \tenrm (3.20)}}}
\global\def\_@eqn@ConditionsOnFPEC{\relax \unhbox \voidb@x \hbox {{\relax \tenrm (3.21)}}}
\global\def\_@eqn@SolutionSpacePEC{\relax \unhbox \voidb@x \hbox {{\relax \tenrm (3.22)}}}
\global\def\_@eqn@KLaplacian{\relax \unhbox \voidb@x \hbox {{\relax \tenrm (3.23)}}}
\global\def\_@proc@LowerBoundOnTrDKPEC{3.4.1}
\global\def\_@proc@DeltaAPEC{3.4.2}
\global\def\_@proc@EstimateOfDeltaPinchingFactorPEC{3.4.3}
\global\def\_@proc@UpperBoundOfShapeOperatorPEC{3.4.4}
\global\def\_@proc@PropernessPEC{3.4.5}
\global\def\_@proc@LowerBoundOfShapeOperatorPEC{3.4.6}
\global\def\_@proc@CountingPEC{3.4.7}
\global\def\_@subhead@SpecialLagrangianCurvature{3.5}
\global\def\_@proc@DiameterBoundPSLC{3.5.1}
\global\def\_@proc@ShapeOperatorBoundPSLC{3.5.2}
\global\def\_@proc@PropernessPSLC{3.5.3}
\global\def\_@proc@CountingPSLC{3.5.4}
\global\def\_@proc@CountingPSLCTWODIM{3.5.5}
\global\def\_@proc@CountingPSLCTHREEDIM{3.5.6}
\global\def\_@subhead@ExtrinsicCurvatureInTwoDimensions{3.6}
\global\def\_@eqn@ChristoffelSymbol{\relax \unhbox \voidb@x \hbox {{\relax \tenrm (3.24)}}}
\global\def\_@proc@ChristoffelSymbol{3.6.1}
\global\def\_@eqn@ScalarCurvature{\relax \unhbox \voidb@x \hbox {{\relax \tenrm (3.25)}}}
\global\def\_@eqn@TraceOfDerivativeOfShapeOperator{\relax \unhbox \voidb@x \hbox {{\relax \tenrm (3.26)}}}
\global\def\_@eqn@FirstTerm{\relax \unhbox \voidb@x \hbox {{\relax \tenrm (3.27)}}}
\global\def\_@eqn@SecondTerm{\relax \unhbox \voidb@x \hbox {{\relax \tenrm (3.28)}}}
\global\def\_@eqn@IntegralInequalityPEC{\relax \unhbox \voidb@x \hbox {{\relax \tenrm (3.29)}}}
\global\def\_@proc@BoundedIntegralMeanCurvaturePECTwoDim{3.6.4}
\global\def\_@proc@BoundedDiameterPECTwoDim{3.6.5}
\global\def\_@proc@PropernessPECTwoDim{3.6.6}
\global\def\_@proc@CountingPECTwoDimI{3.6.7}
\global\def\_@head@WeaklySmoothMaps{A}
\global\def\_@proc@ChainRuleI{A.0.1}
\global\def\_@proc@ChainRuleII{A.0.2}
\global\def\_@head@PrimeImmersions{B}
\global\def\_@proc@ConditionForDiffeomorphism{B.0.1}
\global\def\_@proc@Stability{B.0.2}
\global\def\_@proc@ConditionForIdentity{B.0.3}
\global\def\_@proc@Closeness{B.0.4}
\global\def\_@proc@Equicontinuity{B.0.5}
\global\def\_@head@Bibliography{C}
 }%
\else{\openout\references=references.tex }%
\fi%
%
%
\newcount\headno%
\global\headno=0%
\def\headinfo{\ifinappendices\alphanum\headno\else\the\headno\fi}%
\def\nextheadno{\global\advance\headno by 1 \global\subheadno=0 \global\procno=0 \global\eqnno=0 \headinfo}%
\makecounter{head}%
%
%
\newcount\subheadno%
\global\subheadno=0%
\def\subheadinfo{\headinfo.\the\subheadno}%
\def\nextsubheadno{\global\advance\subheadno by 1 \global\procno=0 \subheadinfo}%
\makecounter{subhead}%
%
%
\newcount\procno%
\global\procno=0%
\def\procinfo{\subheadinfo.\the\procno}%
\def\nextprocno{\global\advance\procno by 1 \procinfo}%
\makecounter{proc}%
%
%
\newcount\figno%
\global\figno=0%
\def\figinfo{\subheadinfo.\the\figno}%
\def\nextfigno{\global\advance\figno by 1 \figinfo}%
\makecounter{fig}%
%
%
\newcount\eqnno%
\global\eqnno=0%
\def\eqninfo{\text{{\rm (\headinfo.\the\eqnno)}}}%
\def\nexteqnno[#1]{\global\advance\eqnno by 1 \eqninfo\hbox{\eqnlabel{#1}}}%
\makecounter{eqn}%
%
%
%
%
%
\def\gobbleeight#1#2#3#4#5#6#7#8{}%
\newcount\citationno%
\global\citationno=0%
\def\citationinfo{\the\citationno}%
\makecounter{citation}%
\newwrite\biblio%
\def\newref#1#2{%
\def\temptext{#2}%
\edef\bibliotextoutput{\expandafter\gobbleeight\meaning\temptext}%
\global\advance\citationno by 1\citationlabel{#1}%
\ifmakebiblio%
    \edef\fileoutput{\write\biblio{\noindent\hbox to 0pt{\hss$[\the\citationno]$}\hskip 0.2em\bibliotextoutput\medskip}}%
    \fileoutput%
\fi}%
\def\cite#1{%
$[\citationref{#1}]$%
\ifmakebiblio%
    \edef\fileoutput{\write\biblio{#1}}%
    \fileoutput%
\fi%
}%
%
%
%
%
\let\mypar=\par%
\edef\Pagetitle={Blank}\headline={\hfil\Pagetitle\hfil}%
\edef\Pagefooter={Blank}\footline={\hfil\Pagefooter\hfil}%
%
%
\newcount\showpagenumflag%
\global\showpagenumflag=0 %
\def\nextoddpage%
{\newpage\ifodd\pageno%
\else\global\showpagenumflag=0 %
\null\vfil\eject%
\global\showpagenumflag=1 %
\fi}%
%
%
\font\headfont=cmb12%
\def\newhead#1[#2]%
{\ifhmode\mypar\fi%
\ifnum\headno=0 \else\goodbreak\bigskip\fi%
{\headfont\noindent\nextheadno\ - #1.}\headlabel{#2}%
\nobreak\medskip}%
%
%
\def\newsubhead#1[#2]%
{\ifhmode\mypar\fi%
\ifnum\subheadno=0 \else\goodbreak\medskip\fi%
{\bf\noindent\nextsubheadno\ - #1.\ }\subheadlabel{#2}}%
%
%
\newif\ifinproclaim%
\global\inproclaimfalse%
\def\proclaim#1{%
\goodbreak\medskip
\bgroup\inproclaimtrue%
\noindent{\bf #1}%
\nobreak\medskip\sl}%
\def\noskipproclaim#1{%
\goodbreak\medskip%
\bgroup\inproclaimtrue%
\noindent{\bf #1}\nobreak\sl}%
\def\endproclaim{\mypar\egroup\nobreak\medskip\ignorespaces}%
%
%
%
\newcount\xpos\newcount\ypos
\def\makelabelgrid{%
\xpos=-5 \ypos=-5 %
\loop\ifnum\xpos<6 %
{\loop\ifnum\ypos<6 %
\def\labeltext{x}%
\ifnum\xpos=0\def\labeltext{+}\fi%
\ifnum\ypos=0\def\labeltext{+}\fi%
\placelabel[\xpos][\ypos]{\labeltext}%
\advance\ypos by 1 %
\repeat}%
\advance\xpos by 1 %
\repeat}%
\def\placelabel[#1][#2]#3{{%
\setbox10=\hbox{\raise #2cm \hbox{\hskip #1cm #3}}%
\ht10=0pt \dp10=0pt \wd10=0pt \box10}}%
\def\placefigure#1#2#3{%
\medskip%
\midinsert%
\vbox{\line{\hfil#2\epsfbox{#3}#1\hfil}%
\vskip 0.3cm%
\line{\noindent\hfil\sl Figure \nextfigno\hfil}}%
\medskip%
\endinsert}%
%
%
\def\myitem#1{\noindent\hbox to .5cm{\hfill#1\hss}}%
%
%
%
%
%
%
%
%
%
\font\sansseriften=cmss10%
\font\sansserifseven=cmss7%
\font\sansseriffive=cmss5%
\newfam\sansseriffam%
\textfont\sansseriffam=\sansseriften%
\scriptfont\sansseriffam=\sansserifseven%
\scriptscriptfont\sansseriffam=\sansseriffive%
\def\mathsf{\fam\sansseriffam}%
%
%
%
\font\boldten=cmb10%
\font\boldseven=cmb7%
\font\boldfive=cmb5%
\newfam\mathboldfam%
\textfont\mathboldfam=\boldten%
\scriptfont\mathboldfam=\boldseven%
\scriptscriptfont\mathboldfam=\boldfive%
\def\mathbf{\fam\mathboldfam}%
%
%
%
\font\mycmmiten=cmmi10%
\font\mycmmiseven=cmmi7%
\font\mycmmifive=cmmi5%
\newfam\mycmmifam%
\textfont\mycmmifam=\mycmmiten%
\scriptfont\mycmmifam=\mycmmiseven%
\scriptscriptfont\mycmmifam=\mycmmifive%
\def\hexa#1{\ifcase #1 0\or 1\or 2\or 3\or 4\or 5\or 6\or 7\or 8\or 9\or A\or B\or C\or D\or E\or F\fi}%
\mathchardef\mathi="7\hexa\mycmmifam7B%
\mathchardef\mathj="7\hexa\mycmmifam7C%
%
%
\font\mymsbmten=msbm10 at 8pt%
\font\mymsbmseven=msbm7 at 5.6pt
\font\mymsbmfive=msbm5 at 4pt%
\newfam\mymsbmfam%
\textfont\mymsbmfam=\mymsbmten%
\scriptfont\mymsbmfam=\mymsbmseven%
\scriptscriptfont\mymsbmfam=\mymsbmfive%
\mathchardef\mybeth="7\hexa\mymsbmfam69%
\mathchardef\mygimmel="7\hexa\mymsbmfam6A%
\mathchardef\mydaleth="7\hexa\mymsbmfam6B%
%
%
%
%
\def\proof{{\noindent\bf Proof:\ }}%
\def\remark{{\noindent\bf Remark:\ }}%
\def\qed{~$\square$}%
\def\makeop#1{\global\expandafter\def\csname op#1\endcsname{{\text{#1}}}}%
\def\makeopsmall#1{\global\expandafter\def\csname op#1\endcsname{{\text{\lowercase{#1}}}}}%
%
%
\def\munion{\mathop{\cup}}%
\def\minter{\mathop{\cap}}%
%
%
\makeop{Ext}%
\makeop{Int}%
\makeop{Dist}%
\makeop{Diam}%
\makeop{Length}%
%
%
%
%
%
%
%
\def\mliminf{\mathop{{\text{LimInf}}}}%
\def\msup{\mathop{{\text{Sup}}}}%
\def\minf{\mathop{{\text{Inf}}}}%
%
%
\makeop{Dim}%
\makeop{Ker}%
\makeop{Coker}%
\makeop{Tr}%
\makeop{Adj}%
\makeop{Det}%
\makeop{End}%
\makeop{Lin}%
\makeop{Symm}%
\makeop{Mult}%
%
%
\makeop{dx}%
\makeop{dy}%
\makeop{dz}%
\makeop{dt}%
\makeop{dVol}%
\makeop{dArea}%
\makeop{Supp}%
\makeop{Hess}%
\makeop{Lip}%
%
%
\makeop{Re}%
\makeop{Im}%
\makeop{Arg}%
\makeop{Log}%
\makeop{Exp}%
%
%
\makeopsmall{Cos}%
\makeopsmall{Sin}%
\makeopsmall{Tan}%
\makeopsmall{Sec}%
\makeopsmall{Cosec}%
\makeopsmall{Cot}%
\makeopsmall{ArcCos}%
\makeopsmall{ArcSin}%
\makeopsmall{ArcTan}%
\makeopsmall{ArcSec}%
\makeopsmall{ArcCosec}%
\makeopsmall{ArcCot}%
%
%
\makeopsmall{Cosh}%
\makeopsmall{Sinh}%
\makeopsmall{Tanh}%
\makeopsmall{ArcCosh}%
\makeopsmall{ArcSinh}%
\makeopsmall{ArcTanh}%
%
%
\makeop{Vol}%
\makeop{Area}%
\makeop{Riem}%
\makeop{Ric}%
\makeop{Scal}%
\makeop{Euc}%
\makeop{Imm}%
\makeop{Emb}%
%
%
\makeop{Id}%
\makeop{Ad}%
\makeop{O}%
\makeop{SO}%
\makeop{SL}%
\makeop{GL}%
\makeop{Conf}%
\makeop{Homeo}%
\makeop{Diff}%
\makeop{Isom}%
%
%
\makeop{Ind}%
\makeop{Sig}%
\makeop{Spec}%
%
%
\makeop{Conv}%
\makeop{Max}%
\makeop{Min}%
\makeop{Mod}%
\makeop{Deg}%
\makeop{loc}%
%
%
%
%
%
%
%
%
%
%
%
%
%
 %
%
%
%
%
%
\makeop{Simp}
\makeop{Smooth}
\makeop{imm}
\makeop{Gr}
\makeop{Sign}
\makeop{equiv}
\makeop{index}
\makeop{sig}
\makeop{simp}
\makeop{ArgSup}
\def\opTI{{\Cal{T}\Cal{I}}}

\def\opTIi{{\Cal{T}_{[i]}\Cal{J}}}
\def\optildeTI{{\tilde{\Cal{T}}\Cal{I}}}
\makeop{WS}
\makeop{BWS}
\makeop{Crit}
\makeop{prime}
\def\mR{{\overline{R}}}
\def\curvR[#1#2#3#4]{{R_{#1#2#3}{}^{#4}}}
\def\curvmR[#1#2#3#4]{{\overline{R}_{#1#2#3}{}^{#4}}}
\makeop{N}

\def\mRic{{\overline{\text{R}}\text{ic}}}
\makeop{R}
\def\curvAR[#1#2#3#4]{{R^A{}_{#1#2#3}{}^{#4}}}
\newref{Aronszajn}{Aronszajn N., A unique continuation theorem for solutions of elliptic partial differential equations or inequalities of second order, {\sl J. Math. Pures Appl.}, {\bf 36}, (1957), 235--249}
\newref{Berger}{Berger M., {\sl A panoramic view of Riemannian geometry}, Springer-Verlag, Berlin, (2003)}
\newref{CaffNirSprV}{Caffarelli L., Nirenberg L., Spruck J., Nonlinear second-order elliptic equations. V. The Dirichlet problem for Weingarten hypersurfaces, {\sl Comm. Pure Appl. Math.}, {\bf 41}, (1988), no. 1, 47--70}
\newref{ColdingDeLellis}{Colding T. H., De Lellis C., The min-max construction of minimal surfaces, {\sl Surv. Differ. Geom.}, {\bf VIII}, Int. Press, Somerville, MA, (2003), 75--107}
\newref{CrandhallIshiiLyons}{Crandall M. G., Ishii H., Lions P.-L., User's guide to viscosity solutions of second order partial differential equations, {\sl Bull. Amer. Math. Soc.}, {\bf 27}, (1992), no. 1, 1--67}
\newref{ElworthyTromba}{Elworthy K. D., Tromba A. J., Degree theory on Banach manifolds, in {\sl Nonlinear Functional Analysis (Proc. Sympos. Pure Math.)}, Vol. {\bf XVIII}, Part 1, Chicago, Ill., (1968), 86--94, Amer. Math. Soc., Providence, R.I.}
\newref{EspRos}{Espinar J. M., Rosenberg H., When strictly locally convex hypersurfaces are embedded,
to appear in {\sl Math. Zeit.}}
\newref{GilbTrud}{Gilbarg D., Trudinger N. S., {\sl Elliptic partical differential equations of second order}, Die Grundlehren der mathemathischen Wissenschaften, {\bf 224}, Springer-Verlag, Berlin, New York (1977)}
\newref{GuillemanPollack}{Guillemin V., Pollack A., {\sl Differential Topology}, Prentice-Hall, Englewood Cliffs, N.J., (1974)}
\newref{Hamilton}{Hamilton R. S., The inverse function theorem of Nash and Moser, {\sl Bull. Amer. Math. Soc.}, {\bf 7}, (1982), no. 1, 65--222}
\newref{Huisken}{Huisken G., Contracting convex hypersurfaces in Riemannian manifolds by their mean curvature, {\sl Invent. Math.}, {\bf 84}, (1986), no. 3, 463--480}
\newref{Kato}{Kato T., {\sl Perturbation theory for linear operators}, Grundlehren der Mathematischen Wissenschaften, {\bf 132}, Springer-Verlag, Berlin, New York, (1976)}
\newref{LabA}{Labourie F., Probl\`emes de Monge-Amp\`ere, courbes holomorphes et laminations, {\sl Geom. Funct. Anal.}, {\bf 7}, (1997), no. 3, 496--534}
\newref{LabB}{Labourie F., Immersions isom\'etriques elliptiques et courbes pseudo-holomorphes, {\sl Geometry and topology of submanifolds} (Marseille, 1987), 131--140, World Sci. Publ., Teaneck, NJ, (1989)}
\newref{MaximoNunesSmith}{M\'aximo D., Nunes I. P., Smith G., Free boundary minimal annuli in convex three-manifolds, to appear in {\sl J. Diff. Geom.}}
\newref{Milnor}{Milnor J. W., {\sl Topology from the differential viewpoint}, Princeton Landmarks in Mathematics, Princeton, (1997)}
\newref{PacardXu}{Pacard F., Xu X., Constant mean curvature spheres in Riemannian manifolds, {\sl Manuscripta Math.}, {\bf 128}, (2009), no. 3, 275--295}
\newref{Robeday}{Robeday A., Masters Thesis, Univ. Paris VII}
\newref{RosSchneid}{Rosenberg H., Schneider M., Embedded constant curvature curves on convex surfaces, to appear in {\sl Pac. J. Math.}}
\newref{SchneiderI}{Schneider M., Closed magnetic geodesics on $S^2$, {\sl J. Differential Geom.}, {\bf 87}, (2011), no. 2, 343--388}
\newref{SchneiderII}{Schneider M., Closed magnetic geodesics on closed hyperbolic Riemann surfaces, arXiv:1009.1723}
\newref{SimonSmith}{Smith F. R., On the existence of embedded minimal $2$-spheres in the $3$-sphere, endowed with an arbitrary metric, {\sl Bull. Austral. Math. Soc.}, {\bf 28}, (1983), 159--160}
\newref{SmiAAT}{Smith G., An Arzela-Ascoli Theorem for Immersed Submanifolds, {\sl Ann. Fac. Sci. Toulouse Math.}, {\bf 16}, no. 4, (2007), 817--866}
\newref{SmiEC}{Smith G., Constant curvature hyperspheres and the Euler Characteristic, arXiv:1103.3235}
\newref{SmiPPG}{Smith G., The Plateau Problem for General Curvature Functions, arXiv:1008.3545}
\newref{SmiSLC}{Smith G., Special Lagrangian Curvature, to appear in {\sl Math. Ann.}}
\newref{SmiAC}{Smith G., Bifurcation of solutions to the Allen-Cahn equation, to appear in {\sl J. London Math. Soc.}}
\newref{Spivak}{Spivak M., {\sl A comprehensive introduction to differential geometry. Vol. III.}, Publish or Perish, Inc., Wilmington, Del., second edition, (1979)}
\newref{Tromba}{Tromba A. J., The Euler characteristic of vector fields on Banach manifolds and a globalization of Leray-Schauder degree, {\sl Adv. in Math.}, {\bf 28}, (1978), no. 2, 148--173}
\newref{WhiteI}{White B., The space of m-dimensional surfaces that are stationary for a parametric elliptic functional, {\sl Indiana Univ. Math. J.}, {\bf 36}, (1987), no. 3, 567--602}
\newref{WhiteII}{White B., Every three-sphere of positive Ricci curvature contains a minimal embedded torus, {\sl Bull. Amer. Math. Soc.}, {\bf 21}, (1989), no. 1, 71--75}
\newref{WhiteIII}{White B., Existence of smooth embedded surfaces of prescribed genus that minimize parametric even elliptic functionals on 3-manifolds, {\sl J. Differential Geom.}, {\bf 33}, (1991), no. 2, 413--443}
\newref{WhiteIV}{White B., The space of minimal submanifolds for varying Riemannian metrics, {\sl Indiana Univ. Math. J.}, {\bf 40}, (1991), no. 1, 161--200}
\newref{Ye}{Ye R., Foliation by constant mean curvature spheres, {\sl Pacific J. Math.}, {\bf 147}, (1991), no. 2, 381--396}
\catcode`\@=11
\def\multiline#1{\null\,\vcenter{\openup1\jot \m@th %
\ialign{\strut$\displaystyle{##}$\hfil\crcr#1\crcr}}\,}
\def\triplealign#1{\null\,\vcenter{\openup1\jot \m@th %
\ialign{\strut\hfil$\displaystyle{##}\quad$&\hfil$\displaystyle{{}##}$&$\displaystyle{{}##}$\hfil\crcr#1\crcr}}\,}
\def\tripleeqalign#1{\null\,\vcenter{\openup1\jot \m@th %
\ialign{\strut$\displaystyle{##}\quad$\hfil&$\displaystyle{{}##}$\hfil&$\displaystyle{{}##}$\hfil\crcr#1\crcr}}\,}
\catcode`\@=12
\myfontdefault
\global\headno=0
\global\showpagenumflag=1
\def\Pagetitle{}
\def\Pagefooter{}
\null
\vfill
\def\centre{\rightskip=0pt plus 1fil \leftskip=0pt plus 1fil \spaceskip=.3333em \xspaceskip=.5em \parfillskip=0em \parindent=0em}%
\def\textmonth#1{\ifcase#1\or January\or Febuary\or March\or April\or May\or June\or July\or August\or September\or October\or November\or December\fi}
\font\abstracttitlefont=cmr10 at 14pt
{\abstracttitlefont\centre Degree theory of immersed hypersurfaces.\par}
\bigskip
{\centre Harold Rosenberg\par}
\medskip
{\centre Graham Smith\par}
\bigskip
{\centre \the\day\ \textmonth\month\ \the\year\par}
\bigskip
\noindent{\bf Abstract:\ }We develop a degree theory for compact immersed hypersurfaces of prescribed $K$-curvature immersed in a compact, orientable Riemannian manifold, where $K$ is any elliptic curvature function. We apply this theory to count the (algebraic) number of immersed hyperspheres in various cases: where $K$ is mean curvature, extrinsic curvature and special Lagrangian curvature, and we show that in all these cases, this number is equal to $-\chi(M)$, where $\chi(M)$ is the Euler characteristic of the ambient manifold $M$.
\bigskip
\noindent{\bf Key Words:\ }Degree Theory, Immersions, Convexity, Prescribed Curvature, Non-Linear Elliptic PDEs.
\bigskip
\noindent{\bf AMS Subject Classification:\ }58D10 (58B05, 58C40, 58J05)
%
%
\par
\nextoddpage
\global\showpagenumflag=0
\global\pageno=1
\def\Pagetitle{\hfil Degree theory of immersed hypersurfaces.\hfil}
\def\Pagefooter{\hfil{\myfontdefault\folio}\hfil}
\newhead{Introduction}[Introduction]
\newsubhead{General}[General]
Let $M:=M^{d+1}$ be a closed, orientable, $(d+1)$ dimensional riemannian manifold. Our aim is to establish the existence of closed, orientable, immersed hypersurfaces in $M$ of prescribed $K$-curvature, where $K$ is an elliptic curvature function, such as mean curvature, extrinsic curvature, and so on. We achieve this by constructing a degree theory for the space of unparametrised immersions and by then relating this degree to the topology of $M$. Our work is inspired by the beautiful degree theory developed by Brian White for the study of hypersurfaces of prescribed mean curvature and solutions of other parametric elliptic functionals (c.f \cite{WhiteI}, \cite{WhiteII}, \cite{WhiteIII} and \cite{WhiteIV}).
\par
The theory that will be developed in the sequel may be summarised as follows. Let $\Sigma:=\Sigma^d$ be a $d$ dimensional manifold, let $\Cal{I}$ be an open set of unparametrised immersions of $\Sigma$ into $M$, let $\Cal{O}$ be an open set of smooth functions defined over $M$, and let $K$ be an elliptic curvature function. We will be interested in those immersions, $e:\Sigma\rightarrow M$, whose $K$-curvature is prescribed by some function, $f\in\Cal{O}$, that is, such that $K_e=f\circ e$, where $K_e$ here denotes the $K$-curvature of the immersion, $e$. We thus define the solution space,
$$
\Cal{Z} := \left\{ ([e],f)\in\Cal{I}\times\Cal{O}\ |\ K_e = f\circ e\right\},
$$
and we define $\Pi:\Cal{Z}\rightarrow\Cal{O}$ to be the projection onto the second factor. A fundamental hypothesis will be that the projection, $\Pi$, is a proper map. Under these conditions, that is, when $K$ is elliptic and when $\Pi$ is proper, we will show that $\Pi$ has a well defined, integer valued mapping degree. More precisely, we will show that $\Cal{O}$ has a dense subset, $\Cal{O}'$, such that for all $f\in\Cal{O}'$, the preimage, $\Pi^{-1}(\left\{f\right\})$, is finite, and each element, $([e],f)$, of this set has a well defined signature, equal to plus or minus one. For all such $f$ we then define
$$
\opDeg(\Pi;f) := \sum_{([e],f)\in\Pi^{-1}(\left\{f\right\})}\opsig([e],f).
$$
By proving that $\opDeg(\Pi;f)$ is independent of $f$, we obtain a mapping degree for $\Pi$, and by calculating this degree in certain special cases, we obtain existence results for hypersurfaces of prescribed $K$-curvature in various settings.
\newsubhead{Background}[Background]
In order to gain a clear intuition for our construction, it is worth reviewing the analogous finite dimensional theory. Thus, let $X$ be an open subset of $\Bbb{R}^d$, let $Y$ be a finite dimensional manifold, and let $f:X\times Y\rightarrow\Bbb{R}^d$ be a smooth map, which we think of as a family of vector fields over $X$ parametrised by $Y$. Let $Z:=f^{-1}(\left\{0\right\})$ be the zero set of $f$, and let $\Pi:Z\rightarrow Y$ be the projection onto the second factor. When $f$ is a submersion and when $\Pi$ is proper, classical differential topology (c.f. \cite{Milnor} and \cite{GuillemanPollack}) shows that, for generic $y$, the preimage, $\Pi^{-1}(\left\{y\right\})$ is finite. Furthermore, each element, $(x,y)$, of this set has a well defined signature, equal to plus of minus one, and the sum of these signatures is equal to the mapping degree of $\Pi$. However, for any fixed $y$, the set $\Pi^{-1}(\left\{y\right\})$ is also the zero set of the vector field, $f_y:=f(\cdot,y)$. In particular, for generic $y$, this vector field has finitely many zeroes, the number of which, when counted with appropriate sign, is equal to the mapping degree of $\Pi$.
\par
There are therefore two alternative approaches to studying the degree. The first involves determining the topological mapping degree of $\Pi$, and the second involves counting the zeroes of the vector field, $f_y$, for generic data, $y$. Although these two approaches are equivalent, yielding exactly the same values, they lead to complementary formal developments, especially in the infinite dimensional setting.
\par
White's degree theory mirrors the first approach. Here, the space, $\Cal{I}$, of unpara\-metrised immersions plays the role of $X$, an open subset of the space of riemannian metrics of the underlying manifold plays the role of $Y$, and the mean curvature functional plays the role of $f$. The objects of interest are minimal surfaces, which are zeroes of the mean curvature functional. In fact, White views minimal surfaces as critical points of the area functional, but since the mean curvature functional is the $L^2$-gradient of the area functional, this really comes down to the same thing. White then proves that the solution space, $Z$, has the structure of a Banach manifold and that the projection, $\Pi$, is a smooth Fredholm map of Fredholm index $0$. Having established this framework, White then applies the Sard-Smale theorem together with spectral properties of the Jacobi operator to construct a $\Bbb{Z}$-valued degree theory, which yields existence results for minimal surfaces in various settings.
\par
The theory that we will present here is similar to that of White, but mirrors the second approach. Here, the space, $\Cal{I}$, of unparametrised immersions continues to play the role of $X$, but - for the most part - the space, $\Cal{O}$, of functions plays the role of $Y$, and more general curvature functionals play the role of $f$. A notable technical difficulty of this theory, which tends to generate unpleasant notational complexity, is that $\Cal{I}$ does not actually carry the structure of a Banach manifold. In White's work, this difficulty is circumvented by restricting attention to the solution space, $Z$, which, on the contrary, does possess such a structure. In the current paper we propose an alternative resolution of this problem by introducing the elementary concept of ``weakly smooth'' manifolds (c.f. Appendix \headref{WeaklySmoothMaps}). It is straightforward to show that $\Cal{I}$ belongs to this class (c.f. Section \subheadref{TheManifoldOfImmersions}), and we may therefore use the familiar terminology of manifold theory without problems.
\par
A more fundamental difficulty of the present case is, however, that the vector fields of interest to us are no longer necessarily $L^2$-gradients of real valued functionals. This makes the definition of the signature more involved than in White's case. Its construction, which will be carried out in Sections \subheadref{DefiningTheSignatureFD} and \subheadref{DefiningTheSignatureID}, is one of the main contributions of the current paper: using the spectral theory of non-self adjoint operators (c.f. \cite{Kato}), we extend the concept of the sign of the determinant to operators acting on infinite dimensional spaces. This yields a well defined signature for zeroes of certain types of vector fields over $\Cal{I}$. Upon summing these signatures, we obtain the desired $\Bbb{Z}$-valued degree, and this in turn yields existence results in various settings.
\par
Section \headref{DegreeTheory} of this paper will be devoted to constructing this degree. We first introduce the manifold structure of the space of (prime) unparametrised immersions of $\Sigma$ into $M$, showing, in particular, how curvature functionals define vector fields over this manifold. Next, using the spectral properties of the linearised curvature operators, we show how the zeroes of these vector fields have well defined signatures, and, by counting these signatures, we obtain the mapping degree for $\Pi$.
\par
This paper was initially motivated by the question of the existence of embeddings of constant mean curvature of $\Sigma=S^n$ into $(S^{n+1},g)$. Indeed, we conjecture that for any $c\geqslant 0$ and any metric $g$ on $S^3$ of positive sectional curvature, there is an embedding of $S^2$ into $S^3$ of constant mean curvature $c$. This result is known for $c=0$ (c.f. \cite{SimonSmith} and \cite{ColdingDeLellis}) and also for large values of $c$ (c.f. \cite{Ye}). However, even when $n=1$ and $g$ is a positive curvature metric on $S^2$, we do not know if there exist embeddings of $S^1$ having any prescribed, positive geodesic curvature. Nonetheless, Anne Robeday (c.f. \cite{Robeday}) and, independently, Matthias Schneider (c.f. \cite{SchneiderI}) have proven that, under these conditions on the metric there does exist an Alexandrov embedding (i.e. a not necessarily embedded immersion that nonetheless extends to an immersion of the disk) of $S^1$ into $(S^2,g)$ of any given constant geodesic curvature. Anne Robeday's approach used the degree theory of Brian White to prove this result whilst Schneider developed a different degree theory for his proof. Schneider's theory, which applies to immersions of the circle into any riemannian surface, has yielded many other interesting results (c.f. \cite{SchneiderII}). Additionally, in \cite{RosSchneid}, the first author and M. Schneider have also proven that given a metric of positive curvature on $S^2$, there is an $\epsilon>0$, such that for any $c\in]0,\epsilon[$, there are at least two embeddings of $S^1$ into $(S^2,g)$ of geodesic curvature equal to $c$.
\newsubhead{Applications}[Applications]
In Section \headref{Applications} we give applications of our degree theory. First, given a family $\Cal{G}$ of riemannian metrics over $M$, we say that a property holds for {\sl generic} $g\in\Cal{G}$ whenever it holds over a subset $\Cal{G}_0$ of $\Cal{G}$ which itself contains a countable union of dense open sets. We recall that such a subset is said to be of the second category in the sense of Baire and that, by the Baire Category Theorem, it is dense in $\Cal{G}$. We now prove four theorems which count for generic metrics, and under appropriate hypotheses, the algebraic number of immersions of constant curvature of $S^d$ into $(M,g)$. In each case, we will see that this algebraic number (which is the degree of $\Pi$) is equal to $-\chi(M)$, where $\chi(M)$ is the Euler characteristic of $M$. That this should be the case can be readily explained as follows. In the case where the scalar curvature function $R$ of $M$ is a Morse function, Ye proved in \cite{Ye} that a punctured neighbourhood of each critical point of $R$ is foliated by a family of constant mean curvature spheres $\Sigma(h)$ where $h$ varies over an interval $]h_0,\infty[$, where $h_0$ is large and depends on $(M,g)$. Ye's result readily extends to general notions of curvature and it can be shown that the signature of a sphere constructed in this manner is equal to $(-1)^d$ times the signature of the corresponding critical point of $R$. Thus, if $h\geq h_0$ is a regular value of $\Pi$, and if Ye's spheres account for every element of $\Pi^{-1}(h)$, then $\opDeg(\Pi;H)=(-1)^d\chi(M)=-\chi(M)$. This argument will be made precise in Section \subheadref{CalculatingTheDegree}, and as one may expect, the main difficulty then lies in showing that $\Pi$ is a proper map.
\par
Our first result concerns the case where $K=H$ is mean curvature. Consider first the polynomial
$$
p(t) := \frac{t(t-1)(t-d)}{(d-1)}.
$$
Denote
$$\eqalign{
\Gamma &:= \msup_{t\in[0,1]}p(t),\ \text{and}\cr
c_0 &:= \opArgSup_{t\in[0,1]}p(t),\cr}
$$
and let $\phi:]0,\Gamma[\rightarrow]0,c_0[$ be the inverse of $\phi$ over the interval $]0,c_0[$. For a given metric $g$ over $M$ denote
$$
h_0 := h_0(g) := \minf\left\{ h\ \left|\ \frac{2\|\overline{\nabla}\overline{R}\|}{h^3} + \frac{4\|\overline{R}^0\|}{h^2} + \frac{\|\overline{R}\|}{h^2}<\Gamma\right.\right\},
$$
where $\overline{R}$ is the Riemann curvature tensor of $(M,g)$, $\overline{\nabla}\overline{R}$ is its covariant derivative and $\overline{R}^0$ is its trace free component. Define $c:]h_0,\infty[\rightarrow]0,c_0[$ by
$$
c(h) := c(h,g) := \phi\left(\frac{2\|\overline{\nabla}\overline{R}\|}{h^3} + \frac{4\|\overline{R}^0\|}{h^2} + \frac{\|\overline{R}\|}{h^2}\right).
$$
Now, an immersion $e:S^d\rightarrow M^{d+1}$ is said to be pointwise $c(h)$-pinched whenever, for each $p\in S^d$,
$$
\lambda_1(p) > \frac{c(h)}{d}(\lambda_1 + ... + \lambda_d)(p) = c(h)H(p),
$$
where $\lambda_1\leq...\leq\lambda_d$ are its principal curvatures. In particular any such immersion is always locally convex in the sense that its shape operator is non-negative definite. We now show (c.f. Theorem \procref{ExistencePMCI})
\proclaim{Theorem \nextprocno}
\noindent For a generic riemannian metric $g$ over $M$, and for all $h>h_0(M)$, the algebraic number of $c(h)$-pinched embedded hyperspheres in $M$ of constant mean curvature equal to $h$ is itself equal to $-\chi(M)$.
\endproclaim
\proclabel{MainTheoremMeanCurvature}
\remark To see that Theorem \procref{MainTheoremMeanCurvature} constitutes an improvement over Ye's result, observe that in the case of \cite{Ye}, the number $h_0$ will depend on the second derivative of the Riemann curvature function, wheras in our case it only depends on the first.
\medskip
Our next theorem concerns the case where $K:=K_e$ is extrinsic curvature (also referred to as Gauss-Kr\"onecker curvature). We say that $(M,g)$ is {\sl pointwise $1/2$-pinched} whenever
$$
\sigma_\opMax(p) < 2\sigma_\opMin(p)
$$
for every point $p\in M$, where here $\sigma_\opMax(p)$ and $\sigma_\opMin(p)$ denote respectively the maximum and minimum sectional curvatures of $M$ at $p$. Recall also that $(M,g)$ is said to be {\sl $1/4$-pinched} whenever
$$
\msup_{p\in M}\sigma_\opMax(p) < 4\minf_{p\in M}\sigma_\opMin(p).
$$
We prove (c.f. Theorem \procref{CountingPEC})
\proclaim{Theorem \nextprocno}
\noindent For generic $g$ such that $(M,g)$ is $1/4$-pinched and pointwise $1/2$-pinched, and for all $k>0$, the algebraic number of locally strictly convex embedded hypersurfaces in $M$ of constant extrinsic curvature equal to $k$ is itself equal to $-\chi(M)$.
\endproclaim
\proclabel{MainTheoremExtrinsicCurvature}
Our third application concerns the case where $K=R_\theta$ is {\sl special lagrangian curvature} (c.f. \cite{SmiSLC} for a detailled definition). We prove (c.f. Theorem \procref{CountingPSLC})
\proclaim{Theorem \nextprocno}
\noindent If $\theta\in[(d-1)\pi/2,d\pi/2[$, then, for generic $g$ such that $(M,g)$ is $1/4$-pinched, and for all $r>0$, the algebraic number of locally strictly convex embedded hyperspheres in $M$ of constant special lagrangian curvature equal to $r$ is itself equal to $-\chi(M)$.
\endproclaim
\proclabel{MainTheoremSpecialLagrangianCurvature}
\noindent In the special case where $d=3$, Theorem \procref{MainTheoremSpecialLagrangianCurvature} becomes
\proclaim{Corollary \nextprocno}
\noindent If $M$ is $4$-dimensional, then for generic $g$ such that $(M,g)$ is $1/4$-pinched, and for all $r>0$, the algebraic number of locally strictly convex embedded hyperspheres in $M$ such that $K_e=rH$ is itself equal to $-2$, where here $K_e$ and $H$ denote respectively the extrinsic curvature and the mean curvature functions of the hypersphere.
\endproclaim
Our final theorem (and which, for us, is the deepest result) concerns the case where $K=K_e$ is the extrinsic curvature of a surface in a $3$-dimensional manifold. Let $M$ be a $3$-dimensional manifold, and for a given metric $g$ over $M$, denote
$$
k_0^2 = \frac{1}{2}\left(\left|\sigma_{\opMin}^-\right| + \sqrt{\left|\sigma_{\opMin}^-\right| + \|T\|_O^2}\right),
$$
where $T$ is the trace free Ricci curvature of $M$, $\|T\|_O$ is its operator norm, viewed as an endomorphism of $TM$, and $\sigma_\opMin^-$ is defined to be $0$, or the infimum of the sectional curvatures of $M$, whichever is lower. We now obtain (c.f. Theorem \procref{CountingPECTwoDimI})
\proclaim{Theorem \nextprocno}
\noindent If $M$ is $3$-dimensional, then for generic $g$ and for all $k>k_0$, the algebraic number of locally strictly convex immersed spheres in $(M,g)$ of constant extrinsic curvature equal to $k$ is itself equal to $0$.
\endproclaim
\proclabel{MainTheoremExtrinsicCurvatureTwoD}
\newsubhead{Acknowledgements}[Acknowledgements]
We would like to thank Bernard Helffer for detailled correspondence concerning spectral theory. The second author would like to thank IMPA, Rio de Janeiro, and CRM, Barcelona, Spain for providing the conditions required to prepare this work. An earlier draft of this paper was prepared whilst the second author was benefitting from a Marie Curie postdoctoral fellowship.
\newhead{Degree theory}[DegreeTheory]
\newsubhead{The manifold of immersions and its tangent bundle}[TheManifoldOfImmersions]
Let $M:=M^{d+1}$ be a compact, oriented, $(d+1)$ dimensional, riemannnian manifold. Let $\Sigma:=\Sigma^d$ be a compact, oriented, $d$ dimensional manifold without boundary. Consider a smooth immersion, $i:\Sigma\rightarrow M$. When $\Sigma$ is simply connected, we say that $i$ is {\sl prime} whenever there exists no non-trivial diffeomorphism, $\alpha:\Sigma\rightarrow\Sigma$, such that $i\circ\alpha=i$, that is, $i$ is prime whenever it is not a multiple cover. More generally, let $\tilde{\Sigma}$ be the universal cover of $\Sigma$ and let $p:\tilde{\Sigma}\rightarrow\Sigma$ be the canonical projection. We say that $i$ is {\sl prime} whenever the composition, $\tilde{\mathi}:=i\circ p$, is invariant only under the action of deck transformations of $\tilde{\Sigma}$, that is, if $\alpha:\tilde{\Sigma}\rightarrow\tilde{\Sigma}$ is a diffeomorphism such that $\tilde{\mathi}\circ\alpha=\tilde{\mathi}$, then $\alpha\in\pi_1(\Sigma)$. We introduce prime immersions in order to ensure smoothness of the manifold of immersions. We refer the reader to Appendix \headref{PrimeImmersions} for a brief discussion of some of their properties.
\par
Let $C^\infty_\opprime(\Sigma,M)$ denote the space of smooth, prime immersions from $\Sigma$ into $M$. Let $\opDiff^+(\Sigma)$ denote the group of smooth, orientation preserving diffeomorphisms of $\Sigma$. This group acts on $C^\infty_\opprime(\Sigma,M)$ by composition, and we denote the quotient space by $\Cal{I}:=\Cal{I}(\Sigma,M)$, that is,
$$
\Cal{I} := C^\infty_\opprime(\Sigma,M)/\opDiff^\infty(\Sigma).\eqnum{\nexteqnno[DefinitionOfI]}
$$
We furnish this space with the quotient topology, and we call it the space of {\sl unparametrised (prime) immersions}. Given an element $i\in C^\infty_\opprime(\Sigma,M)$, we denote by $[i]$ its equivalence class in $\Cal{I}$. Throughout the sequel, we will often identify the equivalence class, $[i]$, with its representative element, $i$.
\par
For $[i]\in\Cal{I}$, let $N_i:\Sigma\rightarrow TM$ be the unit, normal vector field over $i$ which is compatible with the orientation, and define $\Cal{E}_i:\Sigma\times\Bbb{R}\rightarrow M$ by
$$
\Cal{E}_i(x,t) := \opExp_{i(x)}(tN_i(x)),
$$
where $\opExp$ is the exponential map of $M$. Let $\epsilon_i>0$ be such that the restriction of $\Cal{E}_i$ to $\Sigma\times]-\epsilon_i,\epsilon_i[$ is an immersion (c.f. Proposition $3.2$ of \cite{MaximoNunesSmith}), and define the open subset, $\Cal{U}_i$, of $C^\infty(\Sigma)$ and the map, $\hat{\Phi}_i:\Cal{U}_i\rightarrow C^\infty_\opimm(\Sigma,M)$, by
$$\eqalign{
\Cal{U}_i &:= \left\{ \phi\in C^\infty(\Sigma)\ |\ \|\phi\|_{C^0} < \epsilon_i\right\},\ \text{and}\cr
\hat{\Phi}_i(\phi)(x) &:= \Cal{E}_i(x,\phi(x)).\cr}
$$
By Appendix \headref{PrimeImmersions}, upon reducing $\epsilon_i$ if necessary, we may suppose that $\hat{\Phi}_i(\phi)$ is prime for all $\phi\in\Cal{U}_i$, so that $\hat{\Phi}_i$ projects down to a map, $\Phi_i$, from $\Cal{U}_i$ into $\Cal{I}$. By Propositions $3.3$ and $3.4$ of \cite{MaximoNunesSmith}, this map is a homeomorphism onto an open subset, $\Cal{V}_i$, of $\Cal{I}$. We call the triplet, $(\Phi_i,\Cal{U}_i,\Cal{V}_i)$, the {\sl graph chart} of $\Cal{I}$ about $i$. The set of all graph charts constitutes an atlas of $\Cal{I}$, all of whose transition maps are homeomorphisms, thus furnishing $\Cal{I}$ with the structure of a topological manifold modelled on $C^\infty(\Sigma)$.
\par
Furnishing $\Cal{I}$ with a smooth manifold structure is a non-trivial matter, as the usual approaches each present their own difficulties. Indeed, on the one hand, $\Cal{I}$ has no Banach manifold structure since its transition maps are differentiable neither with respect to any H\"older norm nor with respect to any Sobolev norm. On the other hand, although $\Cal{I}$ does carry the structure of a smooth, tame Frechet manifold (c.f. \cite{Hamilton}), we find this formalism inelegant in terms of its high technical cost in relation to the problem at hand. For these reasons, we introduce the following terminology which, though elementary, possesses all the structure required for the development of our theory (c.f. Appendix \headref{WeaklySmoothMaps} for a formal treatment).
\par
Consider first a smooth, compact, finite dimensional manifold, $X$. Given another smooth, finite dimensional manifold, $N$, and a map, $f:N\rightarrow C^\infty(X)$, we define $\tilde{f}:N\times X\rightarrow\Bbb{R}$ by $\tilde{f}(p,x):=f(p)(x)$. We say that $f$ is {\sl strongly smooth} whenever $\tilde{f}$ is smooth. Now consider another smooth, compact, finite dimensional manifold, $X'$. Given an open subset, $\Cal{U}$, of $C^\infty(X)$, we say that a map, $\Phi:\Cal{U}\rightarrow C^\infty(X')$, is {\sl weakly smooth} whenever the action of composition, $f\mapsto\Phi\circ f$, sends strongly smooth maps into strongly smooth maps, and is continuous with respect to the $C^\infty_\oploc$ topology. A {\sl weakly smooth manifold}, $\Cal{M}$, modelled on $C^\infty(X)$, is now defined to be a Hausdorff topological space furnished with an atlas, all of whose charts are open subsets of $C^\infty(X)$, and all of whose transition maps are weakly smooth. It is straightforward to verify that the above defined atlas furnishes $\Cal{I}$ with the structure of a weakly smooth manifold modelled on $C^\infty(\Sigma)$.
\par
Weakly smooth manifolds have well defined tangent bundles. Indeed, consider again a compact, finite dimensional manifold, $X$. We say that two strongly smooth curves, $c,c':]-\epsilon,\epsilon[\rightarrow C^\infty(X)$, such that $c(0)=c'(0)=:\phi$ are {\sl equivalent} at $\phi$ whenever $\partial_t\tilde{c}|_{t=0}=\partial_t\tilde{c}'|_{t=0}$. Now consider a weakly smooth manifold, $\Cal{M}$. Since equivalence of strongly smooth curves is preserved by weakly smooth maps, we say that two strongly smooth curves, $c,c':]-\epsilon,\epsilon[\rightarrow\Cal{M}$, such that $c(0)=c'(0)=:x$ are {\sl equivalent} at $x$ whenever they are equivalent in every graph chart. {\sl Tangent vectors} over $\Cal{M}$ are now defined to be equivalence classes of strongly smooth curves, and the {\sl tangent bundle} of $\Cal{M}$, denoted by $\Cal{TM}$, is defined to be the set of all tangent vectors. $\Cal{TM}$ carries the structure of a vector bundle over $\Cal{M}$ with typical fibre $C^\infty(X)$ whose total space is a weakly smooth manifold modelled on $C^\infty(X)\times C^\infty(X)=C^\infty(X\sqcup X)$.
\par
Weakly smooth maps between weakly smooth manifolds also have well defined derivatives. Indeed, consider two weakly smooth manifolds, $\Cal{M}$ and $\Cal{N}$, and a weakly smooth map, $\Phi:\Cal{M}\rightarrow\Cal{N}$. Given $x\in\Cal{M}$, the {\sl derivative} of $\Phi$ at $x$ is given by
$$
D\Phi(x)\cdot[c] := [\Phi\circ c],
$$
where $c:]-\epsilon,\epsilon[\rightarrow\Cal{M}$ is a strongly smooth curve with $c(0)=x$, $[c]$ is the vector defined by its equivalence class at $x$, and $[\Phi\circ c]$ is the vector defined by the equivalence class of $\Phi\circ c$ at $\Phi(x)$. $D\Phi$ defines a weakly smooth map from $\Cal{T}\Cal{M}$ into $\Cal{T}\Cal{N}$ which is linear over each fibre.
\par
In order to define sections of the tangent bundle, $\opTI$, it is useful to see how this bundle also arises as an associated bundle. Indeed, using the action of $\opDiff^+(\Sigma)$ by composition on $C^\infty(\Sigma)$, we define the associated bundle,
$$
\optildeTI := C^\infty_\opprime(\Sigma,M)\times C^\infty(\Sigma)/\opDiff^+(\Sigma).
$$
This is a topological vector bundle over $\Cal{I}$ with typical fibre $C^\infty(\Sigma)$. Given an element, $(i,f)$, of $C^\infty_\opprime(\Sigma,M)\times C^\infty(\Sigma)$, we denote its equivalence class in $\opTI$ by $[i,f]$. Given a graph chart, $(\Phi_i,\Cal{U}_i,\Cal{V}_i)$, of $\Cal{I}$, for any $\phi\in\Cal{U}_i$ and for any $f\in C^\infty(\Sigma)$, we define the push-forward of $(\phi,f)$ through $\Phi_i$ by
$$
\Phi_{i,*}(\phi,f) = [\hat{\Phi}_i(\phi),f].
$$
This yields a homeomorphism from $\Cal{U}_i\times C^\infty(\Sigma)$ into $\optildeTI|_{\Cal{V}_i}$ which is linear over each fibre, that is, a trivialisation of $\optildeTI$. In particular, the collections of all such trivialisations furnishes the total space of $\optildeTI$ with the structure of a weakly smooth manifold.
\par
Now choose $[i]\in\Cal{I}$ and let $(\Phi_i,\Cal{U}_i,\Cal{V}_i)$ be the graph chart of $\Cal{I}$ about $i$. For $f\in C^\infty(\Sigma)$, define the strongly smooth curve, $c_f:]-\epsilon,\epsilon[\rightarrow\Cal{I}$, by
$$
c_f(t) := \Phi_i(tf).
$$
By identifying the element $[i,f]$ of $\optildeTI$ with the tangent vector, $[c_f]$, we obtain a weakly smooth bundle isomorphism from $\optildeTI$ into $\opTI$, and we henceforth identify these two spaces. In geometric terms, the element, $[i,f]$, of $\optildeTI$ thus identifies with the perturbation of $i$ along the normal direction, $fN_i$, where $N_i$ is the unit, normal vector field over $i$ compatible with the orientation, which is the usual way in which functions over $\Sigma$ are identified with tangent vectors of $\Cal{I}$.
\par
Finally, there is a natural covariant derivative for weakly smooth sections of $\opTI$. Indeed, let $\sigma$ be a weakly smooth section of $\opTI$ over $\Cal{I}$ and let $[i,f]$ be a tangent vector to $\Cal{I}$ at the point $[i]$. Let $(\Phi_i,\Cal{U}_i,\Cal{V}_i)$ be the graph chart of $\Cal{I}$ about $i$, let $\sigma_i:\Cal{U}_i\rightarrow C^\infty(\Sigma)$ be the pull-back of $\sigma$ through $\Phi_i$, and let $c:]-\epsilon,\epsilon[\rightarrow\Cal{U}_i$ be a strongly smooth curve such that $\partial_t c|_{t=0}=f$. The {\sl derivative} of $\sigma$ in the direction of $[i,f]$ is now defined by
$$
D\sigma([i])\cdot[i,f] := [i,\partial_t(\sigma_i\circ c_t)|_{t=0}],
$$
and this is well defined, independently of the strongly smooth curve, $c$, chosen.
\newsubhead{Curvature as a vector field}[CurvatureAsAVectorField]
We recall the definition of a curvature function (c.f. \cite{CaffNirSprV} and \cite{SmiPPG}). Let the symmetric group, $S_d$, act on vectors in $\Bbb{R}^d$ by permutations of their components. Let $\Lambda\subseteq\Bbb{R}^d$ be an open, convex cone based on the origin which is invariant under this action. A smooth function, $K:\Lambda\rightarrow]0,\infty[$, is said to be a {\sl curvature function} whenever it is also invariant under this action, is homogeneous of order $1$, and is {\sl elliptic} in the sense that its gradient at any point of $\Lambda$ is given by a vector all of whose components are positive. Observe, in particular, that a curvature function is always defined along with its cone of definition.
\par
When performing computations, it is easier to consider curvature functions as functions defined over subsets of the space, $\opSymm(\Bbb{R}^d)$, of $d$-dimensional symmetric matrices. Thus, let $\Lambda$ also denote the open, convex cone in $\opSymm(\Bbb{R}^d)$ consisting of those matrices whose vectors of eigenvalues are elements of $\Lambda$, and for any such matrix, $A$, with eigenvalues, $\lambda_1,...,\lambda_d$, define
$$
K(A) := K(\lambda_1,...,\lambda_d).
$$
By $S_d$-invariance, both $\Lambda$ and $K$ are well-defined.
\par
Given a curvature function, $K$, with cone of definition $\Lambda$, an immersion, $i$, is said to be {\sl $K$-convex} whenever its shape operator, $A_i$, is at every point an element of $\Lambda$. We henceforth suppose that $\Cal{I}$ consists only of $K$-convex immersions. The {\sl $K$-curvature} of any element, $i$, of $\Cal{I}$ is then defined by
$$
K_i(p):=K(A_i(p)).\eqnum{\nexteqnno[DefinitionOfKi]}
$$
The $K$-curvature defines a weakly smooth section of $\opTI$ over $\Cal{I}$ which we also denote by $K$, that is, for all $[i]\in\Cal{I}$,
$$
K([i]) := [i,K_i].\eqnum{\nexteqnno[DefinitionOfTheSectionK]}
$$
The covariant derivative of this section is given by
$$
DK([i])\cdot[i,f] = [i,J_if],\eqnum{\nexteqnno[CovariantDerivativeOfK]}
$$
where $[i,f]$ is a tangent vector to $\Cal{I}$ at $[i]$, and $J_i$ is the {\sl Jacobi operator} of the $K$-curvature at $i$. Recall that $J_i$ is a linear, second-order partial differential operator with smooth coefficients. The ellipticity of $K$ now implies the following key property.
\medskip
{\noindent\bf Ellipticity:\ }For every $[i]\in\Cal{I}$, the Jacobi operator of $K$ at $i$ is elliptic.
\medskip
The zeroes of $K$ are precisely those immersions of constant $K$-curvature equal to $0$. More generally, let $\Cal{O}\subseteq C^\infty(M)$ be an open subset of smooth functions, and define the evaluation functional, $E:\Cal{I}\times\Cal{O}\rightarrow C^\infty(\Sigma)$, by
$$
E(i,f) := f\circ i.\eqnum{\nexteqnno[DefinitionOfEi]}
$$
$E$ defines a weakly smooth family of weakly smooth sections of $\opTI$ over $\Cal{I}$ parametrised by $\Cal{O}$, which we also denote by $E$, that is, for all $([i],f)\in\Cal{I}\times\Cal{O}$,
$$
E([i],f) = [i,f\circ i].\eqnum{\nexteqnno[DefinitionOfTheSectionE]}
$$
The partial derivatives of $E$ are given by
$$\eqalign{
D_1E([i],f)\cdot [i,g] &:= [i,\langle \nabla f, N_i\rangle g],\ \text{and}\cr
D_2E([i],f)\cdot h &:= [i,h\circ i],\cr}\eqnum{\nexteqnno[CovariantDerivativeOfE]}
$$
where $D_1$ and $D_2$ denote the partial derivatives along $\Cal{I}$ and $\Cal{O}$ respectively, and $N_i$ denotes the unit, normal vector field over $i$ which is compatible with the orientation.
\par
We now define $\hat{K}:\Cal{I}\times\Cal{O}\rightarrow\opTI$ by
$$
\hat{K}([i],f) := K([i]) - E([i],f).\eqnum{\nexteqnno[DefinitionOfKHati]}
$$
Its partial derivatives are given by
$$\eqalign{
D_1\hat{K}([i],f)\cdot [i,g] &= [i,\hat{J}_{i,f}g],\cr
D_2\hat{K}([i],f)\cdot h &= [i,-h\circ i],\cr}\eqnum{\nexteqnno[CovariantDerivativeOfKhati]}
$$
where the operator
$$
\hat{J}_{i,f}g := J_ig - \langle\nabla f,N_i\rangle g,\eqnum{\nexteqnno[JacobiOperatorOfKHat]}
$$
will be called the {\sl Jacobi operator} of $\hat{K}$ at $(i,f)$. As before, $\hat{J}_{i,f}$ is a second-order, partial differential operator with smooth coefficients, and since $\hat{J}_i$ is elliptic, so too is $\hat{J}_{i,f}$.
\medskip
{\bf\noindent The solution space:\ }The {\sl solution space}, $\Cal{Z}\subseteq\Cal{I}\times\Cal{O}$, is now defined to be the zero set of $\hat{K}$, that is,
$$
\Cal{Z} := \hat{K}^{-1}(\left\{0\right\}).\eqnum{\nexteqnno[DefinitionOfZ]}
$$
In other words, $\Cal{Z}$ is the set of all pairs, $([i],f)$, such that the $K$-curvature of $i$ is prescribed by $f$, that is, such that
$$
K_i = f\circ i.
$$
Let $\Pi:\Cal{Z}\rightarrow\Cal{O}$ be the projection onto the second factor. We henceforth suppose
\medskip
{\bf\noindent Properness:\ } $\Pi$ defines a proper map from $\Cal{Z}$ into $\Cal{O}$.
\medskip
\noindent For all $f\in\Cal{O}$, denote
$$
\Cal{Z}_f := \Pi^{-1}(\left\{f\right\}) = \left\{[i]\ |\ K_i = f\circ i\right\}.\eqnum{\nexteqnno[SolutionSetOverF]}
$$
We say that $f$ is a {\sl regular value} of $\Pi$ whenever the Jacobi operator, $\hat{J}_{i,f}$, is invertible at every point $[i]$ of $\Cal{Z}_f$. In this case, $\Cal{Z}_f$ is discrete and, by compactness, it is therefore finite.
\par
We will show presently that almost every $f$ is a regular value of $\Pi$, and the degree will then be defined to be the sum of the parities of the critical points of $\Cal{Z}_f$. Considering this as the number of zeroes counted with parity of the vector field, $K_f:=K(\cdot,f)$, we obtain a theory of degrees for vector fields compatible with that developed by Schneider in \cite{SchneiderI} and \cite{SchneiderII}, which itself builds on the earlier work, \cite{ElworthyTromba} and \cite{Tromba}, of Elworthy and Tromba. Interestingly, this perspective also yields the formal intepretation of the degree as the Euler characteristic of the manifold $\Cal{I}$, even though the latter does not - strictly speaking - have a clear mathematical definition.
\newsubhead{Simplicity}[Simplicity]
Recall that an immersion, $i:\Sigma\rightarrow M$, is said to be {\sl simple} whenever it has the property that, for any two distinct points, $p$ and $q$, of $\Sigma$, and for all sufficiently small neighbourhoods, $U_p$ and $V_q$, of $p$ and $q$ respectively, we have $i(U_p)\neq i(V_q)$. Every simple immersion is trivially prime, but the converse does not, in general, hold. However, in the case at hand, the two are indeed equivalent for immersions that lie in the solution space. To see this, first, for $i\in C^\infty_\opprime(\Sigma,M)$, let $J^\infty(i)$ be its $C^\infty$-jet.
\proclaim{Lemma \nextprocno}
\noindent For all $([i],f)\in\Cal{Z}$, $J^\infty(i)$ is injective.
\endproclaim
\proclabel{PropUniqueContinuation}
\proof Suppose the contrary. Choose $p\neq q\in\Sigma$ such that $J^\infty(i)(p)=J^\infty(i)(q)$. Since $K$ is a second order, non-linear elliptic operator, and since $i$ satisfies $K_i=f\circ i$, by Aronszajn's unique continuation theorem (c.f. \cite{Aronszajn}), there exists a diffeomorphism, $\alpha$, sending a neighbourhood of $p$ to a neighbourhood of $q$ such that $\alpha(p)=q$ and $i\circ\alpha=i$. Let $\tilde{\Sigma}$ be the universal cover of $\Sigma$ and let $\tilde{\mathi}:\tilde{\Sigma}\rightarrow M$ be the lift of $i$. Applying Aronszajn's unique continuation theorem again, we extend $\alpha$ to a diffeomorphism $\tilde{\alpha}$ of $\tilde{\Sigma}$ such that $\tilde{\mathi}\circ\tilde{\alpha}=\tilde{\mathi}$ over the whole of $\tilde{\Sigma}$. However, since $i$ is prime, $\alpha\in\pi_1(\Sigma)$, so that, returning to the quotient, we have $q=\alpha(p)=p$. This is absurd, and the result follows.\qed
\medskip
A point $p\in\Sigma$ is said to be an {\bf injective point} of $i$ if and only if $i(q)\neq i(p)$ for all $q\neq p$.
\proclaim{Lemma \nextprocno}
\noindent For all $([i],f)\in\Cal{Z}$, the set injective points of $i$ is open and dense.
\endproclaim
\proclabel{PropInjectivePointsAreDense}
\proof Denote the set of injective points of $i$ by $\Omega$. For all $p\in\Sigma$ and for all $r>0$, let $B_r(p)$ be the intrinsic ball of radius $r$ about $p$ in $\Sigma$. By compactness, there exists $\epsilon>0$ such that, for all $p\in\Sigma$, the restriction of $i$ to $B_{2\epsilon}(p)$ is injective. Choose $p\in\Omega$ and denote $B:=B_\epsilon(p)$. By definition,
$$
i(p)\notin i(B^c).
$$
Since $B^c$ is compact, there exists a neighbourhood $U$ of $p$ in $B$ such that
$$
i(U)\minter i(B^c) = \emptyset.
$$
Since the restriction of $i$ to $B$ is injective, $U\subseteq\Omega$, and this proves that $\Omega$ is open.
\par
Suppose that $\Omega$ is not dense. Let $U$ be an open subset of $\Omega^c$. Choose $p\in U$ and denote $B:=B_\epsilon(p)$. Since $i$ is everywhere locally injective, the set of points distinct from $p$ but having the same image as $p$ is discrete and therefore finite. Let $q_1,...,q_n\in B^c$ be these points and define $V\subseteq\Sigma$ by
$$
V := \munion_{1\leqslant k\leqslant n}B_\epsilon(q_k).
$$
In particular,
$$
i(p) \notin i((B\munion V)^c).
$$
Since $(B\munion V)^c$ is compact, there exists a neighbourhood, $W$ of $p$ in $U$ such that
$$
i(W)\minter i((B\munion V)^c)=\emptyset.
$$
For each $k$, $i(B)\minter i(B_{2\epsilon}(q_k))$ does not contain any open subset of $i(B)$. Indeed, otherwise, $i$ would have the same $C^\infty$-jet at two distinct points, which is absurd by Lemma \procref{PropUniqueContinuation}. Thus, for each $k$, there exists a dense subset $\tilde{B}_{k}\subseteq B$ such that
$$
i(\tilde{B}_{k})\minter i(B_{2\epsilon}(q_k)) = \emptyset,
$$
and there therefore exists an open dense subset $B_k\subseteq B$ such that
$$
i(B_k)\minter i(B_\epsilon(q_k))= \emptyset.
$$
Define $B_0\subseteq B$ by
$$
B_0 := \minter_{1\leqslant k\leqslant n}B_{k}.
$$
$B_0\minter W$ is a non-trivial subset of $U$ consisting of injective points of $i$. This is absurd, and the result follows.\qed
\medskip
\noindent In particular, prime immersions in the solution space are simple.
\proclaim{Corollary \nextprocno}
\noindent For all $([i],f)\in\Cal{Z}$, for all $p\neq q\in\Sigma$, and for all sufficiently small neighbourhoods, $U_p$ and $V_q$, of $p$ and $q$ respectively
$$
i(U_p)\neq i(V_q).
$$
In particular, the immersion, $i$, is simple.
\endproclaim
\proclabel{CorEquivalenceOfNotionsOfSimplicity}
\newsubhead{Surjectivity}[Surjectivity]
We now show that the derivative, $D\hat{K}$, is surjective at every point of $\Cal{Z}$. This, together with properness, is fundamental to the development of a working topological degree theory. First, for $[i]\in\Cal{I}$, define the pairing $\langle\cdot,\cdot\rangle_i:C^\infty(\Sigma)\times C^\infty(M)\rightarrow\Bbb{R}$ by
$$
\langle g,h\rangle_i := \int_\Sigma g(x) h(i(x))\opdVol,
$$
where $\opdVol$ is the volume form of the metric induced over $\Sigma$ by the immersion $i$.
\proclaim{Lemma \nextprocno}
\noindent For $([i],f)\in\Cal{Z}$, and for any finite-dimensional subspace, $E$, of $C^\infty(\Sigma)$, there exists a finite-dimensional subspace, $F$, of $C^\infty(M)$ such that the restriction of $\langle\cdot,\cdot\rangle_i$ to $E\times F$ is non-degenerate. Furthermore, we may suppose that every element, $h$, of $F$ satisfies
$$
\langle\nabla h,N_i\rangle = 0,
$$
where $N_i$ is the unit normal vector field over $i$ compatible with the orientation.
\endproclaim
\proclabel{NonDegeneracyI}
\proof Denote $m:=\opDim(E)$. Let $\Omega\subseteq\Sigma$ be the set of injective points of $i$. By Lemma \procref{PropInjectivePointsAreDense}, $\Omega$ is open and dense, and so there exists a finite family, $(p_k)_{1\leq k\leq m}$, of points in $\Omega$ such that the map, $A:E\rightarrow\Bbb{R}^m$, given by
$$
A(g)_k := g(p_k),
$$
is a linear isomorphism. For each $k$, choose $\phi_k\in C_0^\infty(\Omega)$ and define the map $A_\phi:E\rightarrow\Bbb{R}^m$ by
$$
A_\phi(g)_k := \int g(x)\phi_k(x)\opdVol.
$$
For each $k$, let $\delta_k$ be the Dirac delta function supported on $p_k$. As $(\phi_1,...,\phi_m)$ converges to $(\delta_1,...,\delta_m)$ in the distributional sense, $A_\phi$ converges to $A$, and there therefore exists $(\phi_1,...,\phi_m)$ such that $A_\phi$ is also a linear isomorphism.
\par
Let $\pi:M\rightarrow i(\Sigma)$ be the nearest point projection. Since the restriction of $i$ to $\Omega$ is an embedding, $\pi$ is smooth near $i(\Omega)$, and so, for all $1\leq k\leq m$, there exists $h_k\in C^\infty(M)$ such that, for all $x$ near $i(\Sigma)$,
$$
h_k(x) = \phi_k\circ\pi(x).
$$
In particular,
$$\eqalign{
h_k\circ i &= \phi_k, \text{and}\cr
\langle\nabla h_k,N_i\rangle &= 0.\cr}
$$
Now let $F$ be the linear span of $h_1,...,h_m$. If $\langle g,h\rangle_i=0$ for all $h\in F$, then $A_\phi(g)=0$, and so $g=0$. On the other hand, if $h=\sum_{i=1}^m a^i h_i$ is such that $\langle g,h\rangle=0$ for all $g$, then $\sum_{i=1}^m a^i A_\phi(g)_i =0$ for all $g$, and since $A_\phi$ is a linear isomorphism, it follows that $a^i=0$ for all $i$, so that $h=0$. We conclude that the restriction of $\langle\cdot,\cdot\rangle_i$ to $E\times F$ is non-degenerate, and this completes the proof.\qed
\proclaim{Lemma \nextprocno}
\noindent $D\hat{K}$ is surjective at every point of $\Cal{Z}$.
\endproclaim
\proclabel{PropGeneralSurj}
\proof Choose $([i],f)\in\Cal{Z}$, and observe that
$$
D\hat{K}([i],f)\cdot([i,\phi],\psi) = [i,L(\phi,\psi)],
$$
where
$$
L(\phi,\psi) := \hat{J}_{i,f}\phi - \psi\circ i.
$$
It thus suffices to show that $L$ is surjective. Let $\hat{J}^*_{i,f}$ be the $L^2$-dual of $\hat{J}_{i,f}$. By elliptic regularity, $\opKer(\hat{J}^*_{i,f})$ is finite-dimensional and only consists of smooth functions. It follows by Lemma \procref{NonDegeneracyI} that there exists a finite-dimensional subspace, $E$, of $C^\infty(M)$ such that the restriction of the pairing $\langle\cdot,\cdot\rangle_i$ to $\opKer(\hat{J}_{i,f}^*)\times E$ is non-degenerate. Denote $F:=\left\{ \psi\circ i\ |\ \psi\in E\right\}$, and let $\pi:F\rightarrow\opKer(\hat{J}^*_{i,f})$ be the orthogonal projection with respect to the $L^2$ inner-product of $\Sigma$. Since $\langle\cdot,\cdot\rangle_i$ is non-degenerate, $\pi$ is a linear isomorphism, so that
$$
\opKer(\hat{J}_{i,f}^*)^\perp\minter F = \left\{0\right\}.
$$
Thus,
$$
C^\infty(\Sigma) = \opKer(\hat{J}_{i,f}^*)^\perp\oplus\opKer(\hat{J}^*_{i,f})
=\opKer(\hat{J}^*_{i,f})^\perp\oplus F
=\opIm(\hat{J}_{i,f})\oplus F,
$$
and surjectivity follows.\qed
\newsubhead{Finite dimensional sections}[FiniteDimensionalSections]
Let $X$ be a finite dimensional manifold, let $f:X\rightarrow\Cal{O}$ be a strongly smooth map, and define $\hat{K}_f:\Cal{I}\times X\rightarrow\opTI$ by
$$
\hat{K}_f([i],x) := \hat{K}([i],f_x),\eqnum{\nexteqnno[DefinitionOfKf]}
$$
define $\Cal{Z}_f\subseteq\Cal{I}\times X$ by
$$
\Cal{Z}_f := \hat{K}_f^{-1}(\left\{0\right\}),\eqnum{\nexteqnno[DefinitionOfZf]}
$$
and let $\Pi_f:\Cal{Z}_f\rightarrow X$ be the projection onto the second factor.
\par
Observe that $\hat{K}_f$ is weakly smooth over $\Cal{I}\times X$. We say that the pair $(X,f)$ is {\sl non-degenerate} whenever $D\hat{K}_f$ is surjective at every point of $\Cal{Z}_f$. In this section we show that for all such pairs, $\Cal{Z}_f$ carries the structure of a smooth, finite dimensional manifold. We first verify a technical property required of finite dimensional manifolds.
\proclaim{Lemma \nextprocno}
\noindent $\Cal{Z}_f$ is separable.
\endproclaim
\proof Since $X$ is a finite dimensional manifold, it is separable and locally compact, and the result now follows by the properness of $\Pi_f$.\qed
\medskip
The manifold structure of $\Cal{Z}_f$ is obtained by applying H\"older space theory within each graph chart. First, let $(\Phi_i,\Cal{U}_i,\Cal{V}_i)$ be a graph chart of $\Cal{I}$ and let $\hat{K}_i$ be the pull-back of $\hat{K}$ through $\Phi_i$. Define $\hat{K}_{i,f}:\Cal{U}_i\times X\rightarrow C^\infty(\Sigma)$ by
$$
\hat{K}_{i,f}(\phi,x) := \hat{K}_i(\phi,f_x),\eqnum{\nexteqnno[DefinitionOfRestrictionOfKf]}
$$
and define
$$
\Cal{Z}_{i,f} := (\hat{K}_{i,f})^{-1}(\left\{0\right\}).\eqnum{\nexteqnno[DefinitionOfRestrictionOfZf]}
$$
Observe that $\hat{K}_{i,f}$ is weakly smooth over $\Cal{U}_i\times X$ and that $\Cal{Z}_{i,f}$ is simply the pre-image of $\Cal{Z}_f$ under $\Phi_i$.
\proclaim{Lemma \nextprocno}
\noindent If $(X,f)$ is non-degenerate, then $D\hat{K}_{i,f}$ is surjective at every point of $\Cal{Z}_{i,f}$.
\endproclaim
\proclabel{SurjectivityInGraphChart}
\proof Choose $(\phi,x)\in\Cal{Z}_{i,f}$. Denote $j:=\hat{\Phi}_i(\phi)$. Let $D_1\hat{K}_{i,f}(\phi,x)$ and $D_2\hat{K}_{i,f}(\phi,x)$ denote respectively the partial derivatives of $\hat{K}_{i,f}$ at the point $(\phi,x)$ with respect to the first and second components. Likewise, let $D_1\hat{K}_f([j],x)$ and $D_2\hat{K}_f([j],x)$ denote respectively the partial derivatives of $\hat{K}_f$ at the point $([j],x)$ with respect to the first and second components. Finally, let $\hat{J}_{j,f_x}$ be the Jacobi operator of $\hat{K}$ at the point $(j,f_x)$. On the one hand (c.f. Proposition $3.8$ of \cite{MaximoNunesSmith}), there exists a smooth, positive function $g\in C^\infty(\Sigma)$ such that for all $\psi\in C^\infty(\Sigma)$,
$$
[j,D_1\hat{K}_{i,f}(\phi,x)\cdot\psi] = [j,\hat{J}_{j,f_x}(g\psi)] = D_1\hat{K}_f([j],x)\cdot[j,g\psi].
$$
On the other hand, for all $\xi_x\in T_xX$,
$$
[j,D_2\hat{K}_{i,f}(\phi,x)\cdot\xi_x] = D_2\hat{K}_f([j],x)\cdot\xi_x.
$$
It follows that $D\hat{K}_{i,f}(\phi,x)$ is surjective whenever $D\hat{K}_f([j],x)$ is, as desired.\qed
\medskip
For all $(k,\alpha)$, let $C^{k,\alpha}(\Sigma)$ be the space of $(k+\alpha)$-times H\"older differentiable functions over $\Sigma$ and define
$$
\Cal{U}^{k,\alpha}_i := \left\{\phi\in C^{k,\alpha}(\Sigma)\ |\ \|\phi\|_0<\epsilon_i\right\}.
$$
Since $\hat{K}_{i,f}$ is constructed by a finite combination of multiplication, addition, differentiation, and composition by smooth functions, it extends to a map, $\hat{K}^{2,\alpha}_{i,f}:\Cal{U}^{2,\alpha}\times X\rightarrow C^{0,\alpha}(\Sigma)$, which is smooth in the sense of maps between open subsets of Banach spaces. Furthermore, by elliptic regularity, every element of $(\hat{K}^{2,\alpha}_{i,f})^{-1}(\left\{0\right\})$ is smooth, so that
$$
\Cal{Z}_{i,f} = (\hat{K}^{2,\alpha}_{i,f})^{-1}(\left\{0\right\}).
$$
\proclaim{Lemma \nextprocno}
\noindent If $(X,f)$ is non-degenerate, then $\Cal{Z}_{i,f}$ is a smooth, embedded, $\opDim(X)$-dimensional submanifold of $\Cal{U}^{2,\alpha}\times X$. Furthermore, the canonical embedding, $e:\Cal{Z}_{i,f}\rightarrow\Cal{U}\times X$, is strongly smooth, and at every point $(\phi,x)$ of $\Cal{Z}_{i,f}$, its tangent space is given by
$$
T_{(\phi,x)}\Cal{Z}_{i,f}=\opKer(D\hat{K}_{i,f}).
$$
\endproclaim
\proof By Lemma \procref{SurjectivityInGraphChart}, $D\hat{K}_{i,f}$ is surjective at every point of $\Cal{Z}_{i,f}$. Now let $D\hat{K}^{2,\alpha}_{i,f}$ be the derivative of $\hat{K}^{2,\alpha}_{i,f}$ with respect to the H\"older norm. By elliptic regularity, $D\hat{K}^{2,\alpha}_{i,f}$ is also surjective at every point of $\Cal{Z}_{i,f}$. However, $D\hat{K}^{2,\alpha}_{i,f}$ is at every point a Fredholm map of Fredholm index $\opDim(X)$. It follows by the implicit function theorem for smooth maps between Banach manifolds that $\Cal{Z}_{i,f}$ is a smooth, embedded, $\opDim(X)$-dimensional submanifold of $\Cal{U}^{2,\alpha}\times X$. By elliptic regularity again, the canonical embedding, $e:\Cal{Z}_{i,f}\rightarrow\Cal{U}\times X$, is strongly smooth, and, furthermore, for all $(\phi,x)\in\Cal{Z}_{i,f}$,
$$
T_{(\phi,x)}\Cal{Z}_{i,f} = \opKer(D\hat{K}_{i,f}^{2,\alpha}) = \opKer(D\hat{K}_{i,f}),
$$
as desired.\qed
\proclaim{Lemma \nextprocno}
\noindent If $(X,f)$ is non-degenerate, and if $(\Phi_i,\Cal{U}_i,\Cal{V}_i)$ and $(\Phi_j,\Cal{U}_j,\Cal{V}_j)$ are graph charts of $\Cal{I}$, then the transition map $T_{ji}:=\Phi_j^{-1}\circ\Phi_i$ defines a smooth diffeomorphism from an open subset of $\Cal{Z}_{i,f}$ into an open subset of $\Cal{Z}_{j,f}$.
\endproclaim
\proof Indeed, let $Y$ be a smooth, finite dimensional manifold and let $\alpha:Y\rightarrow\Cal{Z}_{i,f}$ be a smooth map. In particular, $\alpha$ is a strongly smooth map from $Y$ into $\Cal{U}_i\times X$. Since strong smoothness is preserved by transition maps between graph charts, the composition $T_{ji}\circ\alpha$ is also a strongly smooth map from $Y$ into $\Cal{U}_j\times X$. In particular, it is smooth as a map from $Y$ into $\Cal{U}_j^{2,\alpha}\times X$. Letting $\alpha$ be the canonical embedding, we deduce that $T_{ji}$ defines a smooth map from an open subset of $\Cal{Z}_{i,f}$ into $\Cal{Z}_{j,f}$, and since its inverse is also smooth, the result follows.\qed
\medskip
\noindent Combining these results yields
\proclaim{Theorem \nextprocno}
\noindent If $(X,f)$ is non-degenerate, then $\Cal{Z}_f$ has the structure of a smooth, $\opDim(X)$-dimensional manifold. Furthermore, the canonical embedding $e:\Cal{Z}_f\rightarrow\Cal{I}\times X$ is strongly smooth, and, for all $([i],x)\in\Cal{Z}_f$, the tangent space to $\Cal{Z}_f$ at this point is given by
$$
T_{([i],x)}\Cal{Z}_f = \opKer(D\hat{K}_f).\eqnum{\nexteqnno[TangentSpaceOfZ]}
$$
\endproclaim
\proclabel{NonDegeneracyYieldsSmoothness}
\proof It only remains to determine the tangent space of $\Cal{Z}_f$. However, let $([i],x)$ be a point of $\Cal{Z}_f$. Let $(\Phi_i,\Cal{U}_i,\Cal{V}_i)$ be the graph chart of $\Cal{I}$ about $i$. Then,
$$
T_{([i],x)}\Cal{Z}_f
=(\Phi_i)_*T_{(0,x)}\Cal{Z}_{i,f}
=(\Phi_i)_*\opKer(\hat{K}_{i,f}(0,x))
=\opKer(\hat{K}_f([i],x)),
$$
as desired.\qed
\newsubhead{Extensions}[Extensions]
Let $Y$ be a finite dimensional subspace of $C^\infty(M)$ and for all $\delta>0$, let $Y_\delta$ be the closed ball of radius $\delta$ about $0$ in $Y$ with respect to the $C^0$-norm. For $\delta$ sufficiently small, define the strongly smooth map, $\tilde{f}:X\times Y_\delta\rightarrow\Cal{O}$, by
$$
\tilde{f}_{x,y} = f_x + y.
$$
We call $(X\times Y_\delta,\tilde{f}_{x,y})$ the {\sl extension} of $(X,f)$ along $Y_\delta$.
\proclaim{Lemma \nextprocno}
\noindent If $X$ is compact, then there exists a finite-dimensional subspace $Y\subseteq C^\infty(M)$ such that for all sufficiently small $\delta$, the extension of $(X,f)$ along $Y_\delta$ is non-degenerate.
\endproclaim
\proclabel{Surjectivity}
\proof Choose $p:=([i],x)\in\Cal{Z}_f$. By ellipticity, the image of $D\hat{K}_f([i],x)$ has finite codimension in $\opTIi$. By Lemma \procref{PropGeneralSurj}, there exists a finite dimensional subspace, $Y_p\subseteq C^\infty(M)$, such that if $Y$ contains $Y_p$, then the derivative, $D\hat{K}_{\tilde{f}}([i],x,0)$, of the extension is surjective. Since surjectivity of elliptic operators is an open property, and since, by properness, $\Cal{Z}_f$ is compact, there exist $p_1,...,p_m\in\Cal{Z}_f$ such that if $Y:=Y_{p_1}+...+Y_{p_m}$, then $D\hat{K}_{\tilde{f}}([i],x,0)$ is surjective for all $([i],x)\in\Cal{Z}_f$. Finally, by properness, $\Cal{Z}_{\tilde{f}}$ converges to $\Cal{Z}_f$ in the Hausdorff sense in $\Cal{I}\times X\times Y$ as $\delta$ tends to $0$. Thus, using again the fact that surjectivity of elliptic operators is an open property, we find that, for all sufficiently small $\delta$, $D\hat{K}_{\tilde{f}}$ is surjective at every point of $\Cal{Z}_{\tilde{f}}$, as desired.\qed
\medskip
We now modify $(X,f)$ in such a manner as to make it non-degenerate. To this end, if $(X\times Y_\delta,\tilde{f})$ is an extension of $(X,f)$, then, for a smooth function, $g:X\rightarrow Y_\delta$, we denote by $X_g$ its graph in $X\times Y_\delta$, and by $f_g$ the restriction of $\tilde{f}$ to this graph. We recall that a subset of any complete metric space is said to be in the second category in the sense of Baire whenever it is the intersection of a countable family of dense, open sets. A given property is then said to hold {\sl generically} whenever it holds for all elements of some such subset.
\proclaim{Theorem \nextprocno}
\noindent Let $(X\times Y_\delta,\tilde{f})$ be a non-degenerate extension of $(X,f)$. Let $X_0\subseteq X$ be the set of all points $x$ such that $f_x$ is a regular value of $\Pi$. Let $\Cal{G}_0$ be the set of all smooth functions $g:X\rightarrow Y_\delta$ which vanish over $X_0$. Then, for generic $g\in\Cal{G}_0$, $(X_g,f_g)$ is also non-degenerate.
\endproclaim
\proclabel{Transversality}
\proof Let $\Pi_{\tilde{f}}:\Cal{Z}_{\tilde{f}}\rightarrow X\times Y_\delta$ be the canonical projection. By classical transversality theory (c.f. \cite{GuillemanPollack}), it suffices to show that $X_g$ has the desired properties whenever it is transverse to $\Pi_{\tilde{f}}$. Suppose therefore that $X_g$ is transverse to $\Pi_{\tilde{f}}$. Let $([i],x,y)$ be a point of $\Cal{Z}_{\tilde{f}}$ above some point, $(x,y)$, of $X_g$. We need to show that the restriction of $D\hat{K}_{\tilde{f}}([i],x,y)$ maps $\opTIi\times T_{(x,y)}X_g$ surjectively onto $\opTIi$. Thus, let $\nu$ be an element of $\opTIi$. Since $(X\times Y_\delta,\tilde{f})$ is non-degenerate, there exists $(\psi_1,\xi_1,\eta_1)\in\opTIi\times T_xX\times T_yY$ such that
$$
D\hat{K}_{\tilde{f}}([i],x,y)\cdot(\psi_1,\xi_1,\eta_1) = \nu.
$$
However, since $X_g$ is transverse to $\Pi_{\tilde{f}}$ at $(x,y)$, there exists $(\psi_2,\xi_2,\eta_2)\in T_{[i],x,y}\Cal{Z}_{\tilde{f}}$ and $(\xi_3,\eta_3)\in T_{(x,y)}X_g$ such that
$$
(\xi_1,\eta_1) = (\xi_2,\eta_2) + (\xi_3,\eta_3).
$$
Thus, since $T_{[i],x,y}\Cal{Z}_{\tilde{f}}=\opKer(D\hat{K}_{\tilde{f}}([i],x,y))$, setting $\psi_3:=\psi_1-\psi_2$, yields
$$\eqalign{
D\hat{K}_{\tilde{f}}([i],x,y)(\psi_3,\xi_3,\eta_3)
&=D\hat{K}_{\tilde{f}}([i],x,y)(\psi_1,\xi_1,\eta_1) - D\hat{K}_{\tilde{f}}([i],x,y)(\psi_2,\xi_2,\eta_2),\cr
&=D\hat{K}_{\tilde{f}}([i],x,y)(\psi_1,\xi_1,\eta_1),\cr
&=\nu,\cr}
$$
and the result follows.\qed
\newsubhead{Orientation - the finite-dimensional case}[DefiningTheSignatureFD]
Let $E$ and $F$ be oriented vector spaces of finite dimensions equal to $(m+n)$ and $n$ respectively. Let $\opGr_m(E)$ be the grassmannian of $m$-dimensional linear subspaces of $E$, let $\opGr_m^+(E)$ be the grasssmannian of oriented, $m$-dimensional linear subspaces of $E$, and let $\pi:\opGr_m^+(E)\rightarrow\opGr_m(E)$ be the canonical projection. Elements of $\opGr_m^+(E)$ are given by pairs of the form $(K,[\mu])$, where $K$ is an $m$-dimensional linear subspace of $E$ and $[\mu]$ is an $\Bbb{R}^+$-equivalence class of volume forms over $K$. In what follows, we shall make no distinction between a volume form, $\mu$, and its equivalence class, $[\mu]$.
\par
Let $\Cal{M}:=\Cal{M}(E,F)$ be the space of all surjective linear maps from $E$ into $F$, and let $K:\Cal{M}\rightarrow\opGr_m(E)$ be the function which sends a linear map to its kernel.
\proclaim{Lemma \nextprocno}
\noindent There exists a canonical continuous map, $K^+:\Cal{M}\rightarrow\opGr_m^+(E)$, such that $\pi\circ K^+=K$.
\endproclaim
\proclabel{CanonicalLifting}
\remark We refer to $K^+$ as the {\sl canonical} lifting of $K$. Observe that, via reversal of orientation, the existence of one lifting implies the existence of exactly one other. However, $K^+$ is uniquely defined by the orientations of $E$ and $F$. Furthermore, in the case that will be of interest to us, where $E=X\oplus F$, since $F$ is also a summand of the domain, it is sufficient to prescribe an orientation only on $X$ for $K^+$ to be uniquely defined.
\medskip
\proof Choose $A\in\Cal{M}$. Let $L$ be any complementary subspace of $K:=K(A)$ in $E$. Let $\mu_E$ and $\mu_F$ be the orientation forms of $E$ and $F$. There exists a unique orientation form $\mu_K$ over $K$ such that
$$
\mu_E = p^*\mu_K\wedge q^*(A|_L)^*\mu_F,\eqnum{\nexteqnno[FormulaForOrientationForm]}
$$
where $p:E\rightarrow K$ and $q:E\rightarrow L$ are the projections along $L$ and $K$ respectively. Furthermore, up to a positive scalar factor, $\mu_K$ is well-defined independent of $L$. The result follows.\qed
\medskip
Now let $K_0$ and $L_0$ be preferred complementary subspaces of $E$ of dimension $m$ and $n$ respectively, furnished with the respective orientation forms $\mu_{K_0}$ and $\mu_{L_0}$. Let $p_0:E\rightarrow K_0$ and $q_0:E\rightarrow L_0$ be the projections along $L_0$ and $K_0$ respectively.
\proclaim{Lemma \nextprocno}
\noindent For all $A\in\Cal{M}$, the restriction of $p_0$ to $K(A)$ is a linear isomorphism if and only if $A|_{L_0}$ is a linear isomorphism.
\endproclaim
\proof Indeed,
$$
\opKer(p_0|_{K(A)}) = \opKer(p_0)\minter\opKer(A) = \opKer(A|_{L_0}),
$$
and the result follows.\qed
\medskip
\noindent Suppose now that the orientation form, $\mu_E$, of $E$ satisfies
$$
\mu_E = p_0^*\mu_{K_0}\wedge q_0^*\mu_{L_0}.
$$
We then obtain the following explicit formula for $K^+$.
\proclaim{Lemma \nextprocno}
\noindent For all $A$ such that $A|_{L_0}$ is invertible,
$$
K^+(A) = (K(A),\opDet(A|_{L_0})^{-1}p_0^*[\mu_{K_0}]).\eqnum{\nexteqnno[FormulaForOrientation]}
$$
\endproclaim
\proof Indeed, when $A|_{L_0}$ is invertible, $L_0$ is complementary to $K(A)$ and $\mu_{L_0}$ is equal to $\opDet(A|_{L_0})^{-1}(A|_{L_0})^*\mu_F$. Furthermore, if $p:E\rightarrow K(A)$ is the projection along $L_0$, then $p_0p=p_0$, so that $p^*p_0^*\mu_{K_0}=p_0^*\mu_{K_0}$. Combining these results yields
$$
\mu_E = p_0^*\mu_{K_0}\wedge q_0^*\mu_{L_0} = p^*(\opDet(A|_{L_0})^{-1}p_0^*\mu_{K_0})\wedge(A|_{L_0})^*\mu_F,
$$
and the result follows by \eqnref{FormulaForOrientationForm}.\qed
\medskip
It is important to know how extensions affect $K^+(A)$. Thus, let $V$ be a finite-dimensional vector space. Denote $\tilde{E}:=E\oplus V$, $\tilde{F}:=F\oplus V$ and define $\tilde{\Cal{M}}:=\Cal{M}(\tilde{E},\tilde{F})$, $\tilde{K}$ and $\tilde{K}^+$ as before. Let $\pi_1:\tilde{E}\rightarrow E$, $\pi_2:\tilde{E}\rightarrow V$, $\tau_1:\tilde{F}\rightarrow F$ and $\tau_2:\tilde{F}\rightarrow V$ be the canonical projections.
\proclaim{Lemma \nextprocno}
\noindent For $\tilde{A}\in\tilde{\Cal{M}}$, if $\tau_1\circ\tilde{A}|_V=0$ and if $\tau_2\circ\tilde{A}|_V$ is a linear isomorphism, then $A:=\tau_1\circ\tilde{A}|_E$ is surjective and $\pi_1$ restricts to a linear isomorphism from $\tilde{K}(\tilde{A})$ into $K(A)$.
\endproclaim
\proclabel{Extensions}
\proof First observe that
$$
\opKer(\pi_1|_{\tilde{K}(\tilde{A})}) = \opKer(\pi_1)\minter\opKer(\tilde{A}) = \opKer(\tilde{A}|_V),
$$
so that the restriction of $\pi_1$ to $\tilde{K}(\tilde{A})$ is injective. Now choose $\tilde{\xi}\in\tilde{K}(\tilde{A})$ and write $\tilde{\xi}=:(\xi,\eta)$ where $\xi\in E$ and $\eta\in V$. Then
$$
A\xi = (\tau_1\circ\tilde{A})(\xi,0) = (\tau_1\circ\tilde{A})(\xi,\eta) = 0,
$$
so that $\xi\in K(A)$, and the result follows since $\tilde{K}(\tilde{A})$ and $K(A)$ both have the same dimension.\qed
\medskip
\noindent In particular, thinking of $K^+$ as a section of $\opGr_m^+(E)$ over $K$, near any $\tilde{A}\in\tilde{\Cal{M}}$ satisfying the hypotheses of Lemma \procref{Extensions}, the pull-back, $\pi^*K^+$, of this section through $\pi$ is well-defined, and applying \eqnref{FormulaForOrientation} immediately yields
\proclaim{Corollary \nextprocno}
\noindent For $\tilde{A}\in\tilde{\Cal{M}}$, if $\tau_1\circ\tilde{A}|_V=0$ and if $\tau_2\circ\tilde{A}|_V$ is a linear isomorphism, then, denoting $A:=\tau_1\circ\tilde{A}|_E$,
$$
\tilde{K}^+(\tilde{A}) = \opSig(\tau_2\circ\tilde{A}|_V)\pi_1^* K^+(A),
$$
where, for any matrix, $M$, $\opSig(M)$ is the sign of its determinant.
\endproclaim
\proclabel{CorExtensions}
\newsubhead{Orientation - the infinite-dimensional case}[DefiningTheSignatureID]
Let $E$ and $F$ be Hilbert spaces, and let $i:E\rightarrow F$ be a compact injection with dense image. Via this injection, any Fredholm map, $B:E\rightarrow F$, can be understood as a closed map from $F$ to itself with compact resolvent (c.f. \cite{Kato}). Now let $\Cal{A}$ be the space of all Fredholm maps, $B:E\rightarrow F$, such that
$$
\minf_{v\neq 0}\frac{\langle B(v),i(v)\rangle}{\langle i(v),i(v)\rangle} > -\infty.\eqnum{\nexteqnno[RealEigenvaluesBoundedBelow]}
$$
Since this is a convex condition, $\Cal{A}$ is an open subset of $\opLin(E,F)$. Let $A$ be an element of $\Cal{A}$. Given an eigenvalue, $\lambda$, of $A$, we denote by $E_\lambda$ its null-space in $E$, that is, the union of the kernels of $(A-\lambda)^k$ as $k$ varies over all non-negative integers, and we define its {\bf multiplicity} to be the dimension of this space. Since $A$ has compact resolvent, its spectrum is discrete and consists only of eigenvalues all of which have finite multiplicity (c.f. \cite{Kato}). Furthermore, it follows from \eqnref{RealEigenvaluesBoundedBelow} that $A$ has only finitely many strictly negative real eigenvalues. The {\bf index} of $A$, which we denote by $\opInd(A)$, is then defined to be the sum of their multiplicities, and its {\bf signature}, which we denote by $\opSig(A)$, is defined by
$$
\opSig(A) := (-1)^{\opInd(A)}.\eqnum{\nexteqnno[DefinitionOfSignature]}
$$
Observe that, although the index is not stable under small perturbations away from a given invertible map, $A$, the signature is. Indeed, although strictly negative real eigenvalues may perturb to complex eigenvalues, they only do so in conjugate pairs.
\par
Now let $X$ be an oriented, finite-dimensional vector space and let $p:X\oplus E\rightarrow X$ be the canonical projection. Let $\Cal{M}$ be the space of all pairs, $(M,A)$, such that
\medskip
\myitem{(1)} $M:X\rightarrow F$ is linear;
\medskip
\myitem{(2)} $A:E\rightarrow F$ is an element of $\Cal{A}$; and
\medskip
\myitem{(3)} $M+A$ is surjective.
\medskip
\noindent Let $\opGr(X,E)$ be the grassmannian of $\opDim(X)$-dimensional linear subspaces of $X\oplus E$, let $\opGr^+(X,E)$ be the grassmannian of oriented $\opDim(X)$-dimensional linear subspaces of $X\oplus E$, and let $\pi:\opGr^+(X,E)\rightarrow\opGr(X,E)$ be the canonical projection. Let $K:\Cal{M}\rightarrow\opGr(X,E)$ be the map which sends the pair $(M,A)$ to the kernel of $M+A$.
\proclaim{Theorem \nextprocno}
\noindent There exists a canonical continuous map $K^+:\Cal{M}\rightarrow\opGr^+(X,E)$ such that $\pi\circ K^+=K$. Furthermore, if $A$ is invertible, then $p$ restricts to a linear isomorphism from $K(M,A)$ onto $X$, and
$$
K^+(M,A) = (K(M,A),\opSig(A) p^*[\mu_X]),\eqnum{\nexteqnno[DefinitionOfParity]}
$$
where $[\mu_X]$ is the orientation form of $X$.
\endproclaim
\proclabel{LemmaDefinitionOfParity}
\proof Let $\Lambda$ be a finite subset of $\opSpec(A)$ which contains all non-positive real eigenvalues and which is symmetric under complex conjugation. Denote
$$
E_\Lambda := \oplus_{\lambda\in\Lambda}E_\lambda.
$$
Observe that $A$ preserves $E_\Lambda$. Furthermore, by classical spectral theory (c.f. \cite{Kato}), there exists a closed subspace $R\subseteq F$ such that $E=E_\Lambda\oplus(R\minter E)$, $F=E_\Lambda\oplus R$ and $A$ maps $R\minter E$ into $R$.
\par
Consider first $N:X\rightarrow E_\Lambda$ and $B:E_\Lambda\rightarrow E_\Lambda$ such that $N+B$ is surjective and denote by $K(N,B)$ the kernel of $N+B$ in $X\oplus E_\Lambda$. Importantly, since the codomain, $E_\Lambda$, is also a summand of the domain, it is not necessary to prescribe an orientation over this space in order to uniquely define the canonical lifting constructed in Lemma \procref{CanonicalLifting} (c.f. the remark following that result), and we denote this lifting by $K^+(N,B)$.
\par
Now let $P:X\times E\rightarrow X\times E_\Lambda$ and $Q:F\rightarrow E_\Lambda$ be the projections along $R\minter E$ and $R$ respectively. Since $E_\Lambda$ contains the kernel of $A$, it follows as in Lemma \procref{Extensions} that $P$ restricts to a linear isomorphism from $K(M,A)$ into $K(Q\circ M,Q\circ A|_{E_\Lambda})$, and we therefore define
$$
K^+_\Lambda(M,A) := P^*K^+(Q\circ M,Q\circ A|_{E_\Lambda}).
$$
\par
We claim that $K^+_\Lambda(M,A)$ is independent of the choice of $\Lambda$. To show this, observe first that if $\Lambda'$ is another finite subset of $\opSpec(A)$ which contains all non-positive real eigenvalues and which is symmetric under complex conjugation, then so too is $\Lambda\minter\Lambda'$, so that it suffices to consider the case where $\Lambda\subseteq\Lambda'$. Now define $P':X\times E\rightarrow X\times E_{\Lambda'}$ and $Q':F\rightarrow E_{\Lambda'}$ as before. Denote $\Lambda'':=\Lambda'\setminus\Lambda$ and let $\pi_1:X\times E_{\Lambda'}\rightarrow X\times E_\Lambda$, $\pi_2:X\times E_{\Lambda'}\rightarrow E_{\Lambda''}$, $\tau_1:E_{\Lambda'}\rightarrow E_\Lambda$ and $\tau_2:E_{\Lambda'}\rightarrow E_{\Lambda''}$ be the canonical projections. Since $\tau_1\circ A|_{E_{\Lambda''}}=0$, and since $\tau_2\circ A|_{E_{\Lambda''}}$ is an orientation-preserving linear isomorphism, it follows by Lemma \procref{Extensions} and Corollary \procref{CorExtensions} that $\pi_1$ restricts to a linear isomorphism from $K(Q'\circ M,Q'\circ A|_{E_{\Lambda'}})$ into $K(Q\circ M, Q\circ A|_{E_{\Lambda}})$ and
$$
K^+(Q'\circ M,Q'\circ A|_{E_{\Lambda'}}) = \pi_1^*K^+(Q\circ M,Q\circ A|_{E_{\Lambda}}).
$$
Thus, since $P=\pi_1\circ P'$,
$$\eqalign{
K_{\Lambda'}^+(M,A) &= (P')^*K^+(Q'\circ M,Q'\circ A|_{E_{\Lambda'}})\cr
&=(P')^*\pi_1^* K^+(Q\circ M,Q\circ A|_{E_{\Lambda}})\cr
&=P^*K^+(Q\circ M,Q\circ A|_{E_{\Lambda}})\cr
&=K^+_\Lambda(M,A),\cr}
$$
so that $K^+_\Lambda(M,A)$ is indeed independent of $\Lambda$. We now define $K^+(M,A):=K_\Lambda^+(M,A)$, for any suitable $\Lambda$. It follows from the continuous dependence on $A$ of its eigenspace decomposition (c.f. \cite{Kato}) that $K^+(M,A)$ is continuous. Finally, when $A$ is invertible, so too is $Q\circ A|_{E_\Lambda}$, and so, by \eqnref{FormulaForOrientation},
$$\eqalign{
K^+(M,A) &= P^*K^+(Q\circ M,Q\circ A|_{E_{\Lambda}}),\cr
&= P^*(K(Q\circ M,Q\circ A|_{E_{\Lambda}}),\opDet(Q\circ A|_{E_{\Lambda}})^{-1}p^*[\mu_X]),\cr
&= (K(Q,A),\opSig(A)P^*p^*[\mu_X]),\cr
&= (K(Q,A),\opSig(A)p^*[\mu_X]),\cr}
$$
as desired.\qed
\newsubhead{Constructing the degree}[ConstructingTheDegree]
We first control the spectrum of $\hat{J}_{i,f}$ in order to ensure that the orientation is well-defined.
\proclaim{Lemma \nextprocno}
\noindent If $L:f\mapsto -a^{ij}f_{;ij} + b^if_{;i} + cf$ is a generalised Laplacian over $\Sigma$, then there exists $B>0$, which only depends on the metric on $\Sigma$ as well as
\medskip
\myitem{(1)} the $C^1$ norm of $a$; and
\medskip
\myitem{(2)} the $C^0$ norms of $a^{-1}$, $b$ and $c$;
\medskip
\noindent such that the real eigenvalues of $L$ lie in $]-B,+\infty[$.
\endproclaim
\proclabel{LemmaLowerBoundOfRealPart}
\proof At each point $p\in\Sigma$, consider $a^{ij}$ as a scalar product over $T^*_p\Sigma$. Since $L$ is elliptic, this yields a metric, $a_{ij}$, over $\Sigma$. Let $\Gamma^k_{ij}$ be the relative Christophel symbol of the Levi-Civita covariant derivative of this metric with respect to that of the standard metric. Thus, if ``,'' denotes covariant differentiation with respect to the Levi-Civita covariant derivative of $a_{ij}$, then, for all $f$
$$
f_{,ij} = f_{;ij} - \Gamma^k_{ij}f_{;k}.
$$
Observe that $\Gamma$ depends on the first derivative of $a$. Denote
$$
\tilde{b}^i = b^i - {\Gamma^i}_{pq}a^{pq},
$$
so that, for all $f$
$$\matrix
Lf \hfill&= -a^{ij}f_{,ij} + \tilde{b}^if_{,i} + cf\hfill\cr
&= -\Delta^af + \tilde{b}^if_{,i} + cf,\hfill\cr
\endmatrix$$
where $\Delta^a$ is the Laplacian of the metric $a_{ij}$.
\par
Now let $\lambda$ be a real eigenvalue of $L$ and let $f\in C^\infty(\Sigma,\Bbb{R})$ be a corresponding real eigenvector, chosen such that
$$
\|f\|_{L^2}^2 = \int_\Sigma f^2\opdVol^a = 1,
$$
where $\opdVol^a$ here denotes the volume form of the metric $a$. Bearing in mind Stokes' Theorem and the Cauchy-Schwarz Inequality, we obtain
$$\matrix
\lambda\hfill&= \int_\Sigma fLf\opdVol^a\hfill\cr
&=\int_\Sigma -f\Delta^a f + f\tilde{b}^if_{,i} + cf^2\opdVol^a\hfill\cr
&=\int_\Sigma \|\nabla f\|_a^2 + f\tilde{b}^if_{,i} + cf^2\opdVol^a\hfill\cr
&\geqslant \|\nabla f\|_{L^2}^2 - \|\tilde{b}\|_{L^\infty}\|\nabla f\|_{L^2} - \|c\|_{L^\infty}\hfill\cr
&=(\|\nabla f\|_{L^2} - (1/2)\|\tilde{b}\|_{L^\infty})^2 - \|c\|_{L^\infty} - \frac{1}{4}\|\tilde{b}\|_{L^\infty}^2\hfill\cr
&\geqslant  - \|c\|_{L^\infty} - \frac{1}{4}\|\tilde{b}\|_{L^\infty}^2,\hfill\cr
\endmatrix$$
as desired.\qed
\proclaim{Theorem \nextprocno}
\noindent If $(X,f)$ is non-degenerate, then $\Cal{Z}_f$ carries a canonical orientation, $[\mu_Z]$, such that, for $([i],x)\in\Cal{Z}_f$, if $J_{i,x}$ is non-degenerate, then
$$
[\mu_Z] = \opSig(J_{i,x})\Pi_f^*[\mu_X],\eqnum{\nexteqnno[DefinitionOfOrientationForm]}
$$
where $[\mu_X]$ here denotes the orientation form of $X$.
\endproclaim
\proclabel{Orientation}
\proof Choose $p:=([i],x)\in\Cal{Z}_f$. Since $(X,f)$ is non-degenerate, $D\hat{K}_f([i],x)$ is surjective. However, let $D_1\hat{K}_f$ denote the partial derivative of $\hat{K}_f$ with respect to the first component. By \eqnref{CovariantDerivativeOfKhati}, for all $g$,
$$
D_1\hat{K}_f([i],x)[i,g] = [i,\hat{J}_{i,f_x}g],
$$
where $\hat{J}_{i,f_x}$ is the Jacobi operator of $\hat{K}$ at the point $(i,f_x)$. However, by ellipticity, $\hat{J}_{i,f_x}$ defines a Fredholm map of index $0$ from the Sobolev space $H^2(\Sigma)$ into the Hilbert space $L^2(\Sigma)$. Furthermore, by Lemma \procref{LemmaLowerBoundOfRealPart}, the real part of its spectrum is bounded below. It thus follows by Theorem \procref{LemmaDefinitionOfParity} that there exists a canonical orientation $[\mu_Z]$ over $T_{([i],x)}\Cal{Z}_f=\opKer(D\hat{K}_f([i],x))$ which satisfies \eqnref{DefinitionOfOrientationForm} whenever $\hat{J}_{i,f_x}$ is non-degenerate. Finally, since the canonical embedding, $e:\Cal{Z}_f\rightarrow\Cal{I}\times X$, is strongly smooth, $[\mu_Z]$ varies continuously with $p\in\Cal{Z}_f$, and this completes the proof.\qed
\proclaim{Theorem \nextprocno, {\bf Existence of a degree}}
\noindent The set, $\Cal{O}'$, of regular values of $\Pi$ is open and dense in $\Cal{O}$. Furthermore, if for all $f\in\Cal{O}'$, we define
$$
\opDeg(\Pi;f) := \sum_{([i],f)\in\Pi^{-1}(\left\{f\right\})}\opSig(\hat{J}_{i,f}),
$$
then $\opDeg(\Pi;f)$ is independent of the point $f\in\Cal{O}'$ chosen.
\endproclaim
\proclabel{DefinitionOfTheDegree}
\proof The openness of $\Cal{O}'$ follows since $\Pi$ is proper and since surjectivity of elliptic maps is an open property. To see that it is dense, consider a function, $f$, in $\Cal{O}$. Let $X:=\left\{x\right\}$ be the manifold consisting of a single point and think of $f$ as a strongly smooth map from $X$ into $\Cal{O}$. By Lemma \procref{Surjectivity}, there exists an extension, $(X\times Y_\delta,\tilde{f})$, which is non-degenerate and by Theorem \procref{Transversality}, for generic $y\in Y_\delta$, $\tilde{f}(y)$ is a regular value of $\Pi$. Since $f\in\Cal{O}$ is arbitrary, it follows that $\Cal{O}'$ is dense, as desired.
\par
Now consider two functions, $f_0$ and $f_1$, in $\Cal{O}'$. Let $f:[0,1]\rightarrow\Cal{O}$ be a strongly smooth curve such that $f(0)=f_0$ and $f(1)=f_1$. By Lemma \procref{Surjectivity}, again, there exists an extension, $([0,1]\times Y_\delta,\tilde{f})$, which is non-degenerate. Denote $\tilde{X}:=[0,1]\times Y_\delta$ and choose an orientation over this space. By Theorem \procref{Orientation}, $\Cal{Z}_{\tilde{f}}$ also carries a canonical orientation, so that $\Pi_{\tilde{f}}:\Cal{Z}_{\tilde{f}}\rightarrow X$ has a well defined $\Bbb{Z}$-valued mapping degree. However, it follows from the definition that this mapping degree is equal to $\opDeg(\Pi;f_0)$ and $\opDeg(\Pi;f_1)$, so that the two are equal, as desired.\qed
\newsubhead{Varying the metric}[VaryingTheMetric]
We conclude this section by outlining how to adapt the results of the preceding sections to allow for variations of the metric. Let $\Cal{G}$ be an open subset of the space of Riemannian metrics over $M$, let $\hat{\Cal{I}}$ be an open subset of $\Cal{I}\times\Cal{G}$, and for all $g\in\Cal{G}$, let $\Cal{I}_g$ be the fibre of $\hat{\Cal{I}}$ over $g$, that is
$$
\Cal{I}_g = \left\{ [i]\in\Cal{I}\ |\ ([i],g)\in\hat{\Cal{I}}\right\}.
$$
Let $K$ be an elliptic curvature function, and suppose that for all $([i],g)\in\hat{\Cal{I}}$, the immersion, $i$, is $K$-convex with respect to the metric $g$.  Observe, in particular, that, in this case, the fibre will vary with $g$. For all $([i],g)\in\hat{\Cal{I}}$, the {\bf $K$-curvature} of $i$ with respect to $g$ is then defined by
$$
K_{i,g}(p) := K(A_{i,g}(p)),\eqnum{\nexteqnno[DefinitionOfKCurvatureWithVaryingMetric]}
$$
where $A_{i,g}$ is the shape operator of $i$ with respect to $g$. This defines a weakly smooth section of $\Cal{T}\hat{\Cal{I}}$ over $\hat{\Cal{I}}$ which we also denote by $K$, that is, for all $([i],g)\in\hat{\Cal{I}}$,
$$
K([i],g) = [i,K_{i,g}].\eqnum{\nexteqnno[DefinitionOfCurvatureSectionWithVaryingMetric]}
$$
Let $\hat{\Cal{O}}$ be an open subset of $C^\infty(M)\times\Cal{G}$ and for all $g\in\Cal{G}$, let $\Cal{O}_g$ be its fibre over $g$. Since the conditions required on $\Cal{O}_g$ in order to prove properness often depend on $g$, this fibre will also often vary with $g$.
\par
Define
$$
\hat{\Cal{U}} := \left\{([i],f,g)\ |\ ([i],g)\in\hat{\Cal{I}}\ \&\ (f,g)\in\hat{\Cal{O}}\right\},
$$
and for all $g\in\Cal{G}$, let $\Cal{U}_g$ be the fibre of $\Cal{U}$ over $g$, that is
$$
\Cal{U}_g = \Cal{I}_g\times\Cal{O}_g.
$$
Define the evaluation functional, $E:\hat{\Cal{U}}\rightarrow\Cal{T}\hat{\Cal{I}}$, by
$$
E([i],f,g) = [i,f\circ i],\eqnum{\nexteqnno[DefinitionOfEvalutionSectionWithVaryingMetric]}
$$
and define the weakly smooth section, $\hat{K}:\hat{\Cal{U}}\rightarrow\Cal{T}\hat{\Cal{J}}$, by
$$
\hat{K}([i],f,g) = K([i],g) - E([i],f).
$$
The {\bf solution space} $\hat{\Cal{Z}}\subseteq\hat{\Cal{U}}$ is now defined by
$$
\Cal{Z} = \hat{K}^{-1}(\left\{0\right\}),
$$
and $\Pi:\hat{\Cal{Z}}\rightarrow\hat{\Cal{O}}$ is defined to be the projection onto the second and third factors. As before, we will suppose
\medskip
{\bf\noindent Properness:\ }$\Pi$ defines a proper map from $\hat{\Cal{Z}}$ into $\hat{\Cal{O}}$.
\medskip
We leave the reader to verify that the results developed in this paper readily extend to this framework. In particular, it is worth noting how the genericity in our results depends on the surjectivity properties of $D\hat{K}$. Indeed, Theorem \procref{DefinitionOfTheDegree} is valid for generic prescribing functions given a fixed metric. This follows from the fact, proven in Lemma \procref{PropGeneralSurj}, that $D_1\hat{K}\oplus D_2\hat{K}$ is surjective at every point of $\Cal{Z}$. In particular, by showing that $D_1\hat{K}\oplus D_3\hat{K}$ is surjective at every point of $\Cal{Z}$, we show that Theorem \procref{DefinitionOfTheDegree} is also valid for generic metrics given a fixed prescribing function. We proceed as follows.
\proclaim{Lemma \nextprocno}
\noindent Let $([i],f,g)$ be a point of $\Cal{U}$. If $\phi\in C^\infty(M)$ is such that $\langle\nabla\phi,N_i\rangle=0$ over $i$, where $N_i$ is the unit normal vector field over $i$ that is compatible with the orientation. Then,
$$
D_3\hat{K}([i],f,g)\cdot(\phi g) = [i,-\phi K_i/2].
$$
\endproclaim
\proclabel{PreliminaryToSurjectivityForMetrics}
\proof Indeed, consider the family, $g^t$, of metrics given by $g^t:=e^{2 t\phi}g$. Let $\overline{\nabla}$ be the Levi-Civita covariant derivative of $g$, and for all $t$, let $\overline{\nabla}^t$ be the Levi-Civita covariant derivative of $g^t$. For all $t$, let $\Omega^t:=\overline{\nabla}^t-\overline{\nabla}$ be the relative Christoffel symbol of $\overline{\nabla}^t$ with respect to $\overline{\nabla}$. Using the subscript ``$:$'' to denote covariant differentiation of $\overline{\nabla}$ and raising and lowering indices with respect to $g$, by the Koszul formula, we have
$$
(\Omega^t){}^k_{ij} = \delta^k_i\phi_{:j} + \delta^k_j\phi_{:i} - g_{ij}\phi_{:}{}^k.
$$
Let $N_i$ be the unit normal vector field of $i$ with respect to $g$, which is compatible with the orientation, and observe that, for all $t$, $N_i^t:=e^{-t\phi}N_i$ is its unit normal vector field with respect to the metric $g^t$. Let $A_i$ be the shape operator of $i$ with respect to $g$, and for all $t$ let $A_i^t$ be its shape operator with respect to $g^t$. For all $X$ tangent to $i$, bearing in mind the hypotheses on $\phi$,
$$\eqalign{
A_i^t\cdot X &= \nabla_X^tN_i^t\cr
&= \nabla_Xe^{-t\phi}N_i + e^{-t\phi}\Omega^t(X,N_i)\cr
&= e^{-t\phi}\nabla_X N_i + te^{-t\phi}\langle\nabla\phi,N_i\rangle X\cr
&= e^{-t\phi}A_i\cdot X.\cr}
$$
Consequently
$$
\partial_t A^t_i|_{t=0} = -\phi A_i,
$$
so that
$$
\partial_t K(A^t_i)|_{t=0} = DK(A_i)\cdot\partial_t A^t_i|_{t=0} = -\phi DK(A_i)\cdot A_i.
$$
The result now follows, since $K$ is homogeneous of order $1$.\qed
\proclaim{Lemma \nextprocno}
\noindent $D_1\hat{K}\oplus D_3\hat{K}$ is surjective at every point of $\Cal{Z}$.
\endproclaim
\proclabel{PropGeneralSurjMetrics}
\proof Choose $([i],f,g)\in\Cal{Z}$ and let $N_i$ be the unit normal vector field over $i$ with respect to $g$ which is compatible with the orientation. If $\phi\in C^\infty(M)$, and if $\psi\in C^\infty(M)$ satisfies $\langle\nabla\phi,N_i\rangle=0$, then, by Lemma \procref{PreliminaryToSurjectivityForMetrics},
$$
D_1\hat{K}([i],f,g)\cdot[i,\phi] + D_3\hat{K}([i],f,g)\cdot(\psi g) = [i, L(\phi,\psi)],
$$
where
$$
L(\phi,\psi) = \hat{J}_{i,f,g}\phi - \frac{1}{2}\psi K_{i,g},
$$
and the result now follows as in Lemma \procref{PropGeneralSurj}.\qed
\medskip
\noindent It thus follows that Theorem \procref{DefinitionOfTheDegree} is valid for generic metrics given any fixed prescribing function, as asserted. In particular, Theorem \procref{DefinitionOfTheDegree} holds for generic metrics for any fixed, constant prescribing function.
\newhead{Applications}[Applications]
\newsubhead{The generalised Simons' formula}[GeneralisedSimonsFormula]
We now apply the degree theory developed in Section \headref{DegreeTheory} to the study of hypersurfaces of prescribed curvature in various settings. Each of our applications will depend on a corresponding compactness theorem. To begin with, we introduce some relatively straightforward results which follow directly from a generalised Simons' type formula. In order to understand the concepts involved, we find it informative to first review how the generalised Simons' formula serves to prove elementary Hopf type theorems for constant curvature immersions inside ambient spaces of constant curvature. Indeed, the compactness results that will interest us in Sections \subheadref{SpheresOfPrescribedMeanCurvature} to \subheadref{ExtrinsicCurvature}, below, follow essentially by perturbations of the arguments underlying these results.
\par
Let $\Sigma:=S^d$ be the standard $d$-dimensional sphere. Let $M:=M^{d+1}$ be a smooth, riemannian manifold of dimension $(d+1)$. Let $e:\Sigma\rightarrow M$ be an immersion, and let $A$ be its shape operator. We denote by $\overline{g}$, $\overline{\nabla}$ and $\mR$ respectively the metric, the Levi-Civita covariant derivative and the Riemann curvature tensor of the ambient space. We denote by $g$, $\nabla$ and $R$ respectively the metric induced over $\Sigma$ by the immersion $e$, its Levi-Civita covariant derivative, and its Riemann curvature tensor. Throughout the sequel, we adopt the convention that orders the indices of the Riemann curvature tensor so that
$$
\curvR[ijkl]\partial_l = R(\partial_i,\partial_j)\partial_k =
\nabla_{\partial_i}\nabla_{\partial_j}\partial_k - \nabla_{\partial_j}\nabla_{\partial_i}\partial_k - \nabla_{[\partial_i,\partial_j]}\partial_k.
$$
We have
\proclaim{Lemma \nextprocno, {\bf Generalised Simons' formula.}}
\noindent If $M$ has constant sectional curvature, then, for all $p$ and for all $q$,
$$
A_{pp;qq} = A_{qq;pp} + (A_{pp} - A_{qq})(A_{pp}A_{qq} + \mR_{qppq}),
$$
where ``$;$'' here denotes Levi-Civita covariant differentiation over $\Sigma$.\footnote*{Here, and in all that follows, we adopt the convention that orders indices so that, for any tensor $\alpha$, $\alpha_{p_1...p_m;ij}=(\nabla_{\partial_j}\nabla_{\partial_i}\alpha)(\partial_{p_1},...,\partial_{p_m})$.}
\endproclaim
\proclabel{SimonsFormulaI}
\proof Indeed, by the Codazzi-Mainardi equation, $A_{ij;k}$ is symmetric under all permutations of the indices. Thus, recalling the definition of curvature,
$$\eqalign{
A_{ij;kl} &= A_{kj;il}\cr
&= A_{jk;li} + \curvR[ilkp]A_{pj} + \curvR[iljp]A_{kp}\cr
&= A_{lk;ji} + \curvR[ilkp]A_{pj} + \curvR[iljp]A_{kp}.\cr}
$$
However, by Gauss' equation,
$$
R_{ijkl} = \mR_{ijkl} + A_{il}A_{jk} - A_{ik}A_{jl}.
$$
so that
$$\eqalign{
A_{ij;kl} &= A_{lk;ji}\cr
&\qquad + A_{ij}^2A_{lk} - A_{ik}A_{lj}^2 + A_{ik}^2A_{lj} - A_{ij}A_{lk}^2\cr
&\qquad + \curvmR[ilkp]A_{pj} + \curvmR[iljp]A_{kp},\cr}
$$
and the result now follows upon substituting $i=j=p$ and $k=l=q$.\qed
\medskip
We will apply this formula to the principal curvatures of $e$. However, some care is required as these functions are rarely everywhere twice differentiable. We address this by introducing the formalism of viscosity solutions (c.f. \cite{CrandhallIshiiLyons}). Consider a domain, $U\subseteq\Bbb{R}^d$, and a lower semi-continuous\footnote{$\dagger$}{Recall that $f:X\rightarrow\Bbb{R}$ is lower semi-continuous whenever $f(x)\leq\mliminf_{y\rightarrow x} f(y)$ for all $x\in X$. It is a straightforward exercise of point-set topology to show that a function, $f$, is lower semi-continuous if and only if, for every compact subset $K$ of $X$, and for every continuous function $g$, the restriction of $(f-g)$ to $K$ attains its minimum value at some point. In this way, lower semicontinuity shows itself to be well adapted to the theory of viscosity solutions.} function, $f:\overline{U}\rightarrow\Bbb{R}$. The {\sl subdifferential} of $f$ at any point, $x$, of $U$ is given by
$$
J^2f(x) := \left\{ (\phi(x),D\phi(x),D^2\phi(x))\ \left|\
\matrix
\phi\in C^2(U)\ \text{and}\hfill\cr
(f-\phi)\ \text{attains a local minimum value of $0$ at}\ x\hfill\cr\endmatrix\right.\right\}.
$$
In particular, this set may be empty. Consider now a second-order, linear, partial differential operator,
$$
L:=a^{ij}(x)\partial_i\partial_j + b^i(x)\partial_i + c(x),
$$
defined over $U$, and a real number, $d$. We say that $Lf(x)\leq d$ in the {\sl viscosity sense} whenever
$$
a^{ij}(x)(\partial_i\partial_j\phi)(x) + b^i(x)(\partial_i\phi)(x) + c(x)\phi(x) \leq d,
$$
for all $(\phi(x),D\phi(x),D^2\phi(x))\in J^2f(x)$. The superdifferential of $f$ and the meaning of the inequality $Lf(x)\geq d$ in the viscosity sense are defined in an analogous manner. The interest of this formalism lies in the fact that it provides a general framework within which the maximum principle still applies. Indeed, Hopf's maximum principle still holds with only minor modifications to the proof (c.f. Lemma $3.4$ of \cite{GilbTrud}), yielding
\proclaim{Lemma \nextprocno, {\bf Strong maximum principle.}}
\noindent Suppose that $U$ is connected, that $L$ is elliptic with continuous coefficients and that $c=0$. Let $f:\overline{U}\rightarrow\Bbb{R}$ be a lower semi-continuous function. If $f$ solves
$$
Lf\leq 0,
$$
in the viscosity sense, then
$$
\minf_{x\in\overline{U}}f(x) = \minf_{x\in\partial U}f(x).
$$
Furthermore, if $f$ attains its minimum at any interior point of $U$, and if $f$ is differentiable at this point, then it is constant over the whole of $U$.
\endproclaim
In order to use this formalism in the current setting, we introduce the following construction. Consider a point $x\in\Sigma$. Let $\partial_1,...,\partial_d$ be the standard coordinate vectors of an exponential chart about $x$. In particular, for all $i$, $\nabla\partial_i(x)=0$. Let $\lambda_1\leq...\leq\lambda_d$ be the principal curvatures of $e$ at $x$, and let $v_1,...,v_d$ be the corresponding principal directions. We extend $v_1,...,v_d$ to a frame in a neighbourhood of $x$ by parallel transport along geodesics leaving $x$. In particular, for all $i$, and for all $j$,
$$
\nabla_{\partial_i}v_j = \nabla_{\partial_i}\nabla_{\partial_i}v_j=0.
\eqnum{\nexteqnno[VanishingDerivatives]}
$$
For each $i$, we define the function $a_i$ in a neighbourhood of $x$ by
$$
a_i(y) := \langle A(y)v_i(y),v_i(y)\rangle.
\eqnum{\nexteqnno[DefinitionOfFunctionA]}
$$
In particular, if we now consider $\lambda_1$ and $\lambda_d$ as functions defined in a neighbourhood of $x$, then $a_1\geq\lambda_1$, and $(a_1-\lambda_1)$ attains a minimum value of $0$ at $x$, whilst $a_d\leq\lambda_d$, and $(a_d-\lambda_d)$ attains a maximum value of $0$ at $x$.
\proclaim{Lemma \nextprocno}
\noindent Let $f:\Sigma\rightarrow]0,\infty[$ be a smooth, positive function. If $\lambda_1/f$ attains a local minimum at $x$, then $\lambda_1$ is differentiable at this point. Likewise, if $\lambda_d/f$ attains a local maximum at $x$, then $\lambda_d$ is differentiable at this point.
\endproclaim
\proclabel{FirstOrderDifferentiabilityOfLambda}
\proof We only prove the first assertion, since the proof of the second is identical. If $\lambda_1/f$ attains a local minimum at $x$, then, for all $y$ near $x$,
$$
\lambda_1(x)/f(x) \leq \lambda(y)/f(y)\leq a_1(y)/f(y),
$$
and differentiability follows.\qed
\proclaim{Lemma \nextprocno}
\noindent For all $i$, and for all $j$, at $x$
$$
a_{i;jj}(x) = A_{ii;jj}(x).
$$
\endproclaim
\proclabel{SecondDerivativeOfA}
\proof Indeed, using \eqnref{VanishingDerivatives}, we obtain
$$\eqalign{
A_{ii;jj} &= \langle\nabla^2A(v_i;\partial_j,\partial_j),v_i\rangle\cr
&=\langle\nabla_{\partial_j}(\nabla A)(v_i;\partial_j),v_i\rangle\cr
&=\langle\nabla_{\partial_j}\nabla_{\partial_j} A(v_i) - \nabla_{\partial_j}A(\nabla_{\partial_j}v_i),v_i\rangle\cr
&=\langle\nabla_{\partial_j}\nabla_{\partial_j} A(v_i) - A(\nabla_{\partial_j}\nabla_{\partial_j}v_i),v_i\rangle\cr
&=\langle\nabla_{\partial_j}\nabla_{\partial_j} A(v_i),v_i\rangle\cr
&=\partial_j\langle\nabla_{\partial_j}A(v_i),v_i\rangle\cr
&=\partial_j\partial_j\langle A(v_i),v_i\rangle - \partial_j\langle A(v_i),\nabla_{\partial_j}v_i\rangle\cr
&=\partial_j\partial_j\langle A(v_i),v_i\rangle - \langle A(v_i),\nabla_{\partial_j}\nabla_{\partial_j}v_i\rangle\cr
&=a_{i;jj},}$$
as desired.\qed
\medskip
We now obtain the desired Hopf type theorem. Here the case of non-negative curvature reveals the overall simplicity of our arguments, whilst the case of negative curvature provides an indication of how we will proceed in more general settings. Before stating the results, we recall that an immersion $e$ is said to be {\sl locally strictly convex (LSC)} whenever $\lambda_1>0$ at every point, and, given $c\in]0,1]$, it is then said to be {\sl pointwise $c$-pinched} whenever $\lambda_1\geq cH$ at every point, and {\sl pointwise strictly $c$-pinched} whenever $\lambda_1>cH$ at every point.
\proclaim{Theorem \nextprocno}
\noindent Suppose that $M$ has constant sectional curvature equal to $\kappa\in\left\{-1,0,1\right\}$, and suppose that $e$ is LSC and of constant mean curvature equal to $H>0$.
\medskip
\myitem{(1)} If $\kappa\in\left\{0,1\right\}$, then $e$ is a geodesic sphere.
\medskip
\myitem{(2)} If $\kappa=-1$, and if $H>1$, then $e$ is a geodesic sphere provided that it is pointwise $c_0$-pinched, where $c_0$ is the unique root in $]0,1[$ of the quadratic polynomial,
$$
\frac{x(d-x)}{(d-1)} = \frac{1}{H^2}.
$$
\endproclaim
\proclabel{HopfTypeTheorem}
\remark The interest of this theorem lies not so much in the the statement as in the proof, which serves to illustrate the ideas developed in Sections \subheadref{SpheresOfPrescribedMeanCurvature} and \subheadref{ExtrinsicCurvature} below. Indeed, for $\kappa\in\left\{-1,0\right\}$, a well known argument using the Alexandrov reflection principle yields a stronger result. However, we are not aware of a pre-existing proof in the $\kappa=1$ case.
\medskip
\remark The same technique applies to a large class of curvature functions. However, we do not propose to study this further in the current paper.
\medskip
\proof We will show that $e$ is totally umbilic, from which the result follows by the fundamental theorem of hypersurface theory (c.f. \cite{Spivak}). To this end, we will show that, under the above hypotheses, $\Delta\lambda_1\leq 0$ in the viscosity sense, and that this inequality is strict unless $\lambda_1=H$. Thus, let $a_1$ be defined as in \eqnref{DefinitionOfFunctionA}, above. By Lemma \procref{SecondDerivativeOfA}, at the point, $x$,
$$
\Delta a_1 = \sum_{i=1}^d a_{1;ii} = \sum_{i=1}^d A_{11;ii}.
$$
Thus, by the generalised Simons' formula,
$$\eqalign{
\Delta a_1 &= \sum_{i=1}^dA_{ii;11} + \sum_{i=1}^d(\lambda_1 - \lambda_i)(\lambda_1\lambda_i + R_{i11i})\cr
&=dH_{;11} + \sum_{i=1}^d(\lambda_1 - \lambda_i)(\lambda_1\lambda_i + \kappa),\cr}
$$
and the result now follows for the non-negative curvature case, since $H$ is constant and $(\lambda_1-\lambda_i)$ is non-positive for all $i$.
\par
For the negative curvature case, more refined estimates are required. First, we have
$$
\Delta a_1 = dH_{;11} + d\lambda_1^2H - \lambda_1\sum_{i=1}^d\lambda_i^2 - d\lambda_1 + dH.
$$
Letting $c:=\lambda_1/H$ be the pinching factor, the above relation is rewritten as
$$
\Delta a_1 = dH_{;11} + dc^2H^3 - cH\sum_{i=1}^d\lambda_i^2 + d(1-c)H.
$$
However, applying the Cauchy-Schwarz inequality yields
$$
\sum_{i=1}^d\lambda_i^2 = \lambda_1^2 + \sum_{i=2}^d\lambda_i^2 \geq
\lambda_1^2 + \frac{1}{(d-1)}\left(\sum_{i=2}^d\lambda_i\right)^2
=c^2H^2 + \frac{(d-c)^2}{(d-1)}H^2,
$$
and substituting this into the above relation yields
$$
\Delta a_1 \leq dH_{;11} + d(1-c)H^3\left(\frac{1}{H^2} - \frac{c(d-c)}{(d-1)}\right).\eqnum{\nexteqnno[BoundOnDeltaA]}
$$
The result now follows, since $H$ is constant, and the properties imposed on $H$ and $c$ ensure that the remaining term on the right hand side is non-positive.\qed
\medskip
The above proof works by showing that
$$
\Delta\lambda_1\leq F[H](c),\eqnum{\nexteqnno[LaplacianOfFirstEigenvalue]}
$$
in the viscosity sense, where $c:=\lambda_1/H$ is the pinching factor, and $F[H](c)$ is a cubic polynomial in $c$ which is determined by the parameter $H$. Qualitatively, $F[H]$ behaves as shown in Figure \figref{QualitativeBehaviour}. In particular, since $F[H]$ is negative over the interval $]c_0,1[$, $\lambda_1$ cannot have a local minimum at any point where the pinching factor $c$ lies in this interval. It follows that, if $\Sigma$ is already known to be $c_0$-pinched, then the pinching factor must be constant and equal to $1$, and our Hopf-type theorem follows.
\par
\def\placefigure#1#2#3#4#5{%
\medskip%
\midinsert%
\vbox{\line{\hfil#2\epsfxsize=5cm \epsfbox{#3}\qquad\qquad\epsfxsize=5cm \epsfbox{#4}#1\hfil}%
\vskip 0.3cm%
\line{\hfil\sl Figure \nextfigno\ -\ {#5}\hfil}}%
\medskip%
\endinsert}%
\placefigure{%
\placelabel[-10.95][0]{$c_0$}%
\placelabel[-8.85][0]{$1$}%
\placelabel[-4.55][0]{$c_-$}%
\placelabel[-3.05][0]{$c_+$}%
\placelabel[-2.4][0]{$1$}%
}{}{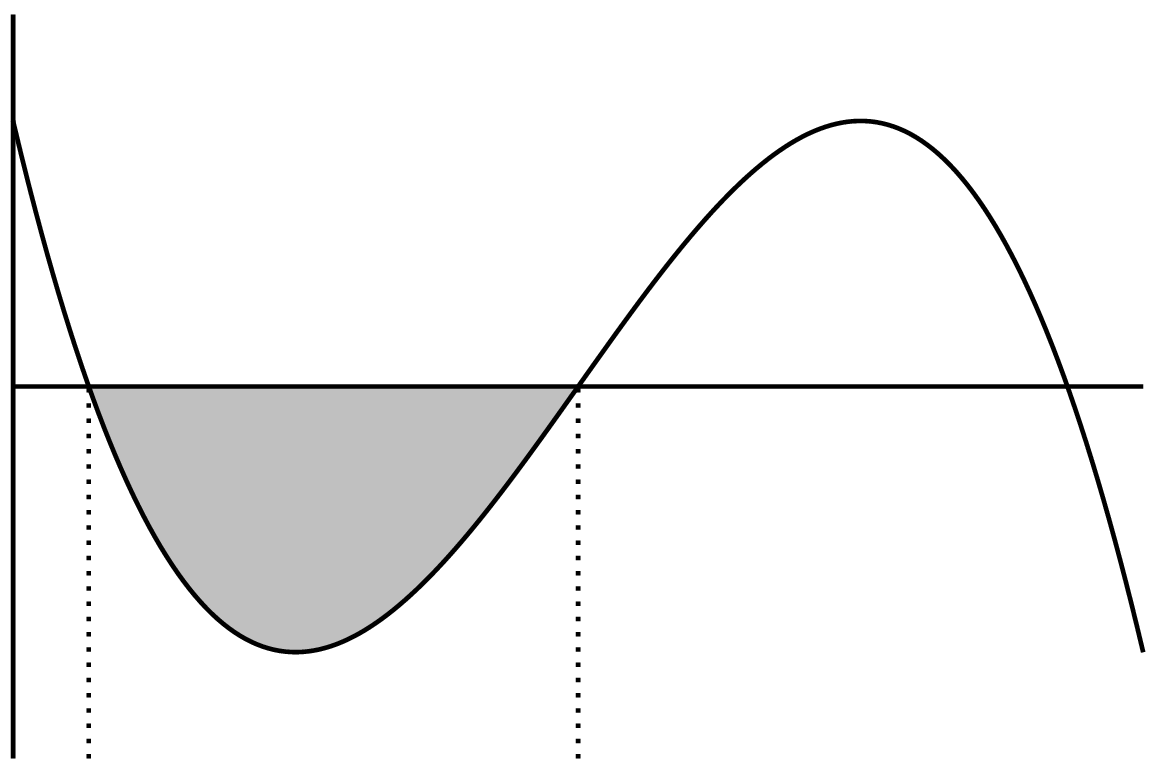}{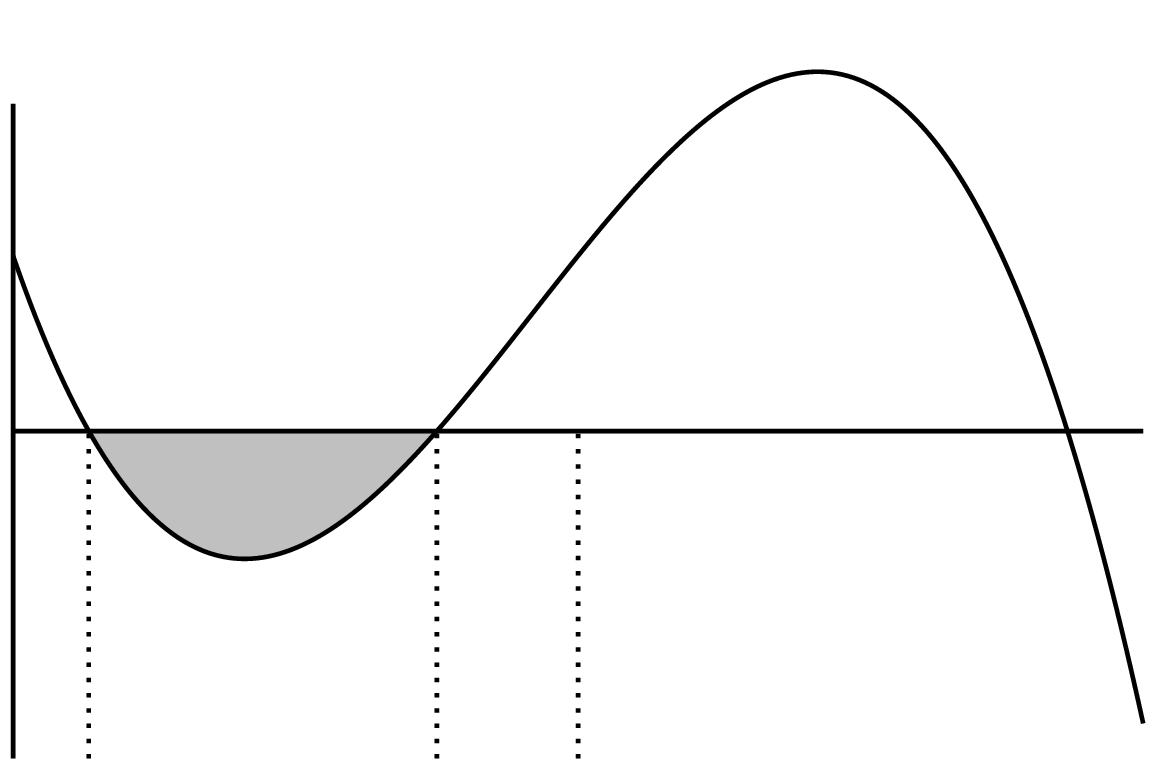}{The qualitative behaviours of $F[H]$ and $F[f,M]$.}
\figlabel{QualitativeBehaviour}
In the more general settings studied in Sections \subheadref{SpheresOfPrescribedMeanCurvature} to \subheadref{ExtrinsicCurvature}, below, we will be interested in hypersurfaces whose curvature is prescribed by some smooth function $f$ inside a general riemannian manifold $M$. Under suitable hypotheses on $f$ and $M$, the relation, \eqnref{LaplacianOfFirstEigenvalue}, perturbs to
$$
\Delta\lambda_1\leq F[f,M](c),
$$
where $F[f,M](c)$ is still a cubic polynomial in $c$, but which now depends on $f$ and $M$, and whose qualitative behaviour is shown in Figure \figref{QualitativeBehaviour}. In particular, since the graph of $F[f,M]$ still dips below $0$ over the subinterval $]c_-,c_+[$ of $]0,1[$, any hypersurface $\Sigma$ of curvature prescribed by $f$ which is $c_-$-pinched, must also be $c_+$-pinched. It is this that will form the basis of the compactness results that allow us to apply the degree theory of Section \headref{DegreeTheory}.
\newsubhead{Prescribed mean curvature}[SpheresOfPrescribedMeanCurvature]
We now study the counting theorems that may be derived from the generalised Simons' formula. We recall the framework of Section \headref{DegreeTheory}. First, let $K$ be mean curvature, that is
$$
K(\lambda_1,...,\lambda_d) := \frac{1}{d}(\lambda_1 + ... + \lambda_d).\eqnum{\nexteqnno[DefnMeanCurvature]}
$$
Let $\Cal{I}$ be the space of LSC unparametrised immersions, where we recall that an immersion is said to be locally strictly convex (LSC) whenever its shape operator is at every point positive definite. Let $\Cal{O}$ be the space of smooth, positive functions $f:M\rightarrow]0,\infty[$ such that, for all $x\in M$,
$$
\frac{2\|\overline{\nabla}^2f(x)\|}{f(x)^3} + \frac{2\|\overline{\nabla}\mR(x)\|}{f(x)^3} + \frac{4\|\mR^o(x)\|}{f(x)^2} + \frac{\|\mR(x)\|}{f(x)^2} < \Gamma,\eqnum{\nexteqnno[ConditionsOnFPMC]}
$$
where
$$
\Gamma :=\msup_{t\in[0,1]}\frac{t(t-1)(t-d)}{(d-1)},\eqnum{\nexteqnno[DefinitionOfM]}
$$
$\mR^o$ denotes the trace free component of $\mR$, that is
$$
\mR^o_{ijkl} := \mR_{ijkl} - \overline{S}(g_{il}g_{jk} - g_{ik}g_{jl}),\eqnum{\nexteqnno[DefinitionOfTraceFreeCurvature]}
$$
and the norm $\|\cdot\|$ of any tensor $\alpha$ is given by
$$
\|\alpha\| = \msup_{\|e_i\|=1\ \forall i}\alpha(e_1,...,e_d),\eqnum{\nexteqnno[DefinitionOfNorm]}
$$
For $f\in\Cal{O}$, define $c_-(f;M)<c_+(f;M)$ to be the two roots in $[0,1]$ of the cubic equation
$$
\frac{c(c-1)(c-d)}{(d-1)} = \msup_{x\in M}\frac{2\|\overline{\nabla}^2f(x)\|}{f(x)^3} + \frac{2\|\overline{\nabla}\mR(x)\|}{f(x)^3} + \frac{4\|\mR^o(x)\|}{f(x)^2} + \frac{\|\mR(x)\|}{f(x)^2},\eqnum{\nexteqnno[KeyPolynomialPEC]}
$$
define the solution space, $\Cal{Z}\subseteq\Cal{I}\times\Cal{O}$, by
$$
\Cal{Z} := \left\{([e],f)\ \left|\
\multiline{H_e = f\circ e;\ \text{and}\cr e\ \text{is strictly $c_-(f;M)$-pinched}\cr}\right.\right\},\eqnum{\nexteqnno[SolutionSpacePMC]}
$$
and let $\Pi:\Cal{Z}\rightarrow\Cal{O}$ be the projection onto the second factor.
\par
Observe that the definition of the solution space used here is subtly different to that given in Section \headref{DegreeTheory}. We therefore define
$$
\hat{\Cal{Z}} := \left\{([e],f)\ |\ H_e = f\circ e\right\},
$$
and in order to apply our degree theory, we first show that $\Cal{Z}$ is a union of connected components of $\hat{\Cal{Z}}$. We achieve this via a modification of the arguments developed in the preceding section. We first require more general versions of some well known formulae. First, we have the Codazzi-Mainardi equations,
$$
A_{ij;k} = A_{kj;i} + \mR_{ki\nu k},\eqnum{\nexteqnno[CodazziMainardi]}
$$
where here ``$;$'' denotes covariant differentiation over $\Sigma$, and $\nu$ denotes the unit normal direction over $e$ that is compatible with the orientation. Next, by definition of the curvature, for any $1$-form $\xi$,
$$
\xi_{k;ij} = \xi_{k;ji} + \curvR[ijkl]\xi_l.\eqnum{\nexteqnno[Curvature]}
$$
Finally, if $f:M\rightarrow\Bbb{R}$ is a twice differentiable function, then
$$
f_{;ij} = f_{:ij} - f_{\nu}A_{ij},\eqnum{\nexteqnno[HessianOfRestriction]}
$$
where here ``$:$'' denotes covariant differentiation over $M$.
\proclaim{Lemma \nextprocno, {\bf Generalised Simons' formula.}}
\noindent Let $e:\Sigma\rightarrow M$ be an immersion, and let $A$ be its shape operator. If $A$ is diagonal at $x$, then, at this point, for all $p$ and for all $q$,
$$\eqalign{
A_{pp;qq} &= A_{qq;pp} + (A_{pp} - A_{qq})(A_{pp}A_{qq} + \mR_{qppq})\cr
&\qquad + A_{pp}(\mR_{pqqp} - \mR_{q\nu\nu q})\cr
&\qquad - A_{qq}(\mR_{pqqp} - \mR_{p\nu\nu p})\cr
&\qquad + \mR_{qp\nu q:p} + \mR_{qp\nu p:q}.}
$$
\endproclaim
\proof Indeed, for all $i$, $j$, $k$ and $l$,
$$\eqalign{
A_{ij;kl} &= A_{kj;il} + (\mR_{ki\nu j})_{;l}\cr
&= A_{jk;li} + \curvR[ilkr]A_{rj} + \curvR[iljr]A_{kr} + (\mR_{ki\nu j})_{;l}\cr
&= A_{lk;ji} + \curvR[ilkr]A_{rj} + \curvR[iljr]A_{kr} + (\mR_{lj\nu k})_{;i} + (\mR_{ki\nu j})_{;l}.\cr}
$$
However, for all $p$, $q$, $r$ and $s$,
$$
(\mR_{pq\nu r})_{;s} = \mR_{pq\nu r:s} - A_{ps}\mR_{\nu q\nu r} - A_{qs}\mR_{p\nu\nu r} + A_s^m\mR_{pqmr}.
$$
Substituting this into the above relation yields,
$$\eqalign{
A_{ij;kl} &= A_{lk;ji} + \curvR[ilkr]A_{rj} + \curvR[iljr]A_{kr}\cr
&\qquad + \mR_{lj\nu k:i} - A_{il}\mR_{\nu j\nu k} - A_{ij}\mR_{l\nu\nu k} + A_i^m\mR_{ljmk}\cr
&\qquad + \mR_{ki\nu j;l} - A_{lk}\mR_{\nu i\nu j} - A_{li}\mR_{k\nu\nu j} + A_l^m\mR_{kimj},\cr}
$$
and the result follows upon substituting $i=j=p$ and $k=l=q$.\qed
\proclaim{Lemma \nextprocno}
\noindent Let $e:\Sigma\rightarrow M$ be an LSC immersion of mean curvature prescribed by $f>0$. At any local minimum of $\lambda_1/f$,
$$\eqalign{
\Delta\left(\frac{\lambda_1}{f}\right)
&\leq
f^2d\left(\frac{(1+c)\|\overline{\nabla}^2f\|}{f^3} + \frac{2\|\overline{\nabla}\mR\|}{f^3} + \frac{2(1+c)\|\mR^o\|}{f^2}\right.\cr
&\qquad\qquad\qquad\qquad + \frac{(1-c)\|\mR\|}{f^2} - \frac{c(c-1)(c-d)}{(d-1)}\Bigg)}
\eqnum{\nexteqnno[KeyInequalityCMC]}
$$
in the viscosity sense, where $c:=\lambda_1/f$ is the pinching factor.
\endproclaim
\proclabel{EstimateOfDeltaPinchingFactorPMC}
\remark Each of the terms in the parenthesis above is scale invariant in the sense that it is unchanged when the metric is multiplied by a constant factor.
\medskip
\proof Let $a_1$ be defined as in \eqnref{DefinitionOfFunctionA}, above. In particular, a local minimum of $\lambda_1/f$ is also a local minimum of $a_1/f$. Now, at any such point,
$$
\nabla\left(\frac{a_1}{f}\right) = \frac{1}{f}\nabla a_1 - \frac{a_1}{f^2}\nabla f = 0,
$$
so that
$$\eqalign{
\Delta\left(\frac{a_1}{f}\right)
&=\frac{1}{f}\Delta a_1 - \frac{a_1}{f^2}\Delta f - \frac{2}{f}\left\langle\frac{1}{f}\nabla a_1 - \frac{a_1}{f^2}\nabla f,\nabla f\right\rangle\cr
&=\frac{1}{f}\Delta a_1 - \frac{a_1}{f^2}\Delta f.\cr}
$$
Repeating the calculation used to obtain \eqnref{BoundOnDeltaA} now yields
$$
\Delta a_1 \leq f^3d\left(\frac{f_{;11}}{f^3} + \frac{2\|\overline{\nabla}\mR\|}{f^3} + \frac{2(1+c)\|\mR^o\|}{f^2} + \frac{(1-c)\|\mR\|}{f^2} - \frac{c(c-1)(c-d)}{(d-1)}\right).
$$
Furthermore, by \eqnref{HessianOfRestriction},
$$
f_{;11} = f_{:11} - f_\nu\lambda_1,
$$
and
$$\eqalign{
\Delta f &= \sum_{i=1}^d f_{;ii}\cr
&= \sum_{i=1}^d(f_{:ii} - f_\nu A_{ii})\cr
&= \left(\sum_{i=1}^df_{:ii}\right) - df_\nu\lambda_1,\cr}
$$
and the result now follows upon combining these relations.\qed
\proclaim{Corollary \nextprocno}
\noindent If $([e],f)\in\Cal{Z}$, then $e$ is $c_+(f;M)$-pinched.
\endproclaim
\proclabel{SolutionsArePinchedPMC}
\proof Indeed, otherwise the pinching factor, $c=\lambda_1/f$, of $[e]$ would attain its minimum value at some point in the interval $]c_-(f;M),c_+(f;M)[$. However, by definition of $c_-$ and $c_+$, the right hand side of \eqnref{KeyInequalityCMC} is negative over this interval, and we have a contradiction by the maximum principle.\qed
\proclaim{Lemma \nextprocno}
\noindent $\Cal{Z}$ is a union of connected components of $\hat{\Cal{Z}}$.
\endproclaim
\proclabel{ConnectedComponentsPMC}
\proof By definition, $\Cal{Z}$ is an open subset of $\hat{\Cal{Z}}$. Consider now a sequence $([e_m],f_m)$ of elements of $\Cal{Z}$ converging towards the element $([e_\infty],f_\infty)$ of $\hat{\Cal{Z}}$. By Corollary \procref{SolutionsArePinchedPMC}, for all $m$, $[e_m]$ is $c_+(f_m,M)$-pinched. Taking limits, it follows that $[e_\infty]$ is $c_+(f_\infty,M)$-pinched, so that $([e_\infty],f_\infty)$ is also an element of $\Cal{Z}$. This set is therefore both open and closed as a subset of $\hat{\Cal{Z}}$, and the result follows.\qed
\proclaim{Lemma \nextprocno}
\noindent If $([e],f)\in\Cal{Z}$, then $[e]$ is prime.
\endproclaim
\proclabel{PrimePMC}
\remark We prove this result using the mean curvature flow developed by Huisken in \cite{Huisken}. The alert reader will notice that although the hypotheses of Huisken's result are stated in terms of the norms of $\mR$ and $\overline{\nabla}\mR$, he is not explicit about {\sl which} norms are being used. Closer examination of the text, and, in particular, lines $11$ to $14$ on $p472$ in the proof of Theorem $4.2$ of that paper, shows, however, that the operator norm introduced above is indeed the correct one.
\medskip
\proof Let $A$ be the shape operator of $e$. Choose a point $x\in\Sigma$ and let $c:=\lambda_1/f$ be the pinching factor of $e$ at this point. Observe that $t\geq t(t-1)(t-d)/(d-1)$ over the interval $[0,1]$. In particular, if we denote by $c_0\in[0,1]$ the unique point of this interval maximising $t(t-1)(t-d)$, then, bearing in mind Corollary \procref{SolutionsArePinchedPMC},
$$
c\geq c_+(f,M)>c_0\geq\frac{c_0(1-c_0)(d-c_0)}{(d-1)}=M.
$$
Thus, since $f$ satisfies \eqnref{ConditionsOnFPMC},
$$
c>\left(\frac{\|\mR\|}{f^2} + \frac{\|\overline{\nabla}\mR\|}{f^3}\right).
$$
Denoting $K:=\|\mR\|$, $L:=\|\overline{\nabla}\mR\|$ and $\hat{H}:=dH$, we thus have
$$\eqalign{
\hat{H}A_{ij} &\geq d\lambda_1 f g_{ij}\cr
&\geq dc f^2 g_{ij}\cr
&\geq \left(dK + \frac{d^2}{\hat{H}}L\right)g_{ij},\cr}
$$
which is exactly the condition given by Huisken in \cite{Huisken} for the existence of a unique, smooth mean curvature flow $\hat{e}:\Sigma\times[0,T[\rightarrow M$ such that $\hat{e}_0=e$, and $(\hat{e}_t)_{t\in[0,T[}$ is asymptotic to a family of round spheres about a point in $M$ as $t\rightarrow T$.
\par
Suppose now that $e=\hat{e}_0$ is not prime. Then there exists a non-trivial diffeomorphism $\alpha:\Sigma\rightarrow\Sigma$ such that $e\circ\alpha=e$. By uniqueness, for all $t\in[0,T[$, $\hat{e}_t\circ\alpha=\hat{e}_t$. However, for $t$ sufficiently close to $T$, $\hat{e}_t$ is embedded. This is absurd, and the result follows.\qed
\proclaim{Lemma \nextprocno}
\noindent For all $f\in\Cal{O}$, there exists $D>0$ such that if $([e],f)\in\Cal{Z}$, then
$$
\opDiam([e]) \leq D.
$$
Furthermore, $D$ can be chosen to vary continuously with $f$.
\endproclaim
\proclabel{UniformDiameterBoundsPMC}
\proof Consider a point $x\in\Sigma$ and let $c=\lambda_1/f$ be the pinching factor of $e$ at $x$. The Ricci curvature of $\Sigma$ at this point satisfies
$$\eqalign{
\opRic_{ik} &\geq \left(\frac{1}{(d-1)}\lambda_1(\lambda_2+...+\lambda_d) - \|\mR\|\right)g_{ik}\cr
&=\left(\frac{c(d-c)}{(d-1)}f^2 - \|\mR\|\right)g_{ik}.\cr}
$$
Observe that the function $t(d-t)$ is increasing over the interval $[0,1]$, and is not less than $t(t-1)(t-d)$ at every point of this interval. In particular, if we denote by $c_0\in[0,1]$ the unique point of this interval maximising $t(t-1)(t-d)$, then, by Corollary \procref{SolutionsArePinchedPMC},
$$
c\geq c_+(f;M)>c_0,
$$
so that
$$
\frac{c(d-c)}{(d-1)} > \frac{c_0(d-c_0)}{(d-1)} > \frac{c_0(c_0-1)(c_0-d)}{(d-1)} = \Gamma.
$$
However, since $f$ satisfies \eqnref{ConditionsOnFPMC},
$$
\|\mR\| < \Gamma f^2 + \epsilon,
$$
for some $\epsilon>0$, so that
$$
\opRic_{ik} \geq (\Gamma f^2 - \|\mR\|) g_{ik} \geq \epsilon g_{ik},
$$
and the result now follows by the Bonnet-Myers theorem.\qed
\proclaim{Lemma \nextprocno}
\noindent The projection $\Pi:\Cal{Z}\rightarrow\Cal{O}$ is a proper map.
\endproclaim
\proclabel{ProperPMC}
\proof Let $(f_m)\in\Cal{O}$ be a sequence converging towards $f_\infty\in\Cal{O}$, and for all $m$, let $[e_m]$ be an immersion such that $([e_m],f_m)\in\Cal{Z}$. By Lemma \procref{UniformDiameterBoundsPMC}, there exists $D>0$ such that, for all $m$,
$$
\opDiam([e_m]) \leq D.
$$
By strict convexity, for all $m$,
$$
\|A_m\| \leq \|f_m\|,
$$
where $A_m$ here denotes the shape operator of $e_m$. It follows by the Arzela-Ascoli theorem for immersed hypersurfaces (c.f. \cite{SmiAAT}) and elliptic regularity, that there exists an unparametrised immersion $[e_\infty]$ towards which $([e_m])$ subconverges. By Corollary \procref{SolutionsArePinchedPMC}, for all $m$, $[e_m]$ is $c_+(f,M)$-pinched. Taking limits, it follows that $[e_\infty]$ is also $c_+(f_\infty,M)$-pinched, so that $([e_\infty],f_\infty)$ is also an element of $\Cal{Z}$. It follows that $\Pi$ is a proper map, as desired.\qed
\newsubhead{Calculating the Degree}[CalculatingTheDegree]
Lemmas \procref{ConnectedComponentsPMC}, \procref{PrimePMC} and \procref{ProperPMC} allows us to apply the degree theory developed in Section \headref{DegreeTheory}. It thus remains only to calculate the degree. This is achieved by considering a special case where hypersurfaces of prescribed curvature can be explicitly counted. To this end, we study the case where the curvature is prescribed by the function $f_t:=(1+t^2f)/t$, where $0<t\ll 1$, and $f:M\rightarrow\Bbb{R}$ is a smooth function whose properties we will describe presently.
\par
The hypersurfaces of curvature prescribed by the function $f_t$ are analysed via a modification of the asymptotic construction developed by Ye in \cite{Ye} which we now describe. The details are presented in full in \cite{SmiEC}. It is first necessary to review Ye's construction in some detail. For simplicity, we first suppose that the function $f$ vanishes and that the scalar curvature function $R$ of the ambient manifold $M$ is of Morse type. The general case will be addressed towards the end of this section. Let $p$ be a point of $M$, let $\Omega$ be a convex, normal neighbourhood of $p$, identify $T_pM$ with $\Bbb{R}^{d+1}$ furnished with the euclidean metric, and let $S^d$ be the unit sphere in this space. Consider the function
$$
E:]0,\infty[\times C^{k+2,\alpha}(S^d)\times\Omega\rightarrow C^{k+2,\alpha}(S^d,M)
$$
defined such that for all $x\in S^d$,
$$
E(t,\phi,q)(x) = \opExp_q(t(1+t^2 \phi(x))T_{q,p}x),
$$
where $\opExp_q$ is the exponential map of $M$ about $q$, and $T_{q,p}$ is the parallel transport along the unique geodesic in $\Omega$ from $p$ to $q$. This is a smooth map between Banach manifolds. Furthermore, for any given $(\phi,q)$, $E(t,\phi,q)$ defines a $C^{k+2,\alpha}$ embedding of $S^d$ into $M$ provided that $t$ is sufficiently small. Throughout the rest of this section, we will suppose that this is the case. For notational convenience, we will denote
$$
e_\mu := E(t,\phi,q),
$$
where $\mu:=(t,\phi,q)$.
\par
Consider now the function
$$
H:]0,\infty[\times C^{k+2,\alpha}(S^d)\times\Omega\rightarrow C^{k,\alpha}(S^d)
$$
defined such that for all $\mu:=(t,\phi,q)$, and for all $x$, $H(t,\phi,q)(x)$ is the mean curvature of the embedding $e_\mu$ at the point $e_\mu(x)$. As before, this is a smooth map between Banach manifolds. In \cite{Ye}, Ye shows that
$$\eqalign{
H(t,\phi,q)(x) &= \frac{1}{t} - \frac{t}{d}\left(d + \Delta\right)\phi(x)\cr
&\qquad -\frac{t}{3}\opRic(q)\left(T_{q,p}x,T_{q,p}x\right)\cr
&\qquad -\left[\frac{t^2}{4}\nabla\opRic(q)\left(T_{q,p}x,T_{q,p}x,T_{q,p}x\right) - \frac{(d+1)}{2(d+3)}\nabla R(q)\left(T_{q,p}x\right)\right]\cr
&\qquad -\frac{t^2(d+1)}{2(d+3)}\nabla R(q)\left(T_{q,p}x\right)\cr
&\qquad +O(t^3)\phantom{\bigg[ },\cr}
\eqnum{\nexteqnno[AsymptoticFormulaForHI]}
$$
where $\Delta$ is the Laplace-Beltrami operator of $S^d$. Here we use the convention that a function $F$ of $t$ and some other variable $\nu$, say, is $O(t^m)$ whenever there exists a smooth function $G:=G(t,\nu)$ defined in a neighbourhood of $\left\{t=0\right\}$ such that $F(t,\nu)=t^mG(t,\nu)$.
\par
The subtleties of Ye's construction can now be understood in terms of the spectrum of $\Delta$. Indeed, $(-d)$ is an eigenvalue of $\Delta$ whose eigenspace $\Lambda$ is the space of all restrictions to $S^d$ of linear functions over $\Bbb{R}^{d+1}$. In particular, $\Lambda$ is the kernel of $(d+\Delta)$, and it follows by self-adjointness and standard elliptic theory (c.f. \cite{GilbTrud}) that its $L^2$ orthogonal complement $\Lambda^\perp$ is the image of this operator. Thus, since the functions in the second and third lines of \eqnref{AsymptoticFormulaForHI} are elements of $\Lambda^\perp$ for all $q$, there exist unique smooth functions $\phi_0,\phi_1:S^d\times\Omega\rightarrow\Bbb{R}$ such that, for all $q$, $\phi_{0,q}:=\phi_0(\cdot,q)$ and $\phi_{1,q}:=\phi_1(\cdot,q)$ are both elements of $\Lambda^\perp$, and
$$\eqalign{
\frac{1}{d}\left(d+\Delta\right)\phi_{0,q}(x) &= -\frac{t}{3}\opRic(q)\left(T_{q,p}x,T_{q,p}x\right),\cr
\frac{1}{d}\left(d+\Delta\right)\phi_{1,q}(x) &= -\left[\frac{t^2}{4}\nabla\opRic(q)\left(T_{q,p}x,T_{q,p}x,T_{q,p}x\right) - \frac{(d+1)}{2(d+3)}\nabla R(q)\left(T_{q,p}x\right)\right].\cr}
$$
Substituting this into \eqnref{AsymptoticFormulaForHI} yields
$$\eqalign{
\frac{1}{t^2}\left[H\left(t,\phi_{0,q}+t\phi_{1,q} + t\phi,q\right) - \frac{1}{t}\right]
&= -\frac{1}{d}(d+\Delta)\phi(x)\cr
&\qquad -\frac{(d+1)}{2(d+3)}\nabla R(q)\left(T_{q,p}x\right)\cr
&\qquad +O(t).\phantom{\bigg[ }\cr}\eqnum{\nexteqnno[AsymptoticFormulaForHII]}
$$
Since the term in the second line of \eqnref{AsymptoticFormulaForHII} is an element of $\Lambda$, it cannot be removed in the same way. Thus, in order for this expression to vanish at $\left\{t=0\right\}$, it is necessary that $p$ be a critical point of the scalar curvature function $R$. In this case, since $R$ is of Morse type, $\opHess(R)(p)$ is non-degenerate, and it follows that the partial derivative of \eqnref{AsymptoticFormulaForHII} with respect to the second and third components $(\phi,q)$ is surjective at the point $(0,0,p)$. By the implicit function theorem for smooth functions defined over Banach manifolds, there therefore exist $C,t_0>0$ with the property that, for all $t\in]0,t_0[$, there exists a unique function $\phi_t\in\Lambda^\perp$ and a unique point $q_t\in\Omega$ such that
$$\eqalign{
\|\phi_t\|_{k+2,\alpha}&<Ct,\cr
d(p,q_t)&<Ct,\cr}$$
and
$$
H\left(t,\phi_{0,q_t} + t\phi_{1,q_t} + t\phi_t,q_t\right) = \frac{1}{t}.
$$
In other words, the embedded sphere $e_{\mu(t)}(S^d)$ has constant mean curvature equal to $1/t$, where
$$
\mu(t):=(t,\phi_{0,q_t} + t\phi_{1,q_t} + t\phi_t,q_t).
$$
In particular, by elliptic regularity, $e_{\mu(t)}$ is smooth for all $t$. This completes the part of Ye's construction that is of relevance to us. In \cite{Ye}, Ye also shows that the family $(e_{\mu(t)})_{t\in]0,t_0[}$ foliates a neighbourhood of $p$ (minus the point, $p$, itself, of course). There is also an elegant variant of Ye's argument, developed by Pacard \& Xu in \cite{PacardXu} which constructs embedded spheres of large constant mean curvature even when $p$ is a degenerate critical point of $R$. However, we will make no use of that result here.
\par
Let $\opCrit(R)$ now denote the set of critical points of $R$. In the present context, Ye's construction defines, for all sufficiently small $t$, a canonical injection $Y_t$ from $\opCrit(R)$ into $\Pi^{-1}(\left\{1/t\right\})$. The next step consists of showing that, for sufficiently small $t$, $Y_t$ is also surjective. To this end, consider a sequence $(t_m)$ of positive real numbers, a sequence $(e_m)$ of smooth embeddings of $S^d$ into $M$, and a sequence $(p_m)$ of points in $M$ such that, for all $m$, $[e_m]$ is an element of $\Pi^{-1}(\left\{1/t_m\right\})$ and $p_m$ lies at a distance of less than $\opDiam([e_m])$ from $[e_m]$.
\proclaim{Lemma \nextprocno}
\noindent For sufficiently large $m$,
$$
\opDiam([e_m]) \leq \frac{\sqrt{2}\pi}{t_m}.
$$
\endproclaim
\proclabel{TheRadiiTendLinearlyToZero}
\proof Indeed, observe that $c_+(t,M)$ tends to $1$ as $t$ tends to infinity. Thus, for sufficiently large $m$, tthe Ricci curvature of $e_m^*g$ satisfies
$$
\opRic_{ik} \geq \frac{t^2}{2}g_{ik}
$$
at every point of $S^d$. Thus, by the Bonnet-Myers Theorem,
$$
\opDiam([e])\leq \frac{\sqrt{2}\pi}{t},
$$
as desired.\qed
\medskip
Denote now the metric of $M$ by $g$, and, for all $m$, denote $g_m:=g/t_m$. Observe that the sequence $(M,g_m,p_m)$ of complete, pointed riemannian manifolds converges in the $C^\infty$ Cheeger-Gromov sense to the complete, pointed riemannian manifold $(\Bbb{R}^{d+1},g_0,0)$ where $g_0$ here denotes the standard euclidean metric over $\Bbb{R}^{d+1}$. Furthermore, for all sufficiently large $m$, the diameter of $[e_m]$ with respect to the rescaled metric $g_m$ is bounded above by $\sqrt{2}\pi$. Thus, as in the proof of Lemma \procref{ProperPMC}, we may suppose that there exists an umparametrised immersion $[e_\infty]$ in $\Bbb{R}^{d+1}$ towards which $([e_m])$ converges in the $C^\infty$ sense. Furthermore, by Hopf's Theorem (or by a suitable adaptation of Corollary \procref{SolutionsArePinchedPMC}, though c.f. also Theorem \procref{HopfTypeTheorem}), $[e_\infty]$ is a unit sphere. Upon modifying the sequence, $(p_m)$, we may suppose that this sphere is centred on the origin, and that, for all sufficiently large $m$, $[e_m]$ is a graph over this sphere of some function in $\Lambda^\perp$. That is, there exists a sequence $(\phi_m)$ of smooth functions in $\Lambda^\perp$, and a sequence $(q_m)$ of points in $\Omega$ such that $(q_m)$ converges to some point $p$ of $M$, $(\phi_m)$ converges to $0$, and, after reparametrisation, for all sufficiently large $m$, and for all $x$,
$$
e_m(x) = \opExp_{q_m}\left(t_m(1 + \phi_m(x))T_{q_m,p}x\right).
$$
\par
It remains to show that $p$ is a critical point of $R$ and that, for all sufficiently large $m$, $\phi_m$ and $q_m$ are of the form given above. However, consider the function
$$
\tilde{E}:]0,\infty[\times C^{k+2,\alpha}(S^d)\times\Omega\rightarrow C^{k+2,\alpha}(S^d,M)
$$
given by
$$
\tilde{E}(t,\phi,q)(x) := \opExp_p(t(1+\phi(x))T_{q,p}x),
$$
and define the function
$$
\tilde{H}:]0,\infty[\times C^{k+2,\alpha}(S^d)\times\Omega\rightarrow C^{k,\alpha}(S^d)
$$
such that, for all suitable $(t,\phi,q)$, $\tilde{H}(t,\phi,q)(x)$ is the mean curvature function of the embedding $\tilde{E}(t,\phi,q)$. Calculating as before the asymptotic expansion of $\tilde{H}$ about $\left\{t=0\right\}$, we now obtain
$$
\tilde{H}(t,\phi,q)(x) - \frac{1}{t} = -\frac{1}{td}(d + \Delta)\phi(x) + \epsilon(t,\phi,q)(x),
$$
where the error term $\epsilon$ satisfies
$$
\|\epsilon(t,\phi,q)\|_{k,\alpha} \leq B_1\left(1 + \|\phi\|_{k+2,\alpha}^2\right),
$$
for some suitable constant $B_1>0$, and for all sufficiently small $(t,\phi)$. Thus, since $e_m$ has constant mean curvature equal to $1/t_m$ for all $m$, since $\phi_m\in\Lambda^\perp$ for all $m$, and since $(\phi_m)$ converges to $0$, it now follows by standard elliptic theory (c.f. \cite{GilbTrud}) that, for all $m$,
$$
\|\phi_m\|_{k+2,\alpha} \leq B_2 t_m,
$$
for some suitable constant $B_2>0$.
\par
In other words, elliptic theory provides an improved estimate of the rate at which the sequence $(\phi_m)$ converges to $0$. This constitutes the first stage of an elliptic bootstrapping argument which determines both the sequences $(\phi_m)$ and $(q_m)$ up to arbitrarily high order in $t_m$ as $m$ tends to infinity. The remaining stages of this argument are similar, and are presented in full detail in \cite{SmiEC} (though c.f. also \cite{SmiAC} for a similar argument in a more straightforward context). They show that $p$ is a critical point of $R$ and that, for sufficiently large $m$, both $\phi_m$ and $q_m$ are of the desired form, so that, for sufficiently small $t$, $Y_t$ does indeed define a surjective map from $\opCrit(R)$ onto $\Pi^{-1}(\left\{1/t\right\})$, as desired.
\par
The final step of calculating the degree consists of calculating the Morse indices of the embedded hyperspheres constructed above. To this end, for all $\mu:=(t,\phi,q)$, let $J_\mu$ be the Jacobi operator of the embedding $e_\mu$. The same asymptotic analysis as before now yields
$$
J_\mu\psi = -\frac{1}{dt^2}(d+\Delta)\psi + O\left(\frac{1}{t}\right).
$$
This formula already allows us to determine most of the spectrum of $J_\mu$. Indeed, the only negative eigenvalue of the operator $-\frac{1}{d}(d+\Delta)$ is $(-1)$, which has unit multiplicity, since its eigenspace consists of all constant functions over $S^d$. Next, as shown above, its kernel is $\Lambda$, which is $(d+1)$-dimensional. Finally, all the rest of its eigenvalues are strictly positive. It thus follows by standard perturbation theory (c.f. \cite{Kato}) that the negative eigenvalues of the operator $J_\mu$ are determined by studying how the degenerate eigenvalue $0$ perturbs.
\par
For all $t$, with $\phi_t$, $q_t$ and $\mu(t)$ as before, consider now $J_{\mu(t)}$, the Jacobi operator of the embedding $e_{\mu(t)}$. By \cite{Kato}, there exists a smooth family $\Lambda_t$ of $(d+1)$-dimensional subspaces of $L^2(S^d)$ such that $\Lambda_0=\Lambda$, and, for all sufficiently small $t$,
$$
J_{\mu(t)}\Lambda_t \subseteq \Lambda_t.
$$
In fact, it follows from elliptic regularity that $\Lambda_t$ actually consists of smooth functions. Consider now a tangent vector $X\in T_p M$ and observe that the function $\langle X,\cdot\rangle$ is an element of $\Lambda$. For all sufficiently small $t$, $\Lambda_t$ is a graph over $\Lambda$, and we therefore define $\psi_{t,X}\in L^2(S^d)$ such that, $\psi_{t,X}\in\Lambda_t$ and, for all $x$,
$$
\Pi_1\left(\psi_{t,X}\right)(x) = \langle X,x\rangle,
$$
where $\Pi_1:L^2(S^d)\rightarrow\Lambda$ is the orthogonal projection. We now define the function
$$
a:]0,t_0[\times\Lambda\times\Lambda\rightarrow\Bbb{R}
$$
by
$$
a_t(X,Y) := a(t,X,Y) := \langle J_{\mu(t)}\psi_{t,X},\psi_{t,Y}\rangle.
$$
It follows by standard perturbation theory (c.f. \cite{Kato}) that those eigenvalues of $J_{\mu(t)}$ which arise through perturbations of the degenerate eigenvalue $0$ of $-(d+\Delta)/dt^2$ coincide with the eigenvalues of $a_t$.
\par
It turns out to be necessary to determine $a_t$ up to and including order $2$ in $t$. A direct calculation of this would require knowledge of $e_{\mu(t)}$ up to and including order $4$. However, this is rather difficult, especially when more general curvature functions are considered later on. In order to bypass this, we therefore consider the following modification of Ye's construction. Observe that, when $p$ is not a critical point of $R$, but when $\opHess(R)(p)$ is nonetheless non-degenerate, Ye's construction still applies to some extent. In this manner, we extend $\mu(t)$ to a smooth function $\mu:]0,t_0[\times\Omega\rightarrow]0,t_0[\times C^{k+2,\alpha}(S^d)\times\Omega$ such that, for all $(t,q)$, the function $e_{\mu(t,q)}$ is a smooth embedding whose mean curvature at the point $e_{\mu(t,q)}(x)$ is equal to
$$
\frac{1}{t} - \frac{t^2(d+1)}{2(d+3)}\nabla R(q)\left(T_{q,p}x\right).\eqnum{\nexteqnno[MeanCurvatureOfNearbySpheres]}
$$
For all $t$, we now define $\overline{\psi}_{t,X}\in C^\infty(S^d)$ by
$$
\overline{\psi}_{t,X}(X) := \left\langle\frac{\partial}{\partial s} e_{\mu(t,\opExp_p(sX))}|_{s=0},N_t(x)\right\rangle.
$$
where $\opExp_p$ is the exponential map of $M$ at $p$, and $N_t$ is the outward pointing unit normal vector field of the embedding $e_{\mu(t)}=e_{\mu(t,p)}$ at the point $e_{\mu(t)}(x)$. By definition of the Jacobi operator, $J_{\mu(t)}\overline{\psi}_{t,X}$ is now determined by differentiating \eqnref{MeanCurvatureOfNearbySpheres} with respect to $q$ at $p$, and we thus obtain, for all $t$,
$$
J_{\mu(t)}\overline{\psi}_{t,X}(x) = -\frac{t^2(d+1)}{2(d+3)}\opHess(R)(X,x).
$$
However, by explicit calculation of the lower order terms of $\psi_{t,X}$ and $\overline{\psi}_{t,X}$ (c.f. Section $4$ of \cite{SmiEC}),
$$
\psi_{t,X}(x) - \overline{\psi}_{t,X}(x) = t^2\langle V_t,x\rangle + O(t^3),
$$
for some tangent vector $V_t\in T_pM$.
\par
We are now able to determine the first non-trivial term in the asymptotic expansion of $a_t$. Indeed, by explicit calculation,
$$
J_{\mu(t)} = -\frac{1}{dt^2}\left(d + \Delta\right) + M_0 + O(t),
$$
where $M_0$ maps $\Lambda$ into $\Lambda^\perp$. From this and the above relations, it readily follows that, for all vectors $X$ and $Y$,
$$
a_t(X,Y) = -\frac{t^2(d+1)}{2(d+3)}\opHess(R)(p)(X,Y) + O(t^3).
$$
We conclude that for every critical point $p$ of $R$, and for all sufficiently small $t$, the Morse index of $Y_t(p)$ is related to the Morse index of $p$ by
$$
\opInd(Y_t(p)) = 1 + \left((d+1)-\opInd(p)\right),
$$
and the signatures are thus related by
$$
\opSig(Y_t(p)) = (-1)^d\opSig(p).
$$
Summing over the whole of $\opCrit(R)$ therefore yields
$$
\opDeg(\Pi) = (-1)^d\chi(M),
$$
where $\chi(M)$ is the Euler characteristic of $M$, and since this vanishes when the ambient space is odd-dimensional, that is, when $d$ is even, this simplifies to
$$
\opDeg(\Pi) = -\chi(M).
$$
\par
We now sketch two straightforward modifications to Ye's construction which allow us to treat more general cases. First, since any curvature function, $K$, coincides with $H$ up to and including order $1$ at the point $(1,...,1)$ most of the preceding analysis applies without modification with $H$ replaced by $K$. Care, however, should be taken in proving the surjectivity of $Y_t$ to ensure that the asymptotic limits exist and that they are indeed round spheres. This will be the case for all the examples studied in this paper, and in this manner, we obtain the degree for all the curvature functions considered in the sequel. Next, when the scalar curvature function $R$ is not of Morse type, we can proceed in one of two ways. Either we can choose to take an arbitrarily small conformal perturbation of $g$ whose scalar curvature function is of Morse type, or we can choose to proceed as follows. Let $f:M\rightarrow\Bbb{R}$ be a smooth function chosen such that
$$
R_f := R + \frac{2(d+3)}{(d+1)}f
$$
is of Morse type, and for all large $t$, consider the function
$$
f_t:=\frac{1}{t}(1 + t^2f).
$$
Ye's construction readily modifies again to show that if $p$ is a critical point of $R_f$ - which is non-degenerate by the Morse property - then there exists $C,t_0>0$ such that for all $t\in]0,t_0[$, there exists a unique function $\phi\in\Lambda^\perp$ and a unique point $q_t\in\Omega$ such that $\|\phi\|_{k+2,\alpha}<Ct$, $d(p,q_t)<Ct$, and such that the smooth embedding
$$
e_t := E(t,\tilde{\phi}_{0,f} + t\tilde{\phi}_{1,f} + t\phi,q_t)
$$
has $K$-curvature prescribed by $f_t$, where the smooth functions $\tilde{\phi}_{0,f}$ and $\tilde{\phi}_{1,f}$ are defined in a similar manner as before. The remainder of the analysis continues unchangeed so that, by taking the sum of the signatures of the elements of $\Pi^{-1}(\left\{f_t\right\})$ for sufficiently small $t$, we recover $\opDeg(\Pi)=-\chi(M)$, as desired.
\par
To summarise, we ave now constructed the degree for the case of embedded hypersurfaces of prescribed mean curvature.
\proclaim{Theorem \nextprocno}
\noindent For generic $f\in\Cal{O}$, the algebraic number of strictly $c_-(f,M)$-pinched, immersed hyperspheres of mean curvature prescribed by $f$ is equal to $-\chi(M)$.
\endproclaim
\proclabel{ExistencePMCI}
\remark Theorem \procref{MainTheoremMeanCurvature} now follows via the discussion of Section \subheadref{VaryingTheMetric}.
\medskip
\proof By Lemmas \procref{ConnectedComponentsPMC}, \procref{PrimePMC} and \procref{ProperPMC}, the degree theory developed in Section \headref{DegreeTheory} applies, and by the preceding discussion,
$$
\opDeg(\Pi) = -\chi(M),
$$
as desired.\qed
\medskip
Finally, by comparing the definitions \eqnref{KeyPolynomialPEC} and \eqnref{SolutionSpacePMC} with \eqnref{KeyInequalityCMC}, we see that the definitions of $\Cal{O}$ and $\Cal{Z}$ are not quite optimal. For this reason, different approaches will yield complementary results that are not entirely contained in Theorem \procref{ExistencePMCI}. For example, fixing the pinching factor $c$ to be equal to $1/2$, we define
$$
H_0 := 4\opMax\left(\|\mR\|^{\frac{1}{2}},\|\overline{\nabla}\mR\|^{\frac{1}{3}}\right),
$$
and we leave the reader to verify that the same reasoning now yields
\proclaim{Theorem \nextprocno}
\noindent For generic $f$ such that
$$\eqalign{
f &> H_0,\ \text{and}\cr
\|\overline{\nabla}^2 f\| &< \frac{3d}{(3d-2)}H_0^2,\cr}
$$
the algebraic number of strictly $(1/2)$-pinched, immersed hyperspheres of prescribed mean curvature equal to $f$ is equal to $-\chi(M)$.
\endproclaim
\proclabel{ExistencePMCII}
\newsubhead{Extrinsic Curvature}[ExtrinsicCurvature]
Now let $K$ be extrinsic curvature, that is
$$
K(\lambda_1,...,\lambda_d) := (\lambda_1\cdot...\cdot\lambda_d)^{\frac{1}{d}}.
$$
For the rest of this section, we will suppose that $M$ is {\sl (1/4)-pinched} in the sense that
$$
\overline{\sigma}_\opMax < 4\overline{\sigma}_\opMin,\eqnum{\nexteqnno[QuarterPinched]}
$$
where $\overline{\sigma}_\opMax$ and $\overline{\sigma}_\opMin$ are respectively the maximum and minimum values of the sectional curvatures of planes tangent to $M$. We suppose, furthermore, that $M$ is {\sl pointwise $(1/2)$-pinched} in the sense that, for all $p\in M$,
$$
\overline{\sigma}_\opMax(p) <  2\overline{\sigma}_\opMin(p),\eqnum{\nexteqnno[PointwiseHalfPinched]}
$$
where $\overline{\sigma}_\opMax(p)$ and $\overline{\sigma}_\opMin(p)$ are respectively the maximum and minimum values of the sectional curvatures of planes tangent to $M$ at $p$.
\par
Each of \eqnref{QuarterPinched} and \eqnref{PointwiseHalfPinched} implies independently that $M$ is diffeomorphic to the standard sphere (c.f. \cite{Berger}). Furthermore, the condition of $(1/4)$-pinching is not actually necessary. Indeed, it is only imposed here in order to ensure embeddedness, which, in particular, implies that all immersions considered are prime. However, in forthcoming work, we will show how to handle the case of non-prime immersions, which will allow this condition to be discarded.
\par
Let $\Cal{I}$ be the space of unparametrised embeddings of $\Sigma$ into $M$ which bound some convex set. Let $\Cal{O}$ be the space of smooth, positive functions $f:M\rightarrow]0,\infty[$ such that, at every point $p$ of $M$,
$$
\frac{\|\overline{\nabla}^2f(p)\|}{f(p)} < 2\overline{\sigma}_\opMin(p) - \overline{\sigma}_\opMax(p).\eqnum{\nexteqnno[ConditionsOnFPEC]}
$$
Let $\Cal{Z}\subseteq\Cal{I}\times\Cal{O}$ be the solution space, that is
$$
\Cal{Z} := \left\{ ([e],f)\ |\ K_e = f\circ e\right\},\eqnum{\nexteqnno[SolutionSpacePEC]}
$$
and let $\Pi:\Cal{Z}\rightarrow\Cal{O}$ be the projection onto the second factor.
\par
Consider an LSC immersion, $e:\Sigma\rightarrow M$. In what follows, we will use the $K$-Laplacian, which is defined over $\Sigma$ as follows. For any twice differentiable function $\phi:\Sigma\rightarrow\Bbb{R}$ and for any point $p$ we have
$$
(\Delta^K\phi)(p) := \sum_{i=1}^d\frac{1}{\lambda_i}f_{;ii}(p),\eqnum{\nexteqnno[KLaplacian]}
$$
where $\lambda_1,...,\lambda_d$ are the principal curvatures of $e$ at this point.
\noskipproclaim{Lemma \nextprocno}
$$
\sum_{i=1}^d\frac{1}{\lambda_i} \geq \frac{d}{f}.
$$
\endproclaim
\proclabel{LowerBoundOnTrDKPEC}
\proof By concavity, for all $0<\lambda_1<...<\lambda_d$,
$$\eqalign{
K(1,...,1) &\leq K(\lambda_1,...,\lambda_d) + DK(\lambda_1,...,\lambda_d)\cdot(1-\lambda_1,...,1-\lambda_d),\cr
&=f + \frac{f}{d}\sum_{i=1}^d\frac{1}{\lambda_i}(1-\lambda_i),\cr
&=\frac{f}{d}\sum_{i=1}^d\frac{1}{\lambda_i}.\cr}
$$
The result follows.\qed
\proclaim{Lemma \nextprocno}
\noindent If $e$ has extrinsic curvature prescribed by $f$, then, at every point of $\Sigma$,
$$
\sum_{i=1}^d\frac{1}{\lambda_i}A_{ii;pp} \geq \frac{df_{:pp}}{f} - \frac{df_{:\nu}\lambda_p}{f},
$$
where here ``$:$'' denotes covariant differentiation over $M$.
\endproclaim
\proclabel{DeltaAPEC}
\proof Indeed, let $A$ be the shape operator of $e$. Differentiating the relation $f=\opDet(A)^{\frac{1}{d}}$ yields
$$
f_{;p} = \frac{f}{d}\opTr(A^{-1}A_{;p}),
$$
and differentiating a second time yields
$$
f_{;pp} = \frac{f}{d}\left(\frac{1}{d}\opTr(A^{-1}A_{;p})^2 - \opTr(A^{-1}A_{;p}A^{-1}A_{;p}) + \opTr(A^{-1}A_{;pp})\right).
$$
However, by the Cauchy-Schwarz inequality, for any symmetric matrix $M$,
$$
\opTr(M^2) \geq \frac{1}{d}\opTr(M)^2.
$$
More generally, for any symmetric matrix $M$, and any symmetric positive definite matrix $N$,
$$
\opTr(MNMN) = \opTr\left(\left(N^{\frac{1}{2}}MN^{\frac{1}{2}}\right)^2\right) \geq
\frac{1}{d}\opTr\left(N^{\frac{1}{2}}MN^{\frac{1}{2}}\right) = \frac{1}{d}\opTr(NM)^2.
$$
Thus, bearing in mind \eqnref{HessianOfRestriction},
$$
\sum_{i=1}^d\frac{1}{\lambda_i}A_{ii;pp} = \opTr(A^{-1}A_{;pp}) \geq \frac{df_{;pp}}{f} = \frac{df_{:pp}}{f} - \frac{d f_{:\nu}\lambda_p}{f},
$$
as desired.\qed
\proclaim{Lemma \nextprocno}
\noindent If the extrinsic curvature of $e$ is prescribed by $f$, then at any local maximum of $\lambda_d/f$,
$$\eqalign{
\Delta^K\left(\frac{\lambda_d}{f}\right)
&\geq f^2\left(\sum_{i=1}^d\frac{1}{\lambda_i}\right)
\left(\frac{c\left(2\overline{\sigma}_\opMin - \overline{\sigma}_\opMax\right)}{f^2} - \frac{c\|\overline{\nabla}^2f\|}{f^3} - \frac{2\|\overline{\nabla}\mR\|}{f^3}\right)\cr
&\qquad\qquad + fd\left(\frac{2\overline{\sigma}_\opMax-\overline{\sigma}_\opMin}{f^2} + \frac{\|\overline{\nabla}^2f\|}{f^3}\right),\cr}
$$
in the viscosity sense, where $c:=\lambda_d/f$ is the pinching factor.
\endproclaim
\proclabel{EstimateOfDeltaPinchingFactorPEC}
\proof Let $a_d$ be defined as in Section \subheadref{GeneralisedSimonsFormula}. In particular, a local maximum of $\lambda_d/f$ is also a local maximum of $a_d/f$. As in the proof of Lemma \procref{EstimateOfDeltaPinchingFactorPMC}, at any such point,
$$
\Delta^K\left(\frac{a_d}{f}\right) = \frac{1}{f}\Delta^K a_d - \frac{a_d}{f^2}\Delta^K f.
$$
Now, by Lemma \procref{SecondDerivativeOfA},
$$
\Delta^K a_d = \sum_{i=1}^d\frac{1}{\lambda_i}a_{d;ii} = \sum_{i=1}^d\frac{1}{\lambda_i}A_{dd;ii}.
$$
Thus, by the generalised Simons' formula,
$$\eqalign{
\Delta^K a_d
&=\sum_{i=1}^d\frac{1}{\lambda_i}A_{ii;dd} + \sum_{i=1}^d\frac{1}{\lambda_i}(\lambda_d-\lambda_i)(\lambda_d\lambda_i + \mR_{iddi})\cr
&\qquad +\sum_{i=1}^d\frac{1}{\lambda_i}\lambda_d(\mR_{iddi} - \mR_{i\nu\nu i})\cr
&\qquad -\sum_{i=1}^d\frac{1}{\lambda_i}\lambda_i(\mR_{iddi} - \mR_{d\nu\nu d})\cr
&\qquad +\sum_{i=1}^d\frac{1}{\lambda_i}(\mR_{id\nu i:d} + \mR_{id\nu d:i}).\cr}
$$
Bearing in mind \eqnref{HessianOfRestriction} and Lemma \procref{DeltaAPEC}, this yields
$$\eqalign{
\frac{1}{f}\Delta^K a_d &\geq
-\frac{f_{\nu}d\lambda_d}{f^2}+\frac{1}{f}\sum_{i=1}^d\lambda_d(\lambda_d-\lambda_i)\cr
&\qquad+f^2\left(\sum_{i=1}^d\frac{1}{\lambda_i}\right)
\left(\frac{(2\overline{\sigma}_\opMin-\overline{\sigma}_\opMax)\lambda_d}{f^3}
- \frac{2\|\overline{\nabla}\mR\|}{f^3}\right)\cr
&\qquad + fd\left(\frac{\overline{\sigma}_\opMin-2\overline{\sigma}_\opMax}{f^2}
- \frac{\|\overline{\nabla}^2f\|}{f^3}\right).\cr}
$$
Using \eqnref{HessianOfRestriction} again, we have
$$\eqalign{
\frac{a_d}{f^2}\Delta^K f &= \frac{\lambda_d}{f^2}\sum_{i=1}^d \frac{1}{\lambda_i}f_{;ii}\cr
&=\frac{\lambda_d}{f^2}\sum_{i=1}^d\frac{1}{\lambda_i}\left(f_{:ii} - f_{:\nu}\lambda_i\right)\cr
&\leq -\frac{df_{:\nu}\lambda_d}{f^2} +f^2\left(\sum_{i=1}^d\frac{1}{\lambda_i}\right)\frac{\|\overline{\nabla}^2f\|\lambda_d}{f^4},\cr}
$$
and the result follows upon combining these relations and substituting $c=\lambda_d/f$.\qed
\proclaim{Lemma \nextprocno}
\noindent For all $f\in\Cal{O}$, there exists $B>0$ such that if $([e],f)\in\Cal{Z}$, then
$$
\|A\| \leq B,
$$
where $A$ is the shape operator of $e$. Furthermore, $B$ can be chosen to vary continuously with $f$.
\endproclaim
\proclabel{UpperBoundOfShapeOperatorPEC}
\proof Let $x\in\Sigma$ maximise the pinching factor $c=\lambda_d/f$. By Lemmas \procref{LowerBoundOnTrDKPEC} and \procref{EstimateOfDeltaPinchingFactorPEC}, there exist $a,b>0$, which only depend on $f$ and $M$ such that, at this point
$$
\Delta^K\left(\frac{\lambda_d}{f}\right) \geq ac - b
$$
in the viscosity sense. It follows by the maximum principle that $c<b/a$, so that
$$
\|A\| \leq B := \frac{b\|f\|}{a},
$$
as desired.\qed
\proclaim{Lemma \nextprocno}
\noindent The projection $\Pi:\Cal{Z}\rightarrow\Cal{O}$ is a proper map.
\endproclaim
\proclabel{PropernessPEC}
\proof Let $(f_m)\in\Cal{O}$ be a sequence converging towards $f_\infty\in\Cal{O}$, and, for all $m$, let $[e_m]$ be an embedding such that $([e_m],f_m)\in\Cal{Z}$. By compactness, $M$ has sectional curvature bounded below by $\epsilon^2>0$. For all $m$, since $e_m$ is convex, the sectional curvature of $e_m^*g$ is also bounded below by $\epsilon^2$, and so, by the Bonnet-Myers theorem,
$$
\opDiam([e_m]) \leq \frac{\pi}{\epsilon}.
$$
By Lemma \procref{UpperBoundOfShapeOperatorPEC}, there exists $B>0$ such that, for all $m$,
$$
\|A_m\| \leq B,
$$
where $A_m$ here denotes the shape operator of $e_m$. It follows by the Arzela-Ascoli Theorem for immersed surfaces (c.f. \cite{SmiAAT}) and elliptic regularity that there exists an unparametrised immersion $[e_\infty]$ towards which $([e_m])$ subconverges. Taking limits, $[e_\infty]$ has extrinsic curvature prescribed by $f_\infty$. In particular, since $f_\infty>0$, $[e_\infty]$ is LSC. Finally, since $M$ is $(1/4)$-pinched, by \cite{EspRos}, $[e_\infty]$ is embedded and bounds a convex set. It follows that $([e_\infty],f_\infty)$ is also an element of $\Cal{Z}$, and this completes the proof.\qed
\proclaim{Lemma \nextprocno}
\noindent For all $f\in C^\infty(M)$, there exists $B,T>0$ such that if $t\geq T$ and if $([e],t(1+t^{-2}f))\in\Cal{Z}$, then
$$
\lambda_1(x) \geq Bf(x),
$$
at every point $x$ of $\Sigma$, where $\lambda_1$ here denotes the least principal curvature of $e$.
\endproclaim
\proclabel{LowerBoundOfShapeOperatorPEC}
\proof For all $t$, denote $f_t:=t(1+t^{-2}f)$, and suppose that $([e],f_t)$ is an element of $\Cal{Z}$. Fix $\epsilon>0$, and let $x\in\Sigma$ maximise the pinching factor $c=\lambda_d/f$ of $[e]$. By Lemmas \procref{LowerBoundOnTrDKPEC} and \procref{EstimateOfDeltaPinchingFactorPEC}, if $t$ is sufficiently large, then at this point,
$$
\Delta^K\left(\frac{\lambda_d}{f_t}\right)
\geq
cf_td\left(\frac{2\overline{\sigma}_\opMin - \overline{\sigma}_\opMax - \epsilon}{f_t^2}\right)
-f_td\left(\frac{2\overline{\sigma}_\opMax - \overline{\sigma}_\opMin + \epsilon}{f_t^2}\right)
$$
in the viscosity sense, so that, by the maximum principle,
$$
c \leq B_1 := \frac{2\overline{\sigma}_\opMax - \overline{\sigma}_\opMin + \epsilon}{2\overline{\sigma}_\opMin - \overline{\sigma}_\opMax - \epsilon}.
$$
Since $[e]$ has extrinsic curvature prescribed by $f$, its least principal curvature therefore satisfies at every point
$$
\lambda_1(x) \geq B_1^{1-d}f(x),
$$
as desired.\qed
\proclaim{Theorem \nextprocno}
\noindent For generic $f\in\Cal{O}$, the algebraic number of convex embedded hyperspheres in $M$ of extrinsic curvature prescribed by $f$ is equal to $-\chi(M)$.
\endproclaim
\proclabel{CountingPEC}
\remark Theorem \procref{MainTheoremExtrinsicCurvature} now follows via the discussion of Section \subheadref{VaryingTheMetric}.
\medskip
\proof Indeed, by Lemma \procref{PropernessPEC}, the degree theory developed in Section \headref{DegreeTheory} applies. It thus remains to calculate this degree. To this end, let $f\in C^\infty(M)$ be a function of Morse type, and for $t>0$, denote $f_t:=t(1+t^{-2}f)$. For sufficiently large $t$, $f_t\in\Cal{O}$. By Lemma \procref{LowerBoundOfShapeOperatorPEC}, there exists $B,T>0$ such that if $t\geq T$, and if $([e],f_t)$ is an element of $\Cal{Z}$, then the Ricci curvature of $e^*g$ satisfies
$$
\opRic_{ik} \geq Bt^2g_{ik}
$$
so that, by the Bonnet-Myers Theorem,
$$
\opDiam([e]) \leq \frac{\pi}{\sqrt{B}t}.
$$
It follows as in Section \subheadref{CalculatingTheDegree} that for sufficiently large $t$, $f_t$ is a regular value of $\Pi$, and the only elements of $\Pi^{-1}(\left\{f_t\right\})$ are those in the image of the canonical map, $Y_t$, constructed by Ye. In particular,
$$
\opDeg(\Pi) = -\chi(M),
$$
as desired.\qed
\newsubhead{Special lagrangian curvature}[SpecialLagrangianCurvature]
We now turn our attention to more sophisticated compactness results. In this section, $K$ will be special Lagrangian curvature, that is,
$$
K_\theta(\lambda_1,...,\lambda_d) := R_\theta(\lambda_1,...,\lambda_d),
$$
where $R_\theta$ is defined as in \cite{SmiSLC}, and $\theta$ is a real parameter taking values in the interval $[(d-1)\pi/2,d\pi/2[$. As in the preceding section, we will suppose that $M$ is $(1/4)$-pinched. We reiterate that this condition is only imposed in order to ensure embeddedness, and that in forthcoming work, by extending our degree theory to the case of non-prime immersions, we will show that it is not actually necessary.
\par
Let $\Cal{I}$ be the space of unparametrised immersions of $\Sigma$ into $M$ which bound some convex set. Let $\Cal{O}$ be the space of all smooth, positive functions $f:M\rightarrow]0,\infty[$. Let $\Cal{Z}\subseteq\Cal{I}\times[(d-1)\pi/2,d\pi/2[\times\Cal{O}$ be the solution space, that is,
$$
\Cal{Z}_\theta := \left\{ ([e],\theta,f)\ |\ K_{\theta,e} = f\circ e\right\},
$$
and let $\Pi:\Cal{Z}\rightarrow[(d-1)\pi/2,d\pi/2[\times\Cal{O}$ be the projection onto the second and third factors.
\proclaim{Lemma \nextprocno}
\noindent There exists $D>0$, which only depends on $M$, such that for all $([e],\theta,f)\in\Cal{Z}_\theta$,
$$
\opDiam([e]) \leq D.
$$
\endproclaim
\proclabel{DiameterBoundPSLC}
\proof Indeed, by compactness, $M$ has sectional curvature bounded below by $\epsilon^2>0$. Thus, since $e$ is convex, its sectional curvature is also bounded below by $\epsilon^2$, and the result now follows by the Bonnet-Myers theorem.\qed
\proclaim{Lemma \nextprocno}
\noindent For all $(\theta,f)\in[(d-1)\pi/2,d\pi/2[\times\Cal{O}$, there exists $B>0$ such that if $([e],\theta,f)\in\Cal{Z}_\theta$, then
$$
\|A\| \leq B,
$$
where $A$ is the shape operator of $e$. Furthermore, $B$ can be chosen to vary continuously with $f$.
\endproclaim
\proclabel{ShapeOperatorBoundPSLC}
\proof Suppose the contrary. There exists a sequence $(\theta_m,f_m)\in[(d-1)\pi/2,d\pi/2[\times\Cal{O}$ converging to $(\theta_\infty,f_\infty)\in [(d-1)\pi/2,d\pi/2[\times\Cal{O}$ such that, for all $m$, there is an element $([e_m],\theta_m,f_m)$ of $\Cal{Z}$ with
$$
\|A_m\| \geq m,
$$
where $A_m$ is the shape operator of $e_m$. Let $UM\subseteq TM$ be the bundle of unit tangent vectors over $M$. For all $m$, Let $x_m$ be the point of $\Sigma_m$ maximising the norm of $A_m$, let $N_m:\Sigma\rightarrow UM$ be the unit normal vector field over $e_m$ that is compatible with the orientation, and think of $N_m$ as an embedding of $\Sigma$ into $UM$ in its own right. By Theorem $1.4$ of \cite{SmiSLC}, there exists a complete, pointed, immersed submanifold $(\Sigma_\infty,N_\infty,x_\infty)$ of $UM$ towards which the sequence, $(\Sigma,N_m,x_m)$, of complete, pointed, immersed submanifolds subconverges in the $C^\infty$-Cheeger-Gromov sense.
\par
Denote $e_\infty:=\pi\circ N_\infty$, where $\pi:UM\rightarrow M$ is the canonical projection. By Theorem $1.3$ of \cite{SmiSLC}, either $e_\infty$ is an immersion at every point, or $N_\infty$ is a covering of the unit, normal sphere bundle, $\opN\Gamma$, of some complete geodesic $\Gamma\subseteq M$. Furthermore, the second case cannot occur unless $\theta_\infty=(d-1)\pi/2$. However, if $e_\infty$ were an immersion at $x_\infty$, then the sequence $(A_m(x_m))$ would converge towards the shape operator of $e_\infty$ at this point, and since we have assumed that this sequence is unbounded, we conclude that the second case holds.
\par
By Lemma \procref{DiameterBoundPSLC}, there exists $D$ such that $\opDiam([e_m])\leq D$ for all $M$. It follows that $\opLength([\Gamma])\leq D$, and, in particular, $\Gamma$ is closed. Since $[e_m]$ is embedded for all $m$, the covering order of $N_\infty$ over $\opN\Gamma$ is equal to $1$. That is, $N_\infty$ is a diffeomorphism, so that, in particular, $\Sigma_\infty$ is compact. Thus, by definition of pointed $C^\infty$-Cheeger-Gromov convergence, $\Sigma_\infty$ is diffeomorphic to $\Sigma$, which is absurd, and the result follows.\qed
\proclaim{Lemma \nextprocno}
\noindent The projection $\Pi:\Cal{Z}\rightarrow\Cal{O}$ is a proper map.
\endproclaim
\proclabel{PropernessPSLC}
\proof Let $(\theta_m,f_m)\in[(d-1)\pi/2,d\pi/2[\times\Cal{O}$ be a sequence converging towards $f_\infty\in\Cal{O}$, and, for all $m$, let $[e_m]$ be an immersion such that $([e_m],\theta_m,f_m)\in\Cal{Z}$. By Lemma \procref{DiameterBoundPSLC}, there exists $D>0$ such that, for all $m$,
$$
\opDiam([e_m])\leq D.
$$
By Lemma \procref{ShapeOperatorBoundPSLC}, there exists $B>0$ such that, for all $m$,
$$
\|A_m\| \leq B,
$$
where $A_m$ here denotes the shape operator of $e_m$. It follows by the Arzela-Ascoli theorem for immersed surfaces (c.f. \cite{SmiAAT}) and elliptic regularity that there exists an unparametrised immersion $[e_\infty]$ towards which $([e_m])$ subconverges. Taking limits, $[e_\infty]$ has $K_{\theta_\infty}$-curvature prescribed by $f_\infty$. In particular, since $f_\infty>0$, $[e_\infty]$ is LSC. Finally, since $M$ is $(1/4)$-pinched, by \cite{EspRos}, $[e_\infty]$ is embedded and bounds a convex set. It follows that $([e_\infty],\theta_\infty,f_\infty)$ is an element of $\Cal{Z}$, and this completes the proof.\qed
\proclaim{Theorem \nextprocno}
\noindent For all $\theta\in[(d-1)\pi/2,d\pi/2[$, and for generic $f\in\Cal{O}$, the algebraic number of convex, embedded hyperspheres in $M$ of $\theta$-special lagrangian prescribed by $f$ is equal to $-\chi(M)$.
\endproclaim
\proclabel{CountingPSLC}
\remark Theorem \procref{MainTheoremSpecialLagrangianCurvature} now follows via the discussion of Section \subheadref{VaryingTheMetric}.
\medskip
\proof Indeed, by Lemma \procref{PropernessPSLC}, the degree theory developed in Section \headref{DegreeTheory} applies. It thus remains to calculate the degree. To this end, fix $\theta\in](d-1)\pi/2,d\pi/2[$. Since $\theta>(d-1)\pi/2$, by Lemma $2.2$ of \cite{EspRos}, there exists $B>0$ such that if $([e],\theta,f)\in\Cal{Z}$, then
$$
\lambda_1(x) \geq (f\circ e)(x)/B,
$$
at every point of $\Sigma$. In particular, the Ricci curvature of $e^*g$ satisfies
$$
\opRic_{ik} \geq \left(\minf_{x\in M}f(x)\right)/B^2,
$$
so that, by the Bonnet-Myers theorem,
$$
\opDiam([e]) \leq B\pi/\left(\minf_{x\in M}f(x)\right).
$$
Now let $f\in C^\infty(M)$ be a function of Morse type, and for $t>0$, denote $f_t:=t(1+t^2f)$. For sufficiently large $t$, $f_t\in\Cal{O}$, and if $([e],\theta,f_t)\in\Cal{Z}$, then $\opDiam([e])\leq B\pi/2t$. It now follows as in Section \subheadref{CalculatingTheDegree} that for sufficiently large $t$, $f_t$ is a regular value of $\Pi$, and the only elements of $\Pi^{-1}(\left\{f_t\right\})$ are those in the image of the canonical map, $Y_t$, constructed by Ye. In particular,
$$
\opDeg(\Pi) = -\chi(M),
$$
as desired.\qed
\medskip
Since the special lagrangian curvature is not a widely known object of study, it is worth considering two special cases in order to correctly understand the significance of Theorem \procref{CountingPSLC}. First, when $d=2$, and $\theta=\pi/2$, the special lagrangian curvature, $K_{\pi/2}$, is none other than the extrinsic curvature, that is
$$
K_{\pi/2}(\lambda_1,\lambda_2) = (\lambda_1\lambda_2)^{1/2},
$$
and we therefore obtain the following refinement of Theorem \procref{CountingPEC} in the $2$-dimensional case.
\proclaim{Theorem \nextprocno}
\noindent Suppose that $g$ is a $(1/4)$-pinched metric over the sphere $S^3$. Then, for generic $f:M\rightarrow]0,\infty[$, the algebraic number of convex embedded spheres in $(S^3,g)$ of extrinsic curvature prescribed by $f$ is equal to $0$.
\endproclaim
\proclabel{CountingPSLCTWODIM}
When $d=3$, and $\theta=\pi$, the special lagrangian curvature $K_\pi$ is the curvature quotient,
$$
K_\pi(\lambda_1,\lambda_2,\lambda_3) = \left(\frac{3(\lambda_1\lambda_2\lambda_3)}{(\lambda_1+\lambda_2+\lambda_3)}\right)^{\frac{1}{2}},
$$
so that, denoting by $K:=\lambda_1\lambda_2\lambda_3$ the extrinsic curvature, and by $H:=\lambda_1+\lambda_2+\lambda_3$ the mean curvature, we obtain,
\proclaim{Theorem \nextprocno}
\noindent Suppose that $g$ is a $(1/4)$-pinched metric over the sphere $S^4$. Then, for generic $f:M\rightarrow]0,\infty[$, the algebraic number of convex embedded spheres $[e]$ in $(S^4,g)$ such that
$$
K_e = (f\circ e)H_e,
$$
is equal to $2$. In particular, there exist at least $2$ distinct, convex embedded spheres with this property.
\endproclaim
\proclabel{CountingPSLCTHREEDIM}
Finally, we recall that in each of Theorems \procref{CountingPSLCTWODIM} and \procref{CountingPSLCTHREEDIM}, the $(1/4)$-pinching condition is only necessary to ensure embeddedness, and in later work we aim to show how it may be discarded.
\newsubhead{Extrinsic curvature in two dimensions}[ExtrinsicCurvatureInTwoDimensions]
We conclude by improving Theorems \procref{CountingPEC} and \procref{CountingPSLC} for the case of surfaces inside $3$-dimensional manifolds. In what follows, we adopt the convention whereby the Ricci and scalar curvatures of an $n$-dimensional manifold are given by
$$\eqalign{
\opRic_{ik} &= -\frac{1}{n-1}\curvR[ijkj],\ \text{and}\cr
\opScal &= \frac{1}{n}\opRic_i^i.\cr}
$$
With these conventions, the Ricci and scalar curvatures of the standard $n$-dimensional sphere are equal to $\delta_{ik}$ and $1$ respectively.
\par
Let $K$ denote (the square root of) the extrinsic curvature. That is
$$
K(\lambda_1,\lambda_2) := (\lambda_1\lambda_2)^{\frac{1}{2}}.
$$
Let $\Cal{I}$ be the space of LSC unparametrised immersions. Let $\Cal{O}$ be the space of smooth, positive functions $f:M\rightarrow]0,\infty[$ such that, for all $x\in M$,
$$
\frac{\|Df(x)\|^2}{f(x)^4} + \frac{\|\mRic_0(x)\|_o^2}{2f(x)^4} + \frac{4\|Df(x)\|\|\mRic_0\|_o}{f(x)^4} < 1 + \frac{\overline{\sigma}_\opMin(x)}{f(x)^2},
$$
where $\overline{\sigma}_\opMin$ is defined as in Section \subheadref{ExtrinsicCurvature}, $\mRic_0$ is the trace-free Ricci curvature, and $\|\mRic_0\|_o$ is its operator norm when considered as an endomorphism of $TM$. Observe again that each of the summands above is scale invariant in the sense that it is unchanged when the metric is multiplied by a constant factor. Let $\Cal{Z}\subseteq\Cal{I}\times\Cal{O}$ be the solution space, that is,
$$
\Cal{Z} := \left\{([e],f)\ |\ K_e = f\circ e\right\},
$$
and let $\Pi:\Cal{Z}\rightarrow\Cal{O}$ be the projection onto the second factor.
\par
As usual, the degree is constructed by first proving that $\Pi$ is proper. In the present case, this result will be based on an estimate, valid in all dimensions, for the scalar curvature of the second fundamental form of the immersion, which, by local strict convexity, also defines a riemannian metric over $\Sigma$. Let $e:\Sigma\rightarrow M$ be an LSC immersion and let $A$ be its shape operator. Since $A$ is positive definite, in particular, it is invertible, and we denote its inverse by $B$, that is
$$
B^{ik}A_{kj} = \delta^i_j.
$$
Define the metric $g^A$ over $\Sigma$ by
$$
g^A(X,Y) := \langle AX,Y\rangle.
$$
Let $\nabla^A$ be its Levi-Civita covariant derivative. Let $\Omega$ be the relative Christoffel tensor of $\nabla$ with respect to $\nabla^A$, that is, for all $i$, $j$, and $k$,
$$
\Omega^k_{ij}\partial_k := \nabla_{\partial_i}\partial_j - \nabla^A_{\partial_i}\partial_j.
$$
\proclaim{Lemma \nextprocno}
\noindent For all $i$, $j$ and $k$
$$
\Omega^k_{ij} = \frac{1}{2}B^{kp}\mR_{pi\nu j} - \frac{1}{2}B^{kp}A_{pi;j},
\eqnum{\nexteqnno[ChristoffelSymbol]}
$$
where $\nu$ here denotes the outward pointing normal direction over $e$, and the subscript ``$;$'' here denotes covariant differentiation with respect to the intrinsic metric of $\Sigma$.
\endproclaim
\proclabel{ChristoffelSymbol}
\proof Observe that if $\hat{\Omega}$ denotes the relative Christoffel symbol of $\nabla^A$ with respect to $\nabla$, then $\hat{\Omega}=-\Omega$. Thus, by the Koszul formula,
$$
\Omega^k_{ij} = -\hat{\Omega}^k_{ij} = -\frac{1}{2}B^{kp}(A_{pi;j} + A_{pj;i} - A_{ij;p}),
$$
so that, by \eqnref{CodazziMainardi},
$$\eqalign{
\hat{\Omega}^k_{ij} &= -\frac{1}{2}B^{kp}(A_{pi;j} + \mR_{ij\nu p} - \mR_{pj\nu i})\cr
&=-\frac{1}{2}B^{kp}(A_{pi;j} + \mR_{ij\nu p} + \mR_{jp\nu i}).\cr}
$$
The result now follows by the first Bianchi identity.\qed
\medskip
Let $\opR^A$, $\opRic^A$ and $\opScal^A$ denote respectively the Riemann curvature tensor, the Ricci curvature tensor and the scalar curvature of the metric $g^A$, using the convention indicated at the beginning of this section.
\proclaim{Lemma \nextprocno}
\noindent Suppose that $\opDet(A)=f^2$, where $f:M\rightarrow]0,\infty[$ is a smooth, positive function. There exists a $1$-form, $\alpha$, over $\Sigma$ such that, for all $\lambda>0$,
$$
\opScal^A \geq \frac{1}{d}B^{ik}\opRic_{ik} + \nabla^A\cdot\alpha - \frac{(1+\lambda)}{f^2d(d-1)}\|Df\|_A^2 - \frac{(1+d\lambda^{-1})}{4d(d-1)}\|\mR_{\cdot\cdot\nu\cdot}\|^2_A,\eqnum{\nexteqnno[ScalarCurvature]}
$$
where $\nabla^A\cdot$ here denotes the divergence operator of $\nabla^A$, and
$$\eqalign{
\|Df\|_A^2 &:= B^{ij}f_{;i}f_{;j},\ \text{and} \cr
\|\mR_{\cdot\cdot\nu\cdot}\|_A^2 &:= B^{ip}B^{jq}B^{kr}\mR_{ij\nu k}\mR_{pq\nu r}.\cr}
$$
\endproclaim
\proof Recall that we use the subscript ``$;$'' to denote covariant differentiation with respect to $\nabla$. In what follows, we will use the subscript ``,'' to denote covariant differentiation with respect to $\nabla^A$. The Riemann curvature tensors of $g$ and $g^A$ are related by
$$
\curvR[ijkl] = \curvAR[ijkl] + \Omega^l_{jk,i} - \Omega^l_{ik,j} + \Omega^l_{im}\Omega^m_{jk} - \Omega^l_{jm}\Omega^m_{ik}.
$$
Contracting this relation yields
$$
B^{ik}\opRic_{ik} = \opScal^A + \frac{1}{d(d-1)}B^{ik}\left(\Omega^p_{pk,i}-\Omega^p_{ik,p}+\Omega^p_{iq}\Omega^q_{pk} - \Omega^p_{pq}\Omega^q_{ik}\right).
$$
We now show that the first and second terms in parentheses on the right hand side combine to yield an exact form. First, differentiating the relation $\opDet(A)=f^2$ yields
$$
B^{ij}A_{ij;k} = \frac{2}{f}f_{;k},\eqnum{\nexteqnno[TraceOfDerivativeOfShapeOperator]}
$$
so that, by \eqnref{ChristoffelSymbol},
$$
\Omega^p_{pk} = \frac{1}{2}B^{pq}(\mR_{qp\nu k} - A_{qp;k}) = - \frac{1}{f}f_{;k} = - \frac{1}{f}f_{,k}.
\eqnum{\nexteqnno[FirstTerm]}
$$
Likewise, bearing in mind \eqnref{CodazziMainardi},
$$\eqalignno{
B^{ik}\Omega^p_{ik}
&=\frac{1}{2}B^{ik}B^{pq}(\mR_{qi\nu k} - A_{qi;k})&\cr
&=-\frac{1}{2}B^{ik}B^{pq}(A_{ki;q} + \mR_{kq\nu i} + \mR_{iq\nu k})&\cr
&=-\frac{1}{f}B^{pq}f_{;q} - B^{ik}B^{pq}\mR_{iq\nu k}.&\nexteqnno[SecondTerm]\cr}
$$
Thus, denoting
$$
\alpha^m := \frac{1}{d(d-1)}B^{ik}B^{mn}\mR_{in\nu k}.
$$
We have,
$$
B^{mk}\Omega^p_{pk} - B^{ik}\Omega^m_{ik} = d(d-1)\alpha^m,
$$
and taking the divergence yields,
$$
B^{ik}\Omega^p_{pk,i} - B^{ik}\Omega^p_{ik,p} = d(d-1)\nabla^A\cdot\alpha,
$$
as asserted.
\par
We now consider the last two terms in parentheses on the right hand side. First, \eqnref{FirstTerm} and \eqnref{SecondTerm} together yield
$$\eqalign{
B^{ik}\Omega^p_{pq}\Omega^q_{ik}
&= \frac{1}{f^2}B^{qr}f_{,q}(f_{,r} + fd(d-1)A_{rs}\alpha^s),\cr
&= \frac{1}{f}\|Df\|_A^2 + \frac{d(d-1)}{f}f_{,p}\alpha^p.\cr}
$$
Likewise, using \eqnref{CodazziMainardi} again,
$$\eqalign{
B^{ik}\Omega^p_{iq}\Omega^q_{pk}
&= B^{ik}\Omega^p_{pr}\Omega^q_{pk}\cr
&= \frac{1}{4}B^{ik}B^{pr}B^{qs}(A_{rq;i} - \mR_{rq\nu i})(A_{sp;k} - \mR_{sp\nu k})\cr
&= \frac{1}{4}B^{ip}B^{jq}B^{kr}(A_{ij;k} - \mR_{ij\nu k})(A_{pq;r} + \mR_{pq\nu r})\cr
&= \frac{1}{4}\|\nabla A\|^2_A - \frac{1}{4}\|\mR_{\cdot\cdot\nu\cdot}\|_A^2.\cr}
$$
Combining these terms yields
$$\eqalign{
\opScal^A
&=\frac{1}{d}B^{ik}\opRic_{ik} + \nabla^A\cdot\alpha + \frac{1}{4d(d-1)}\|\nabla A\|^2_A\cr
&\qquad - \frac{1}{4d(d-1)}\|\mR_{\cdot\cdot\nu\cdot}\|_A^2 - \frac{1}{f} f_{,p}\alpha^p - \frac{1}{f^2 d(d-1)}\|Df\|_A^2.\cr}
$$
However, using the Cauchy-Schwarz inequality, we obtain
$$\eqalign{
\|\alpha\|^2_A &= A_{ij}\alpha^i\alpha^j\cr
&=\frac{1}{d^2(d-1)^2}B^{ij}B^{mn}B^{pq}\mR_{mi\nu n}\mR_{pj\nu q}\cr
&\leq\frac{1}{d(d-1)^2}B^{ij}B^{mn}B^{pq}\mR_{mi\nu p}\mR_{nj\nu q}\cr
&=\frac{1}{d(d-1)^2}\|\mR_{\cdot\cdot\nu\cdot}\|_A^2,\cr}
$$
and applying the Cauchy-Schwarz inequality a second time, we obtain
$$\eqalign{
\|f_{,p}\alpha^p\|
&\leq \frac{d(d-1)\lambda^{-1}}{4\lambda}\|\alpha\|_A^2 + \frac{\lambda}{f^2d(d-1)}\|Df\|_A^2\cr
&=\frac{d\lambda^{-1}}{4d(d-1)}\|\mR_{\cdot\cdot\nu\cdot}\|_A^2 + \frac{\lambda}{f^2d(d-1)}\|Df\|_A^2,\cr}
$$
and the result follows.\qed
\medskip
We now restrict attention to the $2$ dimensional case.
\proclaim{Lemma \nextprocno}
\noindent Suppose that $d=2$ and that $\opDet(A)=f^2$, where $f:M\rightarrow]0,\infty[$ is a smooth, positive function. Then, for all $\lambda:M\rightarrow]0,\infty[$,
$$
\int_\Sigma\left(1
+\frac{\overline{\sigma}_\opMin}{f^2}
-\frac{(1+\lambda)\|Df\|^2}{f^4}
-\frac{(1+2\lambda^{-1})\|\opRic_0\|_o^2}{2f^4}\right)fH\opdVol \leq 4\pi.
\eqnum{\nexteqnno[IntegralInequalityPEC]}
$$
\endproclaim
\proof When $d=2$, \eqnref{ScalarCurvature} yields, for all $\lambda>0$,
$$
\frac{1}{2}B^{ik}\opRic_{ik}
-\frac{(1+\lambda)}{2f^2}\|Df\|_A^2 - \frac{(1+2\lambda^{-1})}{8}\|\mR_{\cdot\cdot\nu\cdot}\|^2_A
\leq\opScal^A - \nabla^A\cdot\alpha.
$$
However, since $\Sigma$ is $2$-dimensional,
$$
\opRic_{ik} = \opScal\delta_{ik},
$$
so that
$$\eqalign{
B^{ik}\opRic_{ik} &= \left(\frac{1}{\lambda_1} + \frac{1}{\lambda_2}\right)\opScal\cr
&= \frac{2H}{f^2}\opScal\cr
&\geq 2H + \frac{2H\overline{\sigma}_\opMin}{f^2}.\cr}
$$
Next,
$$\eqalign{
\|Df\|_A^2 &= \frac{1}{\lambda_1}f_1^2 + \frac{1}{\lambda_2}f_2^2,\cr
&\leq \frac{2H}{f^2}\|Df\|^2.\cr}
$$
Likewise,
$$\eqalign{
\|\mR_{\cdot\cdot\nu\cdot}\|_A^2
&=\frac{2}{\lambda_1\lambda_2}\left(\frac{1}{\lambda_1}\mR^2_{12\nu1} + \frac{1}{\lambda_2}\mR^2_{21\nu2}\right)\cr
&=\frac{2}{f^2}\left(\frac{1}{\lambda_1}\mRic_{2\nu}^2 + \frac{1}{\lambda_2}\mRic_{1\nu}^2\right)\cr
&\leq\frac{4H}{f^4}\|\mRic_0\|_o^2.\cr}
$$
Thus,
$$
\left(1 + \frac{\overline{\sigma}_\opMin}{f^2} - \frac{(1+\lambda)\|Df\|^2}{f^4}
-\frac{(1+2\lambda^{-1})\|\mRic_0\|_o^2}{2f^4}\right)H
\leq \opScal^A - \nabla^A\cdot\alpha.
$$
The result now follows by multiplying by $\opdVol^A = f\opdVol$ and applying both the Gauss-Bonnet theorem and the divergence theorem.\qed
\proclaim{Lemma \nextprocno}
\noindent For all $f\in\Cal{O}$, there exists $B>0$ such that if $([e],f)\in\Cal{Z}$, then
$$
\int_\Sigma H\opdVol \leq B.
$$
Furthermore, $B$ can be chosen to vary continuously with $f$.
\endproclaim
\proclabel{BoundedIntegralMeanCurvaturePECTwoDim}
\proof It suffices to show that for all such $f$, there exists a positive function $\lambda:M\rightarrow]0,\infty[$ such that the coefficient of $fH$ in the integral on the left-hand side of \eqnref{IntegralInequalityPEC} is strictly positive at every point. However, consider the function,
$$\eqalign{
P(\lambda,x) &:= \lambda^2\|Df(x)\|^2  + \|\mRic_0(x)\|_o^2\cr
&\qquad -\left(f(x)^4 + f(x)^2\overline{\sigma}_\opMin(x) - \|Df(x)\|^2 - \frac{1}{2}\|\mRic_0(x)\|^2_o\right)\lambda.\cr}
$$
The condition on $f$ ensures that this polynomial in $\lambda$ has two positive roots at every point, $x$, of $M$, and the result follows.\qed
\proclaim{Lemma \nextprocno}
\noindent For all $f\in\Cal{O}$, there exists $D>0$ such that if $([e],f)\in\Cal{Z}$, then
$$
\opDiam([e]) \leq D.
$$
Furthermore, $B$ can be chosen to vary continuously with $f$.
\endproclaim
\proclabel{BoundedDiameterPECTwoDim}
\proof Indeed, since $f\in\Cal{O}$, there exists $\epsilon>0$ such that
$$
\opScal \geq \overline{\sigma}_\opMin + f^2 > \epsilon,
$$
and the result follows by the Bonnet-Myers theorem.\qed
\proclaim{Lemma \nextprocno}
\noindent When $d=2$, the projection $\Pi:\Cal{Z}\rightarrow\Cal{O}$ is a proper map.
\endproclaim
\proclabel{PropernessPECTwoDim}
\proof Indeed, this follows immediately from Lemmas \procref{BoundedIntegralMeanCurvaturePECTwoDim} and \procref{BoundedDiameterPECTwoDim} and Theorem $2.7$ of \cite{LabA}.\qed
\proclaim{Theorem \nextprocno}
\noindent Let $M$ be a compact, oriented $3$-dimensional manifold. Then, for generic $f\in\Cal{O}$, the algebraic number of LSC immersed spheres in $M$ of extrinsic curvature prescribed by $f$ is equal to $0$.
\endproclaim
\proclabel{CountingPECTwoDimI}
\remark Theorem \procref{MainTheoremExtrinsicCurvatureTwoD} now follows via the discussion of Section \subheadref{VaryingTheMetric}.
\medskip
\proof Indeed, by Lemma \procref{PropernessPECTwoDim}, the degree theory of Section \headref{DegreeTheory} applies. It thus remains to calculate the degree. To this end, let $f\in C^\infty(M)$ be a function of Morse type, and for $t>0$, denote $f_t:=t(1+ t^{-2}f)$. For sufficiently large $t$, $f_t\in\Cal{O}$. Furthermore, for sufficiently large $t$, if $([e],f_t)$ is an element of $\Cal{Z}$, then the scalar curvature of $e^*g$ satisfies
$$
\opScal \geq \frac{t^2}{2}
$$
at every point of $\Sigma$, so that, by the Bonnet-Myers theorem,
$$
\opDiam([e]) \leq \frac{\sqrt{2}\pi}{t}.
$$
It follows as in Section \subheadref{CalculatingTheDegree} that for sufficiently large $t$, $f_t$ is a regular value of $\Pi$, and the only elements of $\Pi^{-1}(\left\{f_t\right\})$ are those in the image of the canonical map, $Y_t$, constructed by Ye. In particular,
$$
\opDeg(\Pi) = -\chi(M),
$$
as desired.\qed
\medskip
\noindent Finally, as outlined in Section \subheadref{VaryingTheMetric}, Theorem \procref{DefinitionOfTheDegree} is also valid for generic metrics, given a fixed prescribing function. In particular, we obtain the following alternative version of Theorem \procref{CountingPECTwoDimI}.
\proclaim{Theorem \nextprocno}
\noindent Let $M$ be a compact, oriented, $3$-dimensional manifold. For a generic riemannian metric $g$ over $M$, and for all $k$ such that
$$
k > \frac{1}{2}\left(\sqrt{\overline{\sigma}_\opMin^2 + \|\mRic_0\|_o^2} - \overline{\sigma}_\opMin\right),
$$
the algebraic number of LSC immersed spheres in $M$ of constant extrinsic curvature equal to $k$ is itself equal to $0$.
\endproclaim
\inappendicestrue
\global\headno=0
\medskip
\goodbreak
\newhead{Weakly smooth maps}[WeaklySmoothMaps]
Different types of manifolds are determined using categories. Indeed, consider a category whose objects are open subsets of some vector space and whose morphisms are continuous maps. For example, the objects could be open subsets of $\Bbb{R}^d$, for some fixed $d$, and the morphisms could be smooth maps. Alternatively, the objects could be open subsets of $\Bbb{C}$, and the morphisms could be holomorphic maps. A {\sl manifold} in this category is defined to be a Hausdorff topological space, $X$, furnished with an atlas whose charts are objects of this category and whose transition maps are morphisms. For example, the first category described above defines smooth, $d$-dimensional manifolds, whilst the second defines Riemann surfaces. In finite dimensions, it is usual to assume in addition that $X$ is second countable. However, this is not always necessary, as in the case of Riemann surfaces, for example, and in infinite dimensions, it is not always viable.
\par
Consider a compact, finite dimensional manifold, $X$. Let $\opWS(X)$ denote the category whose objects are open subsets of $C^\infty(X)$ and whose morphisms are weakly smooth maps (c.f. Section \subheadref{TheManifoldOfImmersions}). A {\sl weakly smooth manifold} modelled on $C^\infty(X)$ is then defined to be a manifold in $\opWS(X)$. In particular, by taking $X$ to be a zero dimensional manifold consisting of a finite number of points, we see that all finite dimensional, smooth manifolds belong to this class.
\par
The remainder of this appendix is devoted to reviewing the basic properties of weakly smooth manifolds. Broadly speaking, the differential geometry of manifolds rests on two principles. The first is the possibility of constructing tangent bundles and of differentiating smooth maps, and the second is the inverse function theorem. We will show that although the weakly smooth category does not have an inverse function theorem, it nonetheless possesses all the requisite properties for the construction of tangent bundles and the differentiation of smooth maps.
\par
We first review how derivatives are defined in $\opWS(X)$. First, consider a compact, finite dimensional manifold, $M$, and a strongly smooth map, $\alpha:M\rightarrow C^\infty(X)$. Recall that the smooth map, $\tilde{\alpha}:M\times X\rightarrow\Bbb{R}$, is given by
$$
\tilde{\alpha}(p,x) := \alpha(p)(x).
$$
For any $p\in M$, the derivative of $\alpha$ at $p$ is given by
$$
D\alpha(p)\cdot\xi_x := \partial_r\tilde{\alpha}(c(r),\cdot)|_{r=0},
$$
where $\xi_x$ is a tangent vector to $M$ at $x$ and $c:]-\epsilon,\epsilon[\rightarrow M$ is a smooth curve such that $\partial_rc(0)=\xi_x$. Consider next an open subset, $\Cal{U}$, of $C^\infty(X)$, another compact, finite dimensional manifold, $X'$, and a weakly smooth map, $\Phi:\Cal{U}\rightarrow C^\infty(X')$. Given $f\in\Cal{U}$, and $g\in C^\infty(X)$, which we consider as a tangent vector to $\Cal{U}$ at $f$, define the strongly smooth map, $\alpha:]-\epsilon,\epsilon[\rightarrow\Cal{U}$, by $\alpha(r):=f+rg$. Since $\Phi$ is weakly smooth, the map $\beta:=\Phi\circ\alpha$ is also strongly smooth, and the derivative of $\Phi$ at $f$ in the direction of $g$ is defined by
$$
D\Phi(f)\cdot g := D\beta(0)\cdot\partial_r.
$$
\par
It follows in the usual manner from the definition that $D\Phi(f)$ maps $C^\infty(X)$ linearly into $C^\infty(X')$. The chain rule, however, is more subtle. There are two cases to be considered.
\proclaim{Lemma \nextprocno, {\bf Chain rule I}}
\noindent Given a compact, finite dimensional manifold, $M$, and a strongly smooth map $F:M\rightarrow\Cal{U}$, the derivative of $\Phi\circ F$ at any point, $p\in M$, satisfies
$$
D(\Phi\circ F)(p) = D\Phi(f)\circ DF(p),
$$
where $f:=F(p)$.
\endproclaim
\proclabel{ChainRuleI}
\proof Let $\xi_p$ be a tangent vector to $M$ at $p$. Let $c:]-\epsilon,\epsilon[\rightarrow M$ be a smooth curve such that $\partial_r c(0)=\xi_p$. Define the strongly smooth map, $\alpha_1:]-\epsilon,\epsilon[\rightarrow C^\infty(X)$, by $\alpha_1 := F\circ c$. Let $g:=D\alpha_1(0)\cdot\partial_r$ and define the strongly smooth map, $\alpha_2:]-\epsilon,\epsilon[\rightarrow C^\infty(X)$, by $\alpha_2(r):=f+rg$. There exists a smooth function, $\tilde{h}:]-\epsilon,\epsilon[\times X\rightarrow\Bbb{R}$, such that
$$
\tilde{\alpha}_2(r,x) = \tilde{\alpha}_1(r,x) + r^2\tilde{h}(r,x).
$$
Define the strongly smooth map, $\alpha:]-\epsilon,\epsilon[\times]-\epsilon^2,\epsilon^2[\rightarrow C^\infty(X)$, by
$$
\alpha(s,t) := \alpha_1(s) + th(s),
$$
where $h:]-\epsilon,\epsilon[\rightarrow C^\infty(X)$ is the strongly smooth map given by $h(r):=\tilde{h}(r,\cdot)$. In particular, $\alpha_1(r)=\alpha(r,0)$ and $\alpha_2(r)=\alpha(r,r^2)$. Since $\Phi$ is weakly smooth, the maps $\beta:=\Phi\circ\alpha$, $\beta_1:=\Phi\circ\alpha_1$ and $\beta_2:=\Phi\circ\alpha_2$ are all also strongly smooth. Furthermore, $\beta_1(r)=\beta(r,0)$ and $\beta_2(r)=\beta(r,r^2)$. However, by definition, $D\Phi(f)\cdot g=D\beta_2(0)\cdot\partial_r$, and so
$$\eqalign{
D\Phi(f)\cdot g &= D\beta(0)\cdot\partial_s\cr
&= D\beta_1(0)\cdot\partial_r\cr
&= D(\Phi\circ F)(p)\cdot\xi_p,\cr}
$$
and since $g=DF(p)\cdot\xi_p$, the result follows.\qed
\medskip
\noindent Consider now an open subset, $\Cal{V}$, of $C^\infty(X')$, another compact, finite dimensional manifold, $X''$, and a weakly smooth map, $\Psi:\Cal{V}\rightarrow C^\infty(X'')$.
\proclaim{Lemma \nextprocno, {\bf Chain rule II}}
\noindent If $\Phi(\Cal{U})\subseteq\Cal{V}$, then the derivative of $\Psi\circ\Phi$ at any point, $f\in\Cal{U}$, satisfies,
$$
D(\Psi\circ\Phi)(f) = D\Psi(f')\circ D\Phi(f),
$$
where $f':=\Phi(f)$.
\endproclaim
\proclabel{ChainRuleII}
\proof Let $g$ be an element of $C^\infty(X)$, which we consider as a tangent vector to $\Cal{U}$ at $f$. Denote $h:=D\Phi(f)\cdot g$. Define the strongly smooth map, $\alpha:]-\epsilon,\epsilon[\rightarrow\Cal{U}$, by $\alpha(r):=f + rg$. Since $\Phi$ is weakly smooth, the map, $\beta:=\Phi\circ\alpha$, is also strongly smooth. Now by definition,
$$
D(\Psi\circ\Phi)(f)\cdot g = D(\Psi\circ\Phi\circ\alpha)(0)\cdot\partial_r = D(\Psi\circ\beta)(0)\cdot\partial_r,
$$
and
$$
h = D\Phi(f)\cdot g = D\beta(0)\cdot\partial_r.
$$
However, by Lemma \procref{ChainRuleI},
$$
D(\Psi\circ\beta)(0)\cdot\partial_r = (D\Psi(f')\circ D\beta(0))\cdot\partial_r = D\Psi(f')\cdot h,
$$
and the result follows upon combining these relations.\qed
\medskip
These results allow us to construct the tangent bundle of any weakly smooth manifold modelled on $C^\infty(X)$. Indeed, given an open subset, $\Cal{U}$, of $C^\infty(X)$, define $\Cal{TU}:=\Cal{U}\times C^\infty(X)$, and given another open subset, $\Cal{V}$, of $C^\infty(X)$ and a weakly smooth map, $\Phi:\Cal{U}\rightarrow\Cal{V}$, define $\Cal{T}\Phi:\Cal{TU}\rightarrow\Cal{TV}$ by
$$
\Cal{T}\Phi(f,g) := (\Phi(f),D\Phi(f)\cdot g).
$$
Trivially, $\Cal{T}\opId=\opId$. Furthermore,
\proclaim{Lemma \nextprocno}
\noindent For all weakly smooth maps, $\Phi:\Cal{U}\rightarrow\Cal{V}$ and $\Psi:\Cal{V}\rightarrow\Cal{W}$, we have
$$
\Cal{T}(\Psi\circ\Phi)=\Cal{T}\Psi\circ\Cal{T}\Phi.
$$
\endproclaim
\proof Indeed, using Lemma \procref{ChainRuleII}, for all $(f,g)\in\Cal{TU}$,
$$\eqalign{
\Cal{T}(\Psi\circ\Phi)(f,g)
&=((\Psi\circ\Phi)(f),D(\Psi\circ\Phi)(f)\cdot g)\cr
&=(\Psi(\Phi(f)),D\Psi(\Phi(f))\cdot D\Phi(f)\cdot g)\cr
&=\Cal{T}\Psi(\Phi(f),D\Phi(f)\cdot g)\cr
&=(\Cal{T}\Psi\circ\Cal{T}\Phi)(f,g),\cr}
$$
as desired.\qed
\medskip
\noindent Finally, observing that $C^\infty(X)\times C^\infty(X)=C^\infty(X\sqcup X)$, we have
\proclaim{Lemma \nextprocno}
\noindent For any weakly smooth map, $\Phi:\Cal{U}\rightarrow\Cal{V}$, $\Cal{T}\Phi$ defines a weakly smooth map from $\Cal{TU}$ into $\Cal{TV}$.
\endproclaim
\proof Let $M$ be a compact, finite dimensional manifold and let $\alpha:=(\alpha_1,\alpha_2):M\rightarrow\Cal{TU}$ be a strongly smooth map. Define the strongly smooth map, $\beta:]-\epsilon,\epsilon[\times M\rightarrow\Cal{U}$ by $\beta(t,p):=\alpha_1(p) + t\alpha_2(p)$. Since $\Phi$ is weakly smooth, $\gamma:=\Phi\circ\beta$ is also strongly smooth. Now define the strongly smooth maps, $\delta_1,\delta_2:M\rightarrow C^\infty(X)$, such that
$$\eqalign{
\tilde{\delta}_1(p,x) &= \tilde{\gamma}(0,p,x),\ \text{and}\cr
\tilde{\delta}_2(p,x) &= \partial_t\tilde{\gamma}(0,p,x).\cr}
$$
In particular, $\delta:=(\delta_1,\delta_2)$ is also strongly smooth. However, by definition, $\delta=\Cal{T}\Phi\circ\alpha$. Since $\alpha$ is arbitrary, and since, in addition, $\delta$ varies continuously with $\alpha$, it now follows that $\Cal{T}\Phi$ is weakly smooth, as desired.\qed
\medskip
Given another compact, finite dimensional manifold, $Y$, we now introduce the category $\text{BWS}(X,Y)$ as follows. Its objects are triplets of the form $(\Cal{U},\Cal{U}\times C^\infty(Y),\Pi)$, where $\Cal{U}$ is an open subset of $C^\infty(X)$, and $\Pi:\Cal{U}\times C^\infty(Y)\rightarrow\Cal{U}$ is the canonical projection. Its morphisms are pairs of the form $(\Phi,\Psi)$, where $\Phi:\Cal{U}\rightarrow\Cal{V}$ and $\Psi:\Cal{U}\times C^\infty(Y)\rightarrow\Cal{V}\times C^\infty(Y)$ are weakly smooth maps such that $\Phi\circ\Pi=\Pi\circ\Psi$ and $\Psi$ is linear over each fibre. We leave the reader to review how this category serves to define weakly smooth vector bundles over weakly smooth manifolds in the same way that the category $\opWS(X)$ serves to define weakly smooth manifolds.
\par
The above lemmas show that the operator, $\Cal{T}$, defines a functor from $\opWS(X)$ into $\text{BWS}(X,X)$. From this, it immediately follows that every weakly smooth manifold, $\Cal{M}$, modelled on $C^\infty(X)$ has a well defined tangent bundle, whose typical fibre is $C^\infty(X)$, and whose total space is a weakly smooth manifold modelled on $C^\infty(X)\times C^\infty(X)=C^\infty(X\sqcup X)$. Furthermore, every weakly smooth map, $\Phi:\Cal{M}\rightarrow\Cal{N}$, between weakly smooth manifolds, has a well defined derivative, $\Cal{T}\Phi:\Cal{TM}\rightarrow\Cal{TN}$, which defines a weakly smooth map from the total space of $\Cal{TM}$ into the total space of $\Cal{TN}$ whose restriction to every fibre is linear. In summary, weakly smooth manifolds possess all the requisite properties for weakly smooth maps between them to be studied in terms of their derivatives.
\newhead{Prime immersions}[PrimeImmersions]
Consider a smooth, simply connected, $d$-dimensional manifold, $\tilde{\Sigma}$. Let $\Gamma$ by a group of diffeomorphisms acting on $\tilde{\Sigma}$ properly discontinuously and cocompactly. Let $i:\tilde{\Sigma}\rightarrow M$ be a $\Gamma$-invariant immersion, and suppose that $i$ is {\sl prime} in the sense that if $\alpha:\tilde{\Sigma}\rightarrow\tilde{\Sigma}$ is a smooth diffeomorphism such that $i=i\circ\alpha$, then $\alpha\in\Gamma$.
\proclaim{Lemma \nextprocno}
\noindent If $\alpha:\tilde{\Sigma}\rightarrow\tilde{\Sigma}$ is a continuous map such that $i=i\circ\alpha$, then $\alpha$ is a smooth diffeomorphism.
\endproclaim
\proclabel{ConditionForDiffeomorphism}
\proof Furnish $\tilde{\Sigma}$ with the metric induced by $i$. In particular, this makes $\alpha$ into a local isometry. Furthermore, since $\Gamma$ is cocompact, $\tilde{\Sigma}$ is complete, and so $\alpha$ is a covering map. Finally, since $\tilde{\Sigma}$ is simply connected, it follows that $\alpha$ is a smooth diffeomorphism, and this completes the proof.\qed
\medskip
Consider the Cartesian product, $\tilde{\Sigma}\times]-\epsilon,\epsilon[$ and let $\pi:\tilde{\Sigma}\times]-\epsilon,\epsilon[\rightarrow\tilde{\Sigma}$ be the projection onto the first factor. As in Section \subheadref{TheManifoldOfImmersions}, define $\Cal{E}:\tilde{\Sigma}\times]-\epsilon,\epsilon[\rightarrow M$ by
$$
\Cal{E}(x,t) := \opExp_{i(x)}(tN_i(x)),
$$
where $\opExp$ is the exponential map of $M$ and $N_i:\tilde{\Sigma}\rightarrow TM$ is the unit, normal vector field over $i$ which is compatible with the orientation. By compactness, we may choose $\epsilon$ such that $\Cal{E}$ is an immersion. Let $C^\infty_\Gamma(\tilde{\Sigma})$ denote the space of smooth, $\Gamma$-invariant functions over $\Sigma$, and for $f\in C^\infty_\Gamma(\tilde{\Sigma})$, define $\hat{f}(x):=\Cal{E}(x,f(x))$. Observe that for $\|f\|_{L^\infty}<\epsilon$, $\hat{f}$ is also an immersion. We now prove that prime immersions are $C^0$-stable in the following sense.
\proclaim{Theorem \nextprocno}
\noindent There exists $0<\epsilon'<\epsilon$ such that if $\|f\|_{L^\infty}<\epsilon'$, then $\hat{f}$ is a prime immersion.
\endproclaim
\proclabel{Stability}
Upon reducing $\epsilon$ further if necessary, we may suppose that for all $p\in\tilde{\Sigma}$, the restriction of $\Cal{E}$ to $B_\epsilon(p)\times]-\epsilon,\epsilon[$ is a diffeomorphism onto its image, and we denote the inverse of this restriction by $\Cal{E}_p^{-1}$.
\proclaim{Lemma \nextprocno}
\noindent If $f\in C^\infty_\Gamma(\tilde{\Sigma})$ and $\alpha:\tilde{\Sigma}\rightarrow\tilde{\Sigma}$ are such that $\|f\|_{L^\infty}<\epsilon$ and $\hat{f}=\hat{f}\circ\alpha$, and if $d(\alpha(p),p)<\epsilon$ for some point, $p\in\tilde{\Sigma}$, then $\alpha=\opId$.
\endproclaim
\proclabel{ConditionForIdentity}
\proof If $\alpha(p)\in B_\epsilon(p)$, then
$$\eqalign{
\alpha(p) &= (\pi\circ\Cal{E}_p^{-1}\circ\Cal{E})(\alpha(p),f(\alpha(p)))\cr
&= (\pi\circ\Cal{E}_p^{-1}\circ\hat{f}\circ\alpha)(p)\cr
&= (\pi\circ\Cal{E}_p^{-1}\circ\hat{f})(p)\cr
&= (\pi\circ\Cal{E}_p^{-1})(p,f(p))\cr
&= p.\cr}
$$
That is, if $d(\alpha(p),p)<\epsilon$, then $d(\alpha(p),p)=0$. In particular, if $U\subseteq\tilde{\Sigma}$ is the set of all points, $p$, of $\Sigma$ such that $d(\alpha(p),p)<\epsilon$, then $U$ is both open and closed, and the result now follows by connectedness.\qed
\medskip
\noindent By compactness, there exists $0<\delta<\epsilon$ such that for all $p\in\tilde{\Sigma}$,
$$
B_\delta(i(p)) \subseteq \Cal{E}(B_\epsilon(p)\times]-\epsilon,\epsilon[).
$$
In particular, $\Cal{E}_p^{-1}$ is well defined over $B_\delta(i(p))$.
\proclaim{Lemma \nextprocno}
\noindent If $f\in C^\infty_\Gamma(\tilde{\Sigma})$ and $\alpha:\tilde{\Sigma}\rightarrow\tilde{\Sigma}$ are such that $\hat{f}=\hat{f}\circ\alpha$ and if $\|f\|_{L^\infty}<\delta/4$, then, for all $p\in\tilde{\Sigma}$ and for all $q\in B_{\delta/4}(p)$,
$$
\alpha(q) = (\pi\circ\Cal{E}_{\alpha(p)}^{-1}\circ\hat{f})(q).
$$
\endproclaim
\proclabel{Closeness}
\proof For all $q\in B_{\delta/4}(p)$,
$$
d(\hat{f}(q),i(\alpha(p)) \leq d(\hat{f}(q),i(q)) + d(i(q),i(p)) + d(i(p),\hat{f}(p)) + d(\hat{f}(p),i(\alpha(p))) < \delta,
$$
so that $\Cal{E}_{\alpha(p)}^{-1}$ is well defined over $\hat{f}(B_{\delta/4}(p))$. Now let $U\subseteq B_{\delta/4}(p)$ be the set of all points such that $\alpha(q)=(\pi\circ\Cal{E}_{\alpha(p)}^{-1}\circ\hat{f})(q)$. Since this set contains $p$, it is non-empty. By continuity, it is relatively closed. To see that it is open, fix $q\in U$. For sufficiently small $\eta>0$, $\alpha(B_\eta(q))\subseteq B_\epsilon(\alpha(p))$, and so, for $r\in B_\eta(q))$,
$$\eqalign{
\alpha(r) &= (\pi\circ\Cal{E}_{\alpha(p)}^{-1}\circ\Cal{E})(\alpha(r),f(\alpha(r))\cr
&=(\pi\circ\Cal{E}_{\alpha(p)}^{-1}\circ\hat{f}\circ\alpha)(r)\cr
&=(\pi\circ\Cal{E}_{\alpha(p)}^{-1}\circ\hat{f})(r),\cr}
$$
as desired. We conclude that $U=B_{\delta/4}(p)$, and this completes the proof.\qed
\proclaim{Lemma \nextprocno}
\noindent Let $(f_m)\in C^\infty_\Gamma(\tilde{\Sigma})$ be a sequence converging to $0$ in the $C^0$ sense, and let $(\alpha_m):\tilde{\Sigma}\rightarrow\tilde{\Sigma}$ be a sequence of continuous maps. If $\alpha_m:\tilde{\Sigma}\rightarrow\tilde{\Sigma}$ such that $\hat{f}_m=\hat{f}_m\circ\alpha_m$ for all $m$, then the sequence $(\alpha_m)$ is equicontinuous.
\endproclaim
\proclabel{Equicontinuity}
\proof Indeed, we may suppose that $\|f_m\|_{L^\infty}<\delta/4$ for all $m$. By Lemma \procref{Closeness}, for all $p\in\tilde{\Sigma}$ and for all $q\in B_{\delta/4}(p)$,
$$
\alpha_m(q) = (\pi\circ\Cal{E}_{\alpha_m(p)}^{-1}\circ\hat{f}_m)(q).
$$
However, by cocompactness, the family $(\pi\circ\Cal{E}_r^{-1})_{r\in\tilde{\Sigma}}$ is equicontinuous. Furthermore, since $(f_m)$ converges to $0$ in the $C^0$ sense, the family $(\hat{f}_m)$ is also equicontinous. Since composition preserves equicontinuity, the result follows.\qed
\medskip
{\bf\noindent Proof of Theorem \procref{Stability}:}\ Suppose the contrary. There exists a sequence $(f_m)\in C^\infty_\Gamma(\tilde{\Sigma})$ converging to $0$ in the $C^0$ sense, and a sequence $(\alpha_m)\in\opDiff(\tilde{\Sigma})\setminus\Gamma$ such that $\hat{f}_m=\hat{f}_m\circ\alpha_m$ for all $m$. Fix an element, $p$, of $\tilde{\Sigma}$. By cocompactness, there exists a sequence, $(\gamma_m)\in\Gamma$, such that the sequence of points, $(\gamma_m\circ\alpha_m(p))$, is contained in some compact subset of $\tilde{\Sigma}$. It follows by Lemma \procref{Equicontinuity} and the Arzela-Ascoli Theorem, that, upon extracting the subsequence if necessary, $(\gamma_m\circ\alpha_m)$ converges to some continuous function, $\alpha_\infty$, say. However, taking limits yields $i\circ\alpha_\infty=i$, so that, by Lemma \procref{ConditionForDiffeomorphism}, $\alpha_\infty$ is a smooth diffeomorphism. Thus, since $i$ is prime, $\alpha_\infty\in\Gamma$, and we may suppose that $\alpha_\infty=\opId$. In particular, for sufficiently large $m$, $d(\gamma_m\circ\alpha_m(p),p)<\epsilon$ for all $p\in\tilde{\Sigma}$, so that, by Lemma \procref{ConditionForIdentity}, $\gamma_m\circ\alpha_m=\opId$. It follows that $\alpha_m\in\Gamma$ for all sufficiently large $m$, which is absurd. This completes the proof.\qed
\goodbreak
\newhead{Bibliography}[Bibliography]
{\leftskip = 5ex \parindent = -5ex
\leavevmode\hbox to 4ex{\hfil \cite{Aronszajn}}\hskip 1ex{Aronszajn N., A unique continuation theorem for solutions of elliptic partial differential equations or inequalities of second order, {\sl J. Math. Pures Appl.}, {\bf 36}, (1957), 235--249}
\medskip
\leavevmode\hbox to 4ex{\hfil \cite{Berger}}\hskip 1ex{Berger M., {\sl A panoramic view of Riemannian geometry}, Springer-Verlag, Berlin, (2003)}
\medskip
\leavevmode\hbox to 4ex{\hfil \cite{CaffNirSprV}}\hskip 1ex{Caffarelli L., Nirenberg L., Spruck J., Nonlinear second-order elliptic equations. V. The Dirichlet problem for Weingarten hypersurfaces, {\sl Comm. Pure Appl. Math.}, {\bf 41}, (1988), no. 1, 47--70}
\medskip
\leavevmode\hbox to 4ex{\hfil \cite{ColdingDeLellis}}\hskip 1ex{Colding T. H., De Lellis C., The min-max construction of minimal surfaces, {\sl Surv. Differ. Geom.}, {\bf VIII}, Int. Press, Somerville, MA, (2003), 75--107}
\medskip
\leavevmode\hbox to 4ex{\hfil \cite{CrandhallIshiiLyons}}\hskip 1ex{Crandall M. G., Ishii H., Lions P.-L., User's guide to viscosity solutions of second order partial differential equations, {\sl Bull. Amer. Math. Soc.}, {\bf 27}, (1992), no. 1, 1--67}
\medskip
\leavevmode\hbox to 4ex{\hfil \cite{ElworthyTromba}}\hskip 1ex{Elworthy K. D., Tromba A. J., Degree theory on Banach manifolds, in {\sl Nonlinear Functional Analysis (Proc. Sympos. Pure Math.)}, Vol. {\bf XVIII}, Part 1, Chicago, Ill., (1968), 86--94, Amer. Math. Soc., Providence, R.I.}
\medskip
\leavevmode\hbox to 4ex{\hfil \cite{EspRos}}\hskip 1ex{Espinar J. M., Rosenberg H., When strictly locally convex hypersurfaces are embedded, {\sl Math. Zeit.}, {\bf 27}, (2012), nos. 3--4, 1075--1090}
\medskip
\leavevmode\hbox to 4ex{\hfil \cite{GilbTrud}}\hskip 1ex{Gilbarg D., Trudinger N. S., {\sl Elliptic partical differential equations of second order}, Die Grundlehren der mathemathischen Wissenschaften, {\bf 224}, Springer-Verlag, Berlin, New York (1977)}
\medskip
\leavevmode\hbox to 4ex{\hfil \cite{GuillemanPollack}}\hskip 1ex{Guillemin V., Pollack A., {\sl Differential Topology}, Prentice-Hall, Englewood Cliffs, N.J., (1974)}
\medskip
\leavevmode\hbox to 4ex{\hfil \cite{Hamilton}}\hskip 1ex{Hamilton R. S., The inverse function theorem of Nash and Moser, {\sl Bull. Amer. Math. Soc.}, {\bf 7}, (1982), no. 1, 65--222}
\medskip
\leavevmode\hbox to 4ex{\hfil \cite{Huisken}}\hskip 1ex{Huisken G., Contracting convex hypersurfaces in Riemannian manifolds by their mean curvature, {\sl Invent. Math.}, {\bf 84}, (1986), no. 3, 463--480}
\medskip
\leavevmode\hbox to 4ex{\hfil \cite{Kato}}\hskip 1ex{Kato T., {\sl Perturbation theory for linear operators}, Grundlehren der Mathematischen Wissenschaften, {\bf 132}, Springer-Verlag, Berlin, New York, (1976)}
\medskip
\leavevmode\hbox to 4ex{\hfil \cite{LabA}}\hskip 1ex{Labourie F., Probl\`emes de Monge-Amp\`ere, courbes holomorphes et laminations, {\sl Geom. Funct. Anal.}, {\bf 7}, (1997), no. 3, 496--534}
\medskip
\leavevmode\hbox to 4ex{\hfil \cite{LabB}}\hskip 1ex{Labourie F., Immersions isom\'etriques elliptiques et courbes pseudo-holomorphes, {\sl Geometry and topology of submanifolds} (Marseille, 1987), 131--140, World Sci. Publ., Teaneck, NJ, (1989)}
\medskip
\leavevmode\hbox to 4ex{\hfil \cite{MaximoNunesSmith}}\hskip 1ex{M\'aximo D., Nunes I. P., Smith G., Free boundary minimal annuli in convex three-manifolds, to appear in {\sl J. Diff. Geom.}}
\medskip
\leavevmode\hbox to 4ex{\hfil \cite{Milnor}}\hskip 1ex{Milnor J. W., {\sl Topology from the differential viewpoint}, Princeton Landmarks in Mathematics, Princeton, (1997)}
\medskip
\leavevmode\hbox to 4ex{\hfil \cite{PacardXu}}\hskip 1ex{Pacard F., Xu X., Constant mean curvature spheres in Riemannian manifolds, {\sl Manu\-scripta Math.}, {\bf 128}, (2009), no. 3, 275--295}
\medskip
\leavevmode\hbox to 4ex{\hfil \cite{Robeday}}\hskip 1ex{Robeday A., Masters Thesis, Univ. Paris VII}
\medskip
\leavevmode\hbox to 4ex{\hfil \cite{RosSchneid}}\hskip 1ex{Rosenberg H., Schneider M., Embedded constant curvature curves on convex surfaces, {\sl Pac. J. Math.}, {\bf 253}, (2011), no. 1, 213--219}
\medskip
\leavevmode\hbox to 4ex{\hfil \cite{SchneiderI}}\hskip 1ex{Schneider M., Closed magnetic geodesics on $S^2$, {\sl J. Differential Geom.}, {\bf 87}, (2011), no. 2, 343--388}
\medskip
\leavevmode\hbox to 4ex{\hfil \cite{SchneiderII}}\hskip 1ex{Schneider M., Closed magnetic geodesics on closed hyperbolic Riemann surfaces, {\sl J. London Math. Soc.}, {\bf 105}, (2012), 424--446}
\medskip
\leavevmode\hbox to 4ex{\hfil \cite{SimonSmith}}\hskip 1ex{Smith F. R., On the existence of embedded minimal $2$-spheres in the $3$-sphere, endowed with an arbitrary metric, {\sl Bull. Austral. Math. Soc.}, {\bf 28}, (1983), 159--160}
\medskip
\leavevmode\hbox to 4ex{\hfil \cite{SmiAAT}}\hskip 1ex{Smith G., An Arzela-Ascoli Theorem for Immersed Submanifolds, {\sl Ann. Fac. Sci. Toulouse Math.}, {\bf 16}, no. 4, (2007), 817--866}
\medskip
\leavevmode\hbox to 4ex{\hfil \cite{SmiEC}}\hskip 1ex{Smith G., Constant curvature hyperspheres and the Euler Characteristic,\break arXiv:1103.3235}
\medskip
\leavevmode\hbox to 4ex{\hfil \cite{SmiPPG}}\hskip 1ex{Smith G., The Plateau Problem for General Curvature Functions, arXiv:1008.3545}
\medskip
\leavevmode\hbox to 4ex{\hfil \cite{SmiSLC}}\hskip 1ex{Smith G., Special Lagrangian Curvature, {\sl Math. Annalen}, {\bf 335}, (2013), no. 1, 57--95}
\medskip
\leavevmode\hbox to 4ex{\hfil \cite{SmiAC}}\hskip 1ex{Smith G., Bifurcation of solutions to the Allen-Cahn equation, to appear in {\sl J. London Math. Soc.}}
\medskip
\leavevmode\hbox to 4ex{\hfil \cite{Spivak}}\hskip 1ex{Spivak M., {\sl A comprehensive introduction to differential geometry. Vol. III.}, Publish or Perish, Inc., Wilmington, Del., second edition, (1979)}
\medskip
\leavevmode\hbox to 4ex{\hfil \cite{Tromba}}\hskip 1ex{Tromba A. J., The Euler characteristic of vector fields on Banach manifolds and a globalization of Leray-Schauder degree, {\sl Adv. in Math.}, {\bf 28}, (1978), no. 2, 148--173}
\medskip
\leavevmode\hbox to 4ex{\hfil \cite{WhiteI}}\hskip 1ex{White B., The space of m-dimensional surfaces that are stationary for a parametric elliptic functional, {\sl Indiana Univ. Math. J.}, {\bf 36}, (1987), no. 3, 567--602}
\medskip
\leavevmode\hbox to 4ex{\hfil \cite{WhiteII}}\hskip 1ex{White B., Every three-sphere of positive Ricci curvature contains a minimal embedded torus, {\sl Bull. Amer. Math. Soc.}, {\bf 21}, (1989), no. 1, 71--75}
\medskip
\leavevmode\hbox to 4ex{\hfil \cite{WhiteIII}}\hskip 1ex{White B., Existence of smooth embedded surfaces of prescribed genus that minimize parametric even elliptic functionals on 3-manifolds, {\sl J. Differential Geom.}, {\bf 33}, (1991), no. 2, 413--443}
\medskip
\leavevmode\hbox to 4ex{\hfil \cite{WhiteIV}}\hskip 1ex{White B., The space of minimal submanifolds for varying Riemannian metrics, {\sl Indiana Univ. Math. J.}, {\bf 40}, (1991), no. 1, 161--200}
\medskip
\leavevmode\hbox to 4ex{\hfil \cite{Ye}}\hskip 1ex{Ye R., Foliation by constant mean curvature spheres, {\sl Pacific J. Math.}, {\bf 147}, (1991), no. 2, 381--396}
\par}
%
%
%
%
\enddocument